\documentclass[a4paper,11pt,reqno]{amsart} 

\usepackage{aamath}
\usepackage{mhequ}
\usepackage[final,notref,notcite]{showkeys}
\usepackage{longtable}

\newcommand{\abbr}[1]{{\small\sc\lowercase{#1}}}
\newcommand{\ie}{{\it i.e.}, }
\newcommand{\eg}{\emph{e.g.},~}
\newcommand{\cal}[1]{{\mathcal{#1}}}
\newcommand\croc[1]{\left\langle #1\right\rangle}

\usepackage{caption}
\usepackage{subcaption}
\usepackage[hidelinks]{hyperref}
\usepackage[foot]{amsaddr}

\usepackage{comment}
\usepackage{fullpage}
\usepackage{microtype}

\usepackage[hidelinks]{hyperref}

\usepackage{bbm}

\def\Kz{K_0}

\usepackage{enumerate}
 \usepackage{tikz}
 \usetikzlibrary{shapes,automata,positioning,arrows,arrows.meta,fit,decorations.pathreplacing}
\usepackage{pgfplots}
\pgfplotsset{compat=1.7}

\usepackage{todonotes}

\usetikzlibrary{shapes,arrows}

\def\G{\mathcal G}
\def\L{\mathcal L}
\def\Z{\mathbb{Z}}
\def\S{\mathcal{S}}
\def\C{\mathcal{C}}
\def\Re{\mathcal{R}}
\def\R{\mathbb{R}}
\def\E{\mathbb{E}}
\def\G{\mathcal{G}}

\def\sg{B_\gamma(K)}
\def\sz{B_0(K)}

\def\ss{B_{\gamma,n}}

\def\epsilon{\varepsilon}

\def\Z{\mathbb Z}

\def\NN{\mathbb Z_{\geq 0}^2}

\def\OO{{\mathbf X}}
\def\XX{{\mathbf Z}^n}

\def\r{\beta}

\title{A proof of the Chemical Recurrence Conjecture in two dimensions}

\author[1]{Andrea Agazzi$^{1,*}$}
\address[1]{Department of Mathematics and Statistics, University of Bern, \href{mailto:andrea.agazzi@unibe.ch}{andrea.agazzi@unibe.ch}}
\author[2]{David F. Anderson$^{2}$}
\address[2]{Department of Mathematics, University of Wisconsin-Madison, \href{mailto:anderson@math.wisc.edu}{anderson@math.wisc.edu},}
\author[3]{Daniele Cappelletti$^{3}$}
\address[3]{Department of Mathematical Sciences, Politecnico di Torino, \href{mailto:daniele.cappelletti@polito.it}{daniele.cappelletti@polito.it}}
\author[4]{Lucie Laurence$^{4}$}
\address[4]{Department of Mathematics and Statistics, University of Bern, \href{mailto:lucie.laurence@unibe.ch}{lucie.laurence@unibe.ch}}
\author[5]{Jonathan~C.~Mattingly$^{5}$}
\address[5]{Department of Mathematics and Department of Statistics, Duke University, \href{mailto:jonm@math.duke.edu}{jonm@math.duke.edu}}
\address[*]{Corresponding author}

\date{
    \today
}

\begin{document}
\begin{abstract}
   At the microscopic level, the dynamics of a network of chemical reactions can be modeled, in the well-mixed case, as a continuous-time Markov chain. This process counts the number of particles in each of the $d$ species as a function of time, and its state-dependent, polynomial jump rates are determined by the structure of the underlying reaction graph. Networks of chemical reactions whose reaction graphs are the union of strongly connected components are called \emph{weakly reversible}. The long-standing \emph{chemical recurrence conjecture} states that if a network is weakly reversible, then the associated Markov chain is positive recurrent. We prove this conjecture for $d=2$.
\end{abstract}
\maketitle
\tableofcontents

\section{Introduction}

Understanding the dynamical properties of large systems of chemical reactions - modeling, e.g., the living cell - is a problem of paramount importance in systems biology.  Arguably the most fundamental and widely applied set of laws connecting a network of chemical reactions to its dynamics is mass action kinetics. 
This family of models associates to each network a dynamical system, which is either deterministic or stochastic, describing the time evolution of the systems' state (respectively the concentration or the number of molecules of each species).

This results in a strong relationship between the structure of a reaction network and its dynamics. The discipline of Chemical Reaction Network (CRN) theory, in which this paper is inscribed, focuses on the study of this relationship.

In the deterministic setting, usually utilized when the number of molecules is large, the 
dynamics of the chemicals' concentrations are typically modeled via a set of ordinary differential equations (ODEs).
Here,
the study of the relationship between the structure and the long-time dynamical properties of CRNs was pioneered by the works \cite{horn72, feinberg72}. In these papers, the authors introduce the concept of \emph{complex balanced} equilibria, generalizing the notion of detailed balance.

Furthermore, these works prove that such complex balanced fixed points must be unique and locally asymptotically stable. 
This introduced the long-standing Global Attractor Conjecture, which states that each complex-balanced equilibrium is, in fact, \textit{globally} asymptotically stable.  The work \cite{anderson11} proved the conjecture in a special case. This work was then extended in \cite{gopal13}, which introduced  the class of \emph{strongly endotactic networks}.

In the stochastic setting, the relationship between complex balancing and invariant measures was explored in \cite{anderson10,cappelletti16,cappelletti18}, where it was shown that the stationary distributions of complex balanced CRNs have a Poisson product structure. Further, \cite{cappelletti23} studies the convergence rate to these invariant measures for a subclass of such models. 

However, the complex balanced condition is quite restrictive, speaking not only to the structure of the network but also requiring very specific choices of model parameters, and is therefore rarely verified in practice. This calls for results requiring  less stringent conditions for existence of stationary distributions (or positive recurrence) for this class of stochastic models. In \cite{agazzi20, anderson2020tier, anderson18} it was shown that conditions for deterministic stability are not sufficient for the existence of a stationary distribution in the stochastic setting. Indeed, there the authors show that strongly endotactic networks, proven to be asymptotically stable in the deterministic setting \cite{gopal13}, can have transient or explosive dynamics when modeled stochastically. Hence, the question of identifying general structural conditions guaranteeing existence of stationary distributions for stochastically modeled CRNs remains widely open.

The following conjecture, well-known in the CRN community, aims to address the question raised above by identifying \emph{weak reversibility} (\dref{d:wr} below) as a sufficient structural condition for positive recurrence. 

\begin{conjecture}	\label{c:main}
	If a \abbr{chemical reaction network} is weakly reversible, then the associated continuous-time Markov chain is positive recurrent.
\end{conjecture}
This paper proves this conjecture for chemical reaction networks with $2$ species, \ie in dimension $d=2$.

\subsection{Notation}\label{sec:model} We denote by $\R$ and $\Z$ the sets of real and integer numbers, respectively. For a given real number $a$, we denote by $\R_{\geq a}$ and $\Z_{\geq a}$ the real and  integer numbers greater than or equal to $a$, respectively.

A \emph{reaction network} is a triple $(\S,\C,\Re)$ where  $\S$, $\C$, and $\Re$ are defined as follows.
\begin{itemize}
    \item $\mathcal S$ is a finite set of \emph{species}, that is a set of $d$ distinct symbols. It is often convenient to denote the species as $\{S_1,\dots, S_d\}$, as this allows us to refer to species via their index. However, because we focus on  two-dimensional  models in this paper, we will often denote our species set as $\{A,B\}$.

\item $\C$ is a finite set of \emph{complexes}.
 We assume each complex is a linear combination of the species over $\Z_{\geq0}$. Here, each coefficient counts the multiplicity of the corresponding species in the complex. Complexes will be regarded as vectors in $\Z_{\geq0}^d$, given that an ordering for the species is chosen.
 \item  $\Re$ is a finite set of \emph{reactions} of the form $y \to y'$ for $y,y' \in \C$. 
 We therefore identify $\Re$ with a finite subset of $\C\times\C$ with the property that for any $y\in\C$ we have $(y,y)\notin\Re$.  
 \end{itemize}
 
 We will enumerate the reactions and denote the ``$r$th'' reaction via $y_r \to y_r'$.  We will also slightly abuse notation and write sums over the reactions as $\sum_{r \in \Re}$.  
 When enumeration is not needed, we will simply denote a generic reaction via $y\to y'$.  Given a reaction $y\to y'\in \Re$ we say that $y$ is the \emph{source complex} and $y'$ is the \emph{product complex} of the reaction $y\to y'$.

Given a reaction network $(\S,\C,\Re)$, one can construct a directed graph with nodes $\C$ and edges $\Re$. This directed graph is called the \emph{reaction graph}. In this paper we assume that all complexes appear in at least one reaction and all species appear in at least one complex. Under this assumption, the reaction graph uniquely determines a reaction network. In fact, reaction networks are usually described by means of their reaction graph, and the same will be done in the present paper.

\begin{example}[Main example] Let $\mathcal S := \{A, B\}$ and consider the reactions
  \begin{equ} r_1~:~0 \to A+B\,, \qquad r_2~:~A+B \to 3A+2B\,, \qquad r_3~:~3A+2B \to 0\,.\end{equ}
 The complexes in this case can be identified with the set $\C := \{(0,0), (1,1), (3,2)\}$ counting the multiplicities of the species in the input/outputs of the reactions above ($\{0A+0B, 1A+1B,3A+2B\}$). The set of reactions $\mathcal R = \{r_1,r_2,r_3\}$ can be identified with the subset of $\C\times \C$ given by $\{((0,0),(1,1)), ((1,1),(3,2)), ((3,2),(0,0))\}$. The reaction graph associated to $(\mathcal S, \C, \Re)$ is obtained by interpreting $\C$ as vertices and $\Re$ as edges:
\begin{align}\label{ex:musterbeispiel}
\begin{split}
    0 &\longrightarrow A+B.\\
    &\nwarrow \hspace{.2in} \swarrow\\
    & 3A + 2B
    \end{split}
\end{align}
   \hfill $\triangle$
\end{example}
\noindent It is easy to see that the reaction graph \eref{ex:musterbeispiel} is \emph{weakly reversible}, as codified in the next definition.

\begin{definition}[\cite{horn72cb, feinberg72}]\label{d:wr} If all the connected components of a network's reaction graph are strongly connected, then the network is called \emph{weakly reversible}.\end{definition}
Weak reversibility is among the oldest structural conditions considered in the mathematical study of chemical reaction networks and is intimately connected with complex balancing:

  \begin{remark}\label{r:cb=wr}
    It was shown in \cite[Theorem 2B]{horn72cb} that every chemical reaction network which permits a complex balanced mass action system is weakly reversible. Namely, the class of weakly reversible CRNs contains the one of complex balanced networks.
  \end{remark}

The \emph{Euclidean embedded graph} \cite{craciun2019polynomial} is a related way to visualize chemical reaction network models.  As mentioned above, given an ordering of the species, we can associate each complex $y$ with a point in $\mathbb Z_{\geq 0}^d$ by interpreting its multiplicities as coordinates. Then, each reaction is identified with a directed edge of the embedded graph, represented by an arc connecting two such points (the input, respectively output, complex of that reaction). Note that this notion of graph embedding is weaker than the one customarily used in graph theory, as edges of this embedded graph \emph{can} intersect with each other.

  \setcounter{example}{0}
\begin{example}[continued] The Euclidean embedded graph associated to \eref{ex:musterbeispiel} is represented in \fref{fig:EuclideanGraphExample1} (A).
  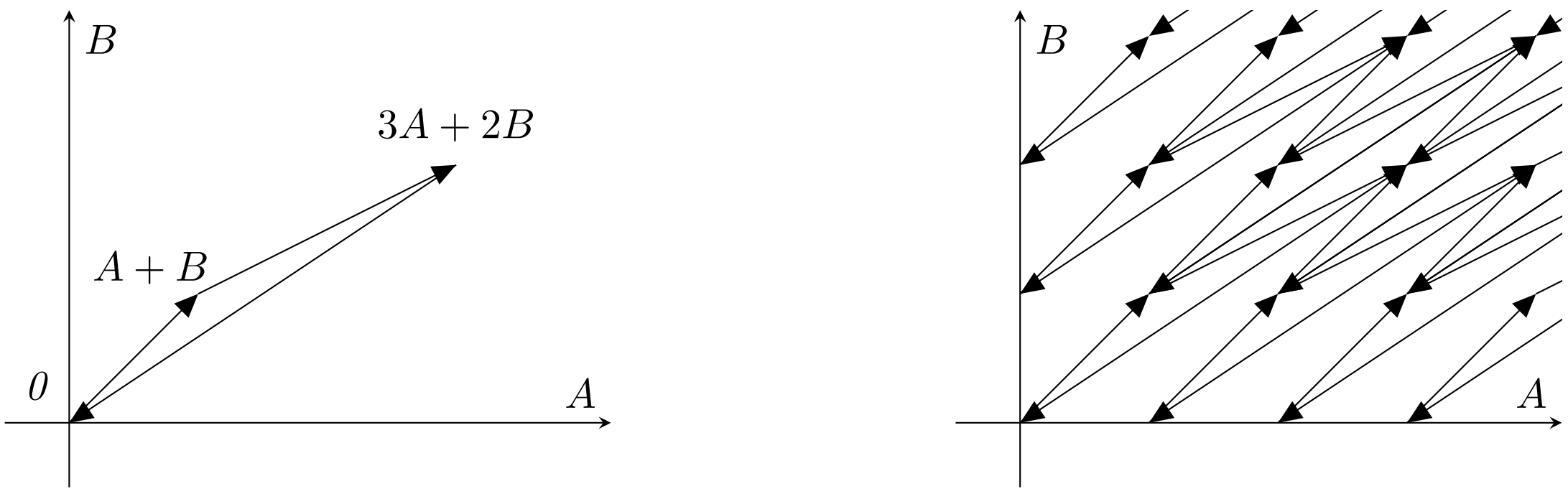
\begin{figure}
    \centering
\begin{subfigure}[b]{0.45\textwidth}
    \centering
    \begin{tikzpicture}
     \begin{axis}[
     axis y line=middle,axis x line=middle,
       xmin=-0.5, xmax=4.2, ymin=-0.5, ymax=3.2,
       unit vector ratio*=1 1 1,
       ticks=none,
       xlabel={$A$}, ylabel={$B$},
       width=1\textwidth]

      \addplot[mark=none, black, -triangle 45]
               coordinates{
                           (0,0)
                           (1,1)
                           };
      \addplot[mark=none, black, -triangle 45]
               coordinates{
                           (1,1)
                           (3,2)
                           };
      \addplot[mark=none, black, -triangle 45]
               coordinates{
                           (3,2)
                           (0,0)
                           };
      \node[label={135:{0}},inner sep=2pt] at (axis cs:0,0) {};
      \node[label={0:{$A+B$}},inner sep=2pt] at (axis cs:0,1.2) {};
      \node[label={90:{$3A+2B$}},inner sep=2pt] at (axis cs:3,2) {};
     \end{axis}
    \end{tikzpicture}
    \caption{Euclidean embedded graph of \eref{ex:musterbeispiel}}
    \label{f:euclideang}
\end{subfigure}
\hfill
\begin{subfigure}[b]{0.45\textwidth}
    \centering
    \begin{tikzpicture}
     \begin{axis}[
     axis y line=middle,axis x line=middle,
       xmin=-0.5, xmax=4.2, ymin=-0.5, ymax=3.2,
       unit vector ratio*=1 1 1,
       ticks=none,
       xlabel={$A$}, ylabel={$B$},
       width=1\textwidth]

      \addplot[mark=none, black, -triangle 45]
               coordinates{
                           (0,0)
                           (1,1)
                           };
      \addplot[mark=none, black, -triangle 45]
               coordinates{
                           (1,1)
                           (3,2)
                           };
      \addplot[mark=none, black, -triangle 45]
               coordinates{
                           (3,2)
                           (0,0)
                           };
           \addplot[mark=none, black, -triangle 45]
               coordinates{
                           (4,3)
                           (1,1)
                           };
     \addplot[mark=none, black, -triangle 45]
               coordinates{
                           (5,3)
                           (2,1)
                           };
     \addplot[mark=none, black, -triangle 45]
               coordinates{
                           (6,3)
                           (3,1)
                           };
     \addplot[mark=none, black, -triangle 45]
               coordinates{
                           (5,4)
                           (2,2)
                           };
     \addplot[mark=none, black, -triangle 45]
               coordinates{
                           (6,4)
                           (3,2)
                           };
     \addplot[mark=none, black, -triangle 45]
               coordinates{
                           (4,4)
                           (1,2)
                           };
     \addplot[mark=none, black, -triangle 45]
               coordinates{
                           (4,5)
                           (1,3)
                           };
     \addplot[mark=none, black, -triangle 45]
               coordinates{
                           (5,5)
                           (2,3)
                           };
     \addplot[mark=none, black, -triangle 45]
               coordinates{
                           (6,5)
                           (3,3)
                           };
     \addplot[mark=none, black, -triangle 45]
               coordinates{
                           (1,0)
                           (2,1)
                           };
     \addplot[mark=none, black, -triangle 45]
               coordinates{
                           (2,0)
                           (3,1)
                           };
     \addplot[mark=none, black, -triangle 45]
               coordinates{
                           (3,0)
                           (4,1)
                           };
     \addplot[mark=none, black, -triangle 45]
               coordinates{
                           (0,1)
                           (1,2)
                           };
     \addplot[mark=none, black, -triangle 45]
               coordinates{
                           (1,1)
                           (2,2)
                           };
     \addplot[mark=none, black, -triangle 45]
               coordinates{
                           (2,1)
                           (3,2)
                           };
     \addplot[mark=none, black, -triangle 45]
               coordinates{
                           (3,1)
                           (4,2)
                           };
     \addplot[mark=none, black, -triangle 45]
               coordinates{
                           (0,2)
                           (1,3)
                           };
     \addplot[mark=none, black, -triangle 45]
               coordinates{
                           (1,2)
                           (2,3)
                           };
     \addplot[mark=none, black, -triangle 45]
               coordinates{
                           (2,2)
                           (3,3)
                           };
     \addplot[mark=none, black, -triangle 45]
               coordinates{
                           (3,2)
                           (4,3)
                           };
      \addplot[mark=none, black, -triangle 45]
               coordinates{
                           (4,2)
                           (1,0)
                           };
      \addplot[mark=none, black, -triangle 45]
               coordinates{
                           (5,2)
                           (2,0)
                           };
      \addplot[mark=none, black, -triangle 45]
               coordinates{
                           (3,3)
                           (0,1)
                           };
      \addplot[mark=none, black, -triangle 45]
               coordinates{
                           (4,3)
                           (1,1)
                           };
      \addplot[mark=none, black, -triangle 45]
               coordinates{
                           (5,3)
                           (2,1)
                           };
      \addplot[mark=none, black, -triangle 45]
               coordinates{
                           (6,3)
                           (3,1)
                           };
     \addplot[mark=none, black, -triangle 45]
              coordinates{
                          (2,1)
                          (4,2)
                          };
      \addplot[mark=none, black, -triangle 45]
               coordinates{
                           (3,1)
                           (5,2)
                           };

     \addplot[mark=none, black, -triangle 45]
              coordinates{
                          (1,2)
                          (3,3)
                          };
      \addplot[mark=none, black, -triangle 45]
               coordinates{
                           (2,2)
                           (4,3)
                           };
     \addplot[mark=none, black, -triangle 45]
              coordinates{
                          (4,1)
                          (6,2)
                          };
      \addplot[mark=none, black, -triangle 45]
               coordinates{
                           (3,2)
                           (5,3)
                           };
       \addplot[mark=none, black, -triangle 45]
                coordinates{
                            (4,2)
                            (6,3)
                            };
        \addplot[mark=none, black, -triangle 45]
                 coordinates{
                             (6,2)
                             (3,0)
                             };
       \addplot[mark=none, black, -triangle 45]
                coordinates{
                            (3,4)
                            (0,2)
                            };
        \addplot[mark=none, black, -triangle 45]
                 coordinates{
                             (7,5)
                             (4,3)
                             };
     \end{axis}
    \end{tikzpicture}
    \caption{Transition graph of \eref{ex:musterbeispiel}}
    \label{f:transitiong}
\end{subfigure}
\caption{The Euclidean embedded graph  (A) 
and transition graph  (B)
for the reaction graph in \eqref{ex:musterbeispiel}.}   \label{fig:EuclideanGraphExample1}
\end{figure}
 \hfill $\triangle$
\end{example}

We will now define the key concept of ``endotactic networks'' (etymologically ``organized inwards'') which aims to capture networks whose dominant reactions asymptotically ``pull'' its dynamics, introduced in the next section, towards a compact set in state space.

In this spirit, 
denoting throughout by $S^{d-1}$ the Euclidean sphere in $\mathbb R^d$, we define the set of $w$-dominant, respectively $w$-nontrivial dominant, reactions as
\begin{align}
\label{e:dominant}
\begin{split}
  \mathcal R_w &:= \arg\max\{\dtp{y,w}~:~y \to y' \in \mathcal R\}\\
  \bar{\mathcal R}_w &:= \arg\max\{\dtp{y,w}~:~y \to y' \in \mathcal R, \dtp{y'-y,w} \neq 0\}.
  \end{split}
\end{align}
We now introduce the classes of \emph{endotactic} and \emph{strongly endotactic} networks:
\begin{definition}
  For any $w \in S^{d-1}$, we call a CRN $(\mathcal S, \C, \Re)$ $w$-\emph{endotactic} if  we have that $\dtp{w, \zeta_r} < 0$ for all $r \in \bar{\mathcal R}_w$. We further call such a CRN $w$-\emph{weakly endotactic} ($w$-WE) if $\bar{\mathcal R}_w \cap {\mathcal R}_w = \emptyset$ and $w$-\emph{strongly endotactic} ($w$-SE) otherwise. Finally, a CRN is \emph{endotactic} if it is $w$-endotactic for all $w \in S^{d-1}$ and \emph{strongly endotactic} if it is $w$-strongly endotactic for all $w \in S^{d-1}$.
\end{definition}
Strong endotacticity has been proven to be a sufficient condition for stability in the deterministic setting \cite{anderson2020tier, gopal13} and has been identified as a key condition for establishing large deviations estimates in the stochastic setting \cite{ agazzi182,agazzi18,agazzi22}.

\begin{remark}\label{r:openc}
  Note that if a network is $w$-strongly endotactic, then, for sufficiently small $\epsilon > 0$, it is $w'$-strongly endotactic for all $w' \in S^{d-1}$ with $\|w'-w\| < \epsilon$.
\end{remark}
\begin{remark}\label{r:endotactic} It has been proven in \cite{gopal13} that weakly reversible CRNs are endotactic, but not necessarily strongly endotactic. \end{remark}

\subsection{Stochastic model}\label{sec:stochastic_model}

Given a reaction network $(\S,\C,\Re)$, a \emph{stochastic kinetics} for $(\S,\C,\Re)$ is an assignment of a \emph{rate function} $\lambda_{y\to y'}:\Z_{\geq0}^d\to \R_{\geq0}$ to each reaction $y\to y'\in \Re$. 
 When an enumeration of the reactions is being utilized, we often write $\lambda_r$ instead of $\lambda_{y_r \to y_r'}$.  We assume throughout that $\lambda_{y\to y'}(x) = 0$ if $x_i<y_i$ for any $i \in \{1,\dots,d\}$.

For $y_r \to y_r'\in \Re$, we denote its \emph{reaction vector} via $\zeta_r = y_r'-y_r\in \Z^d$, which is the vector that determines how the counts of the species change due to a single occurrence of the $r$th reaction. Note that multiple reactions may have the same reaction vector. In this regards, for all $v\in \Z^d$ we define by
\begin{equation*}
 \Re(v):=\{y\to y'\in\Re\,:\,y'-y=v\}
\end{equation*}
the set of reaction with a common reaction vector $v$.

The change in chemical species counts is modeled by means of a continuous-time Markov chain $X$ with state space $\Z_{\geq0}^d$ and transition rates

\begin{align*}
    q(x,x+v)=\sum_{r\in\Re(v)}\lambda_{r}(x), \quad \text{ for $v \in \Z^d$.}
\end{align*}
This process is naturally connected to the associated embedded Markov chain $Y$, a discrete-time Markov chain with transition probabilities
\begin{equ}\label{e:transprob}
    p(x,x+v)=\frac{\sum_{r\in\Re(v)}\lambda_{r}(x)}{\sum_{\substack{r\in\Re}}\lambda_{r}(x)}, \quad \text{ for $v \in \Z^d$.}
\end{equ}
In case of an explosion occurring at a finite-time $T_\infty$, we consider $X(t)=\Delta$ for any $t\geq T_\infty$, where $\Delta$ is a cemetery state not contained in $\Z_{\geq0}^d$.  The infinitesimal generator for the model is the operator $\L$ acting on functions $f:\Z^d_{\ge 0} \to \R$ via
\begin{align}\label{eq:generator}
    \L f(x) = \sum_{r\in \Re} \lambda_{r}(x) ( f(x+\zeta_r)-f(x)).
\end{align}

A popular choice of stochastic kinetics is given by  \emph{stochastic mass-action kinetics}, where for any reaction $y\to y'\in\Re$
\begin{align}\label{e:rates}
\lambda_{y\to y'}(x)=\kappa_{y\to y'} \prod_{i=1}^d \frac{x_i!}{(x_i-y_i)!} \mathbbm{1}_{\{x_i\geq y_i\}},
\end{align}
for some \textit{rate constant} $\kappa_{y\to y'}>0$ and where $\mathbbm{1}_{A}$ is the indicator function on the set $A$. When an enumeration of the reactions is being utilized, we often write $\kappa_r$ instead of $\kappa_{y_r \to y_r'}$. The analysis of the present paper focuses on this choice of kinetics, where $\{\kappa_r\}$ is our arbitrary choice of rate constants. We define
\begin{equation*}
\kappa_+=\max_{r\in\Re}\kappa_r\quad\text{and}\quad\kappa_-=\min_{r\in\Re}\kappa_r.
\end{equation*}

We will consider the portion of an intensity function that is independent of the rate constant.  Hence, for $y \in \C$ we define
\begin{align}\label{e:rates2}
    \lambda_y(x) = \prod_{i = 1}^d \frac{x_i!}{(x_i-y_i)!} \mathbbm{1}_{\{x_i \ge y_i\}}.
\end{align}
\begin{remark}\label{r:active}Note that for any reaction $y\to y'$ we have $\lambda_{y\to y'}(x) = \kappa_{y\to y'} \lambda_y(x)$. Furthermore, note that in the case of mass-action kinetics, a reaction $y\to y'\in \Re$ has a positive probability of occurring at $x$ if and only if
$
    x_i \ge y_i
$
for all $i \in \{1,\dots, d\}$. \end{remark}

The dynamics described above naturally induces a transition graph on $\mathbb Z_{\geq 0}^d$, i.e.\ a graph representing all possible one-step transitions of the model, weighted by their probability of occurrence.
\begin{definition}The \emph{transition graph} $\mathcal G$ associated to a CRN $(\mathcal S, \C, \Re)$ with mass action kinetics $\{\kappa_r\}$ is a weighted directed graph with vertex set $V(\mathcal G) = \mathbb Z_{\geq 0}^d$ and edge set
  \begin{equ}
    E(\mathcal G) := \{x \to x'~:~\exists y \to y',\,x'-x = y'-y, \lambda_y(x) > 0\}.
  \end{equ}
The weight on an edge $x \to x'$ is given by \eref{e:transprob}.
\end{definition}

  Given a directed graph (e.g.\ a transition graph, or a subset of it) $\G'$ and two nodes $v,v'$ we write $v\Rightarrow_{\G'}v'$ if there exists a finite directed path in $\G'$ from $v$ to $v'$.

  \setcounter{example}{0}
  \begin{example}[continued] The transition graph associated to \eref{ex:musterbeispiel} is represented in \fref{fig:EuclideanGraphExample1} (B). \hfill $\triangle$
  \end{example}

\subsection{The main result}\label{s:main}

The main result of this paper is a proof of Conjecture~\ref{c:main} in dimension 2:
\begin{theorem}\label{t:main}
  Let $X$ be the continuous-time Markov chain associated to a weakly reversible {CRN} with two species $(d = 2)$ equipped with stochastic mass action kinetics. Then, $X$ is positive recurrent.
\end{theorem}

Fix a closed irreducible component $S$. It suffices to consider $|S|=\infty$.
Our proof approach combines an analysis of the dynamics near the boundary, where coarse-graining techniques are utilized. In the bulk of the state space, Lyapunov estimates are obtained using standard tools of stochastic calculus.
The difficulty of establishing this result lies on one hand on the delicate character of the averaging estimates it requires and on the other on the generality of the class of reaction networks it aims to apply to. Indeed, establishing the necessary dissipation estimates for any \emph{particular} given reaction network in the class of interest can often be done adopting an \emph{ad-hoc} procedure.  However, in order to treat \emph{all} reaction networks of interest it is required to develop a general procedure/algorithm allowing to identify, for each such network, the first order at which the process displays dissipative behavior and to leverage such understanding to show stability of such process.

The need for such detailed analysis is especially relevant when the dynamics of the continuous-time Markov chain are close to the axes of the state space. Here, the process develops a boundary layer and its behavior is qualitatively different from the one in the \emph{bulk}, i.e.,~the set of states away from the axes. 
While designing a Lyapunov function to establish dissipation in the bulk of the state space can be  achieved in a general way by adapting existing techniques,
when close to the boundary the averaged dynamics of the process need to be investigated in much higher detail.
The main difficulty in this sense is that the order of averaging required to observe dissipation depends on the network being considered and is not uniformly bounded from above on the class of networks of interest. For a given network we therefore sequentially increase the order of the averaging and, at each step, we identify the basins of attraction of the  dynamics. Then, we increase the order of the averaging -- which corresponds to analyzing the dynamics of the process on longer timescales -- and in doing so we merge some of these basins of attraction, analogously to what is done in Wentzell-Freidlin theory. 
This way, we eventually reach an order of averaging where the coarse-grained dynamics are sufficiently simple. Then, we show by contradiction that these simplified dynamics must dissipate. This second part of the proof, the contradiction, relies on the realization that the asymptotic behavior distilled in the previous analysis is independent of the rate constants of the model. Combining this with the existence of a choice of rate constants for which the invariant measure of the model exists and has a particularly simple form concludes the proof.

\subsection{A detailed outline of the proof}
\label{sec:strategy}

In this section, we provide the requisite theory to carry out the strategy outlined in the above paragraph.
Our main result pivots on the following result from \cite[Thm. 8.13]{Robert}, in turn following from \cite{filonov89}, allowing us to capture the dissipation of a process through its \emph{averaged} dynamics.
\begin{lemma}\label{l:filonov}
Let $X$ be an irreducible Markov jump process on $S$ and let $V:S\to\mathbb R_{\geq0}$ have compact sublevel sets. Suppose that every $x\in S$ admits a stopping time $\sigma_x>0$ such that
\[
 \E_x\sigma_x+\E_xV(X(\sigma_x))<\infty,
\]
and that, for some finite $F\subset S$ and $\beta>0$,
\begin{equation}\label{e:filonov-general}
 \E_xV(X(\sigma_x))-V(x)\leq-\beta\E_x\sigma_x,
 \qquad x\in S\setminus F.
\end{equation}
Then $X$ is positive recurrent.
\end{lemma}

Near the boundary we use the following immediate corollary, taking $\sigma_x$ to be the exit time from $A_i$ for $x\in A_i$.
\begin{corollary}\label{l:averagelyapunov}
Let $X$ be a continuous-time Markov chain on $\mathbb Z_{\geq 0}^d$ with bounded jumps, let $S$ be a closed irreducible component, and let $\mathcal A = \{A_i\}_{i= 0}^{\infty}$ be a partition of $S$ into finite sets. Suppose there is a function $V : S \to \R_{\ge 0}$ with compact sublevel sets satisfying $\max_{x \in A_0} |\mathcal{L}V(x)| < \infty$ and for which there is a finite constant $\beta > 0$ such that  for all $i \in \{1,2,3,\dots\}$, we have
\begin{equ}\label{e:mainbound}
    \E_x\left[ \int_0^{\tau_i} \L V(X(s)) ds\right] \le -\beta \cdot \E_x[\tau_i] \qquad \text{ for all } x \in A_i\,,
\end{equ}
 where $\tau_i := \inf\{t>0~:~X(t) \not \in A_i\}$ are integrable. Then $X$ is  positive recurrent.
\end{corollary}

Throughout, we refer to a function $V$ satisfying the requirements of the above results as a \emph{Lyapunov function}. Note that this result is quite similar to the usual Lyapunov method for showing positive recurrence \cite{meyn12}.  There, the goal is typically to show the existence of a non-negative function $V$ with compact sublevel sets satisfying $\max_{x\in \mathcal K} |\L V(x)|<\infty$ and for which there is a constant $\r >0$ so that $\L V(x) \le -\r$ for all $x\notin \mathcal K$, where $\mathcal K$ is a compact set.  
If we can find such a $V$,  $\r$, and $\mathcal K$, and if there are no absorbing states, then enumerating $S\setminus \mathcal K=\{x_i\}_{i \ge 1}$ and choosing  $A_i = \{x_i\}$, $A_0:= \mathcal K$ we get
\begin{align*}
    \E_{x_i} \left[ \int_0^{\tau_i} \L V(X(s))ds\right] 
    = \mathcal L V(x_i)\E_{x_i}[\tau_i]
    \le -  \r  
    \E_{x_i}[\tau_i]\,.
\end{align*}
Hence, the pointwise condition is a special case of the averaged one of \cref{l:averagelyapunov}.

In fact, the paragraphs above give a partial play-book for how we will use \lref{l:filonov} in the present paper. Henceforth we will restrict our attention to the two-species setting, i.e., $d=2$. Near the $x$ and $y$ axes the sets $A_i$ in \cref{l:averagelyapunov} will consist of lumped sets within which the process may get trapped for some period of time. On the other hand, away from the axes, where no such trapping is possible, we will use the stopping times supplied by the bulk estimate.

\begin{figure}[t]
  \centering
  \resizebox{\textwidth}{!}{%
\begin{tikzpicture}[>=triangle 45]

\begin{scope}
  \def\K{3.5}
  \def\L{13}         

  \fill[blue,fill opacity=0.15] (0,0) rectangle (\L,\K);
  \fill[orange,fill opacity=0.22] (0,0) rectangle (\K,\L);

  \draw[->] (0,0) -- (\L+0.6,0) node[right] {$x_1$};
  \draw[->] (0,0) -- (0,\L+0.6) node[above] {$x_2$};

  \draw[blue!55!black,thick] (0,\K) -- (\L,\K);
  \draw[orange!70!black,thick] (\K,0) -- (\K,\L);

  \draw (\K,0.1) -- (\K,-0.1) node[below] {$K$};
  \draw (0.1,\K) -- (-0.1,\K) node[left] {$K$};

  \node[blue!55!black] at (9.0,1.7) {$B_0(K)$};
  \node[orange!70!black,rotate=90] at (1.7,9.0) {$B_0'(K)$};
  \node at (8.5,8.5) {bulk};

  \node at (\L/2,-1.8) {\small (a) $\gamma=0$};
\end{scope}

\begin{scope}[xshift=16cm]
  \def\K{3.5}
  \def\g{1.4}
  \def\L{13}         

  \fill[blue,fill opacity=0.15]
     (0,0) -- (0,\K) -- (1,\K)
     -- plot[domain=1:\L,samples=90,variable=\x] (\x,{\K+\g*ln(\x)})
     -- (\L,0) -- cycle;
  \fill[orange,fill opacity=0.22]
     (0,0) -- (\K,0) -- (\K,1)
     -- plot[domain=1:\L,samples=90,variable=\x] ({\K+\g*ln(\x)},\x)
     -- (0,\L) -- cycle;

  \draw[->] (0,0) -- (\L+0.6,0) node[right] {$x_1$};
  \draw[->] (0,0) -- (0,\L+0.6) node[above] {$x_2$};

  \draw[blue!55!black,thick]
     (0,\K) -- (1,\K) -- plot[domain=1:\L,samples=90,variable=\x] (\x,{\K+\g*ln(\x)});
  \draw[orange!70!black,thick]
     (\K,0) -- (\K,1) -- plot[domain=1:\L,samples=90,variable=\x] ({\K+\g*ln(\x)},\x);

  \draw (\K,0.1) -- (\K,-0.1) node[below] {$K$};
  \draw (0.1,\K) -- (-0.1,\K) node[left] {$K$};

  \node[blue!55!black] at (9.5,1.7) {$B_\gamma(K)$};
  \node[orange!70!black,rotate=90] at (1.7,9.5) {$B_\gamma'(K)$};
  \node at (10.3,10.8) {bulk};
  \node[blue!55!black,anchor=south] at (9.3,7.7)
        {\small $x_2 = K+\gamma\log x_1$};
  \draw[blue!55!black,->,shorten >=2pt]
        (9.0,7.55) to[bend right=12] (8.4,{\K+\g*ln(8.4)});

  \node at (\L/2,-1.8) {\small (b) $\gamma>0$};
\end{scope}

\end{tikzpicture}%
  }
  \caption{\textbf{The boundary layers} in the two regimes of the parameter $\gamma$
  (taken as $\gamma=0$ for $(1,0)$-SE networks and $\gamma>0$, sufficiently large, for
  $(1,0)$-WE networks). \textbf{(a)} For $\gamma=0$ the layers
  $B_0(K)=\{(n,k)\in\mathbb Z_{\ge0}^2:k\le K\}$ and $B_0'(K)$ are flat strips of
  constant thickness $K$ along the two axes, meeting in a $K\times K$ corner.
  \textbf{(b)} For $\gamma>0$ the layers $B_\gamma(K)$ and $B_\gamma'(K)$ widen
  logarithmically, bounded respectively by $x_2=K+\gamma\log x_1$ and its mirror image
  across the diagonal. In either case the bulk is the complement of the union of the
  two layers, which overlap close to the origin.}
  \label{fig:boundary-layers}
\end{figure}
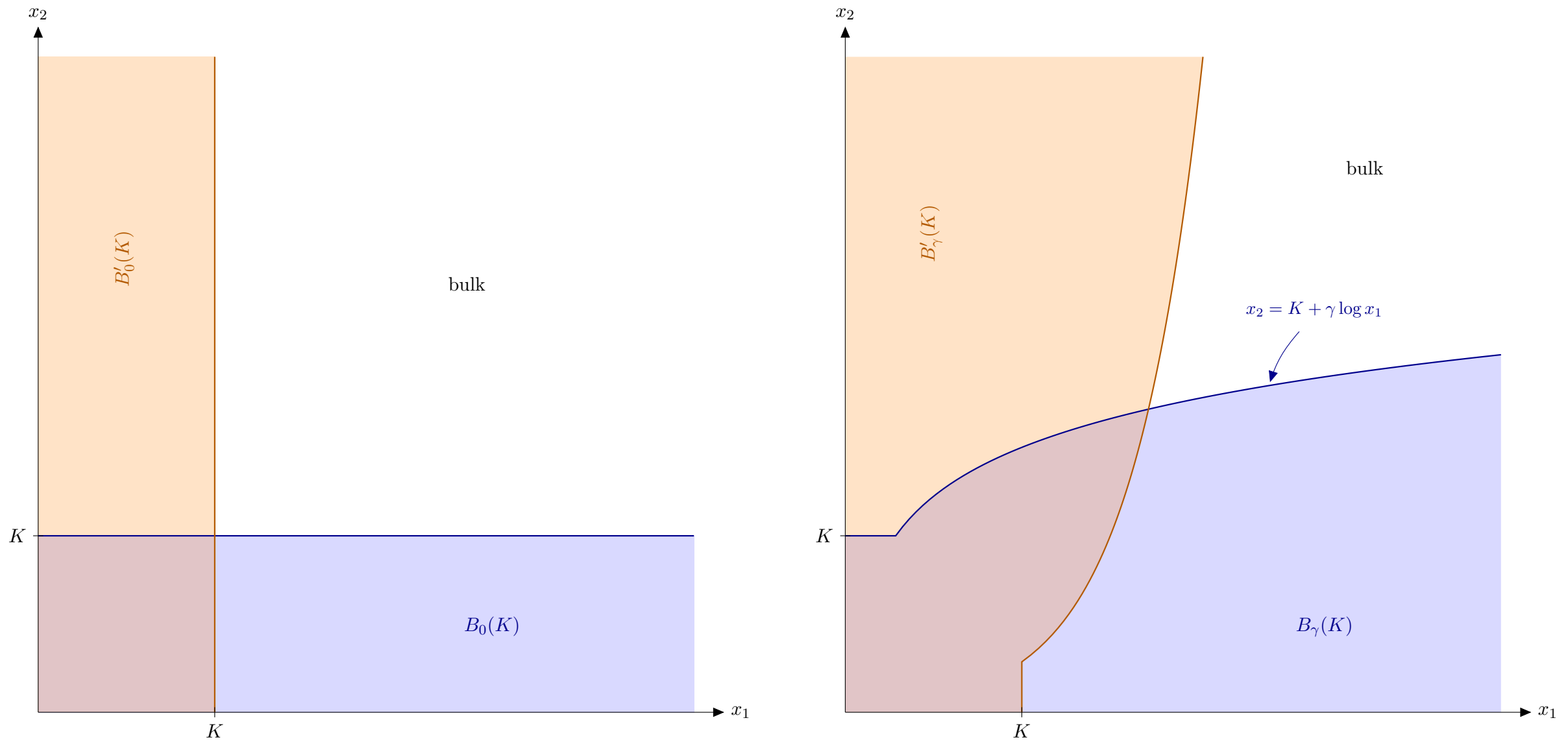

More specifically, to investigate the dynamics of the Markov chain $X$ close to the boundary of the state space, we define the sets
\begin{equs}\label{e:bgammak}
B_\gamma(K) := \{(n,k) \in \mathbb Z_{\geq 0}^2~:~ n\ge1,\ k \leq K + \gamma \log n\}\,,\\
B_\gamma'(K) := \{(k,n) \in \mathbb Z_{\geq 0}^2~:~ n\ge1,\ k \leq K + \gamma \log n\}\,,
\end{equs}
for constants $\gamma \geq 0$, $K>0$ to be chosen later, and define throughout $K_0 = 1+\max_{r}\{\|y_r\|_1+\|\zeta_r\|_1\}$.
These boundary layers and the bulk that separates them are illustrated in \fref{fig:boundary-layers}.

For $K\ge2K_0$, $\gamma_1,\gamma_2\ge0$ and bounded functions $f_i:\mathbb Z_{\ge0}\to\mathbb Z$ with $f_i(k)=0$ for $k\ge K-2K_0$, define for $n\ge1$
  \[
  \begin{gathered}
    h_{i,n}=\lfloor K+\gamma_i\log n\rfloor,\qquad i\in\{1,2\},\\
    \mathcal F_n:=\{(n+f_1(k),k)\in S:0\leq k\leq h_{1,n}\},\\
    \mathcal F'_n:=\{(k,n+f_2(k))\in S:0\leq k\leq h_{2,n}\},
  \end{gathered}
  \]
Choose $M$ sufficiently large that outside $A_0=S\cap[0,M]^2$ the two boundary layers are disjoint and each is partitioned by its respective fibers. For every integer $\alpha\ge1$, set $V=0$ on $A_0$ and define $V$ on $S\setminus A_0$ by
\begin{equ}\label{e:lyapunovfunction}
    V(x)=
    \begin{cases}
        \frac{1}{\alpha+1}(u(x_1)+u(x_2))^{\alpha+1}, &x\notin B_{\gamma_1}(K)\cup B'_{\gamma_2}(K),\\
        \frac{1}{\alpha+1}\pc{u(n)+u(h_{1,n})}^{\alpha+1}, &x\in\mathcal F_n,\\
        \frac{1}{\alpha+1}\pc{u(n)+u(h_{2,n})}^{\alpha+1}, &x\in\mathcal F'_n,
    \end{cases}
  \end{equ}
where the function $u$ is defined as:
\begin{equ}\label{e:entropy-u}
 u(t)=t\log t-t\quad(t>0),\qquad u(0)=0.
\end{equ}

\begin{proposition}\label{p:finalboundary}
For any weakly reversible CRN with two species and all sufficiently large integers $K>2K_0$, there exist $f_1,f_2,\gamma_1,\gamma_2$ as above, with
\[
 \gamma_i=0\quad(e_i\text{-SE}),\qquad \gamma_i>0\quad(e_i\text{-WE}),\qquad e_1=(1,0),\quad e_2=(0,1),
\]
such that for every sufficiently large integer $\alpha$ there exist $\beta>0$ and $n_0$ for which
\begin{equ}\label{e:mainbound2}
 \E_x\left[\int_0^{\tau_n}\L V(X(s))\,ds\right]\le -\beta\E_x\tau_n,
 \qquad x\in\mathcal F_n,\quad n\ge n_0.
\end{equ}
Here $\tau_n=\inf\{t>0:X(t)\notin\mathcal F_n\}$ is integrable. The same estimate holds for $\mathcal F'_n$ and its exit time.
\end{proposition}
In the horizontal arguments below, we write $f=f_1$, $\gamma=\gamma_1$ and $h_n=h_{1,n}$. The corresponding vertical arguments use $f_2$, $\gamma_2$ and $h_{2,n}$.

The corresponding bulk estimate is proved in Section~\ref{s:patch}.

\begin{proposition}\label{p:totalbulk0}
  For each sufficiently large $K$, take $f_1,f_2,\gamma_1,\gamma_2$ and $\alpha$ as in \pref{p:finalboundary}. There exist $\beta>0$ and a finite set $\mathcal K\subset S$ such that every $x\in S\setminus(B_{\gamma_1}(K)\cup B_{\gamma_2}'(K)\cup\mathcal K)$ admits an integrable stopping time $\sigma_x>0$ satisfying
  \begin{equ}\label{e:bulkestimate}
    \E_x V(X(\sigma_x))-V(x)\leq-\beta\E_x\sigma_x.
  \end{equ} 
\end{proposition}

\begin{proof}[Proof of \tref{t:main}]
Choose $K$ sufficiently large for Propositions~\ref{p:finalboundary} and~\ref{p:totalbulk0}. For this $K$, take $f_1,f_2,\gamma_1,\gamma_2$ from \pref{p:finalboundary} and then take $\alpha$ sufficiently large. On each boundary fiber, take $\sigma_x=\tau_n$. By stopped Dynkin's formula, \eref{e:mainbound2} gives \eref{e:filonov-general}. In the bulk use \eref{e:bulkestimate}. With the smaller of the two drift constants, \lref{l:filonov} gives positive recurrence.
\end{proof}

\subsection{Outline of the paper} The remainder of the paper will be organized as follows. In \sref{s:boundary} and \sref{s:boundary2} we prove dissipation close to the axes of the state space, establishing \pref{p:finalboundary}. This entails carrying out the procedure outlined at the end of \sref{s:main} to identify the sets $\mathcal F_n$ and the corresponding exit estimates. 
Finally, in  \sref{s:patch} we the dissipation in the bulk stated in Proposition~\ref{p:totalbulk0}. This last section is independent of Sections \sref{s:boundary} and \sref{s:boundary2}, and use independent notations.

\section{Boundary layer: Large deviations and coarse-graining}
\label{s:boundary}
 In this section we analyze the behavior of the Markov chain in a region close to one of the axes. We will limit our discussion to the $x$-axis, as the same results hold for the $y$-axis by exchanging the two coordinates in the upcoming analysis.
 
 We introduce here the techniques from Large Deviation Theory that we will use in the next section to prove Proposition~\ref{p:finalboundary}. 
 In particular, we will consider the embedded discrete-time Markov chain $Y(s)$,  $s\in\Z_{\geq 0}$, associated to our continuous-time Markov chain $X(t)$, $t\in\R_{\geq 0}$ and provide subtle estimates pertaining to the asymptotics of likely jump directions. To introduce the main ideas behind this strategy let us consider our main example: 
 
\begin{example}[Example 1 continued]\label{ex:mainexample} In the regime $n\gg 1$ the dynamics of $Y(s)$ with $Y(0) = (n,0)$ is restricted, with overwhelming probability, to the set $\OO_n:=\{(n,0),(n+1,1),(n+3,2)\}$, \ie it will remain on the
  periodic path
  \begin{equ}
    g^* := (n,0)\to(n+1,1)\to(n+3,2)\to(n,0)\to \dots
  \end{equ}
  until a low probability jump takes the process outside the set.
This implies that on finite time intervals the system does not appear to dissipate, at least with high probability. Still, each time the process enters the state $(n+1,1)$, there is a low probability that the next  jump is to the state $(n+2,2)$ instead of $(n+3,2)$.  The probability of this transition is
\begin{align*}
         \frac{\lambda_{0 \to A+B}((n+1,1))}{\lambda_{0 \to A+B}((n+1,1)) + \lambda_{A+B\to 3A+2B}((n+1,1))}
         = 
         \frac{\kappa_{0 \to A+B}}{\kappa_{0 \to A+B} + \kappa_{A+B\to 3A+2B}(n+1)} = \mathcal{O}(n^{-1}).
\end{align*}
In this case, we observe a path of the form
\begin{equ}
  g^- :=  \dots \to(n,0) \to(n+1,1)\rightsquigarrow(n+2,2) \to (n-1,0) \to\dots
\end{equ}
where $\leadsto$ denotes that a transition with low probability occurred. This transition ``pulls'' the process $Y(s)$ from $\OO_n$ to $\OO_{n-1}$, i.e., towards the origin. 
We note that it is also possible to ``push'' the process from $\OO_n$ to $\OO_{n+1}$ (away from the origin), e.g., through a path
\begin{equ}
  g^+ :=  \dots \to(n,0) \to(n+1,1)\to(n+3,2)\leadsto(n+5,3) \to (n+2,1) \to (n+4,2) \to (n+1,0) \to\dots
\end{equ}
where again $\leadsto$ denotes an unlikely transition. However, in this case this transition occurs with probability 
\begin{align*}
  \frac{\lambda_{A+B \to 3A+2B}((n+3,2))}{\lambda_{0 \to A+B}((n+3,2)) + \lambda_{A+B \to 3A+2B}((n+3,2)) + \lambda_{3A+2B \to 0}((n+3,2))} = \mathcal{O}(n^{-2}), 
\end{align*}
i.e., an order of magnitude less than the one from $\OO_n$ to $\OO_{n-1}$.

We note that in the example above the transition $\OO_n \to \OO_{n-1}$ is not only more likely than $\OO_n \to \OO_{n+1}$ on a finite time interval, but it also determines the \emph{most} likely exit path of the process from the set $\OO_n$.  In fact,
\begin{equ}
  \lim_{n \to \infty} \mathbb P_{(n,0)}(Y(\sigma) \in \OO_{n-1}) = 1\qquad\text{and}\qquad \lim_{n \to \infty}\mathbb P_{(n,0)}(Y(\sigma) \in \OO_{n+1}) = 0
\end{equ}
where $\sigma := \inf\{s \in \Z_{\geq0}~:~Y(s) \in \OO_{n-1}\cup \OO_{n+1}\}$. 
Now consider a candidate Lyapunov function of the following form
\begin{equ}
  V(x) := \frac1{\alpha + 1}u(n)^{\alpha+1} \qquad \text{for each \, }x \in \OO_n,
\end{equ}
with $u$ as in \eref{e:entropy-u} and some $\alpha\ge1$ to be determined. In this case,
\begin{equs}
  \mathbb E_{x}&[V(Y(\sigma))] - V(x) \\
  & = \left[\frac{u(n-1)^{\alpha+1} -  u(n)^{\alpha+1}}{\alpha + 1}\right]\mathbb P_{x}(Y(\sigma) \in \OO_{n-1})  + \left[\frac{u(n+1)^{\alpha+1} -  u(n)^{\alpha+1}}{\alpha + 1}\right]\mathbb P_{x}(Y(\sigma) \in \OO_{n+1})\\
  & \leq -\frac {1} 2 u(n)^{\alpha}\label{e:dissipation}
\end{equs}
where the last inequality holds for our specific choice of $u(n)$. 
For $x=(n,0)$, take $\sigma_x$ to be the time of the $\sigma$-th jump. Then $X(\sigma_x)=Y(\sigma)$ and $\E_x\sigma_x\leq\kappa_-^{-1}\E_x\sigma$. Estimate \eqref{e:dissipation} alone gives no bound on $\E_x\sigma_x$, as required in \eqref{e:filonov-general}. 
Be aware that in this example we purposely avoided the analysis of further exit paths from the set $\OO_n$, with probability of lower order than the ones discussed above. This serves to highlight the main idea of our technique without introducing too many details. The correction term $u(h_n)$ in~\eqref{e:lyapunovfunction} is used to patch the Lyapunov function with the bulk entropy, and doesn't change the leading order in $n$ of these estimates. \hfill $\triangle$
\end{example}

Our strategy to define a Lyapunov function on the boundary layer generalizes the one adopted in the example  above. There are three main takeaways from the above example, each one corresponding to a lemma in our main proof.
\begin{enumerate}
  \item For $n \gg 1$ the most likely transitions of the network may result in the emergence of ``trapping'' sets -- with high probability at the right timescale -- of the Markov chain $Y(s)$. If such sets exist, then by considering timescales longer than $\mathcal O(1)$ these trapping  sets can be \emph{coarse-grained} into sets $\OO_n$, and the dynamics of the process is well approximated by the (most likely) transitions \emph{between} these sets. We note that this might result in the emergence of \emph{nested} trapping sets (i.e., trapping sets of trapping  sets), and therefore that the coarse-graining process might have to be repeated. This corresponds to considering the dynamics of the process on longer and longer timescales.
  
  \item Repeating the above process sufficiently many times allows us to identify translation-invariant sets $\mathcal F_n$, which we will refer to as \emph{fibers}, each containing a union of coarse-grained sets from the previous point. The fibers will  induce a partial order on points in the boundary layer so that
  the Markov chain $Y(s)$ will exit from $\mathcal F_n$ into  $\mathcal F_{n'}$ with $n' < n$, with probability bounded below away from zero, while the probability of an exit to a larger label tends to zero as $n\to\infty$. 
  
  \item Finally, we define the Lyapunov function $V$ as being constant
  on each $\mathcal F_n$, while growing sufficiently fast as a function of $n$.
  This yields the desired dissipation bound through a computation analogous to \eref{e:dissipation}.
\end{enumerate}

We will expand on these points in more detail below.   To do so, we first need the following definition.

\begin{definition}
    Let $\mathcal G$ be a directed graph. We say that a subset of nodes $\OO$ is a strongly connected set if it is a singleton, or if for all $\omega, \omega'\in\OO$ there is a directed path in $\mathcal G$ from $\omega$ to $\omega'$ which is fully contained in $\OO$ (meaning that every node of the directed path is an element of $\OO$). We say that $\OO$ is a maximal strongly connected set if it is not properly contained in any other strongly connected set.
\end{definition}

 To detail the steps outlined in points (1)-(3), we will associate those steps with operations on graphs. We present these operations below for the case $\gamma = 0$ (covering $(1,0)$-strongly endotactic networks) to build intuition and introduce some notation, and leave the discussion of the general case to the next subsections. \fref{fig:coarsegraining} illustrates the procedure on a network that requires two coarse-grainings before its dissipation becomes visible.
\begin{enumerate}
  \item The coarse-graining procedure is performed as follows. Set $k = 0$ and consider the graph $\mathcal G^{(0)}$ of ``dominant'' transitions of the Markov chain $Y(s)$ for $n \gg 1$: for every state $x\in\sz$ we only consider transitions out of it of the form $x\to x+y'-y$ with $y\to y'$ being a dominant reaction at level $x_2$ (defined formally below in Definition~\ref{d:cost}). Note that by definition, $x\to x'$ is an edge of $\mathcal G^{(0)}$ if and only if $x+(m,0)\to x'+(m,0)$ is an edge of $\mathcal G^{(0)}$, for all $m\in\Z_{\ge 0}$. This implies that the maximal strongly connected sets of $\mathcal G^{(0)}$ (and of $\mathcal G^{(k)}$ as well, as defined below) are translation invariant in the $x$ direction. 
  We then proceed with the recursive step as follows, starting with $k=0$:
  \begin{enumerate}
    \item  We group or \emph{coarse-grain}

    the maximal strongly connected sets of $\mathcal G^{(k)}$ into sets $\OO_{(n,i)}^{(k)}$, where $(n,i)$ denotes the element of the strongly connected set with the lowest second coordinate.

    \item Then, we define the next iterate $\mathcal G^{(k+1)}$ by defining its vertices as the previously identified coarse-grained sets and its directed edges as the transitions with leading order probabilities between these sets.
    
  \item If the graph ${\mathcal G}^{(k+1)}$ is directed acyclic we stop the procedure, else we return to point (a) substituting $k$ with $k+1$.
  \end{enumerate}
  
  \item We now consider the output of the previous iteration, denoted as ${\mathcal G}^{*}$. Here, we prove that we can partially order its vertices, denoted by $\OO_{(n,i)}^*$, by grouping them into \emph{fibers}, i.e., sets $\mathcal F_n$ that are translation invariant in the $x$ direction and contain a node $\OO_{(n,i)}^{*}$ of ${\mathcal G}^{*}$ for each $i$, such that there does not exist a transition $\OO \to \OO'$ in ${\mathcal G}^{*}$ if $\OO \in \mathcal F_n$, $\OO' \in \mathcal F_m$ with $m>n$. 
  \item We finally combine the above construction with our definition \eref{e:lyapunovfunction} and an estimate on the exit times of $Y$ from $\mathcal F_n$ to prove Proposition \ref{p:finalboundary}.
\end{enumerate}
The remainder of this section realizes the above steps. More specifically, we first inscribe the asymptotic behavior of $Y$ in the context of large deviations theory (\sref{s:ldp}), which provides a natural framework to study the exit problem from $\omega$-limit sets, and establish some preparatory bounds (\sref{s:prep}). Then, we establish the required exit time estimates recursively on the operations (1a), (1b) (respectively, Sections~\ref{s:wgraphs1} and \ref{s:wgraphs2}). Combining this with the structure of weakly reversible CRNs (\sref{s:fibers1}) we then show the existence of fibers $\mathcal F_n$ that capture the dissipative behavior of the coarse-grained dynamics (\sref{s:fibers2}), realizing (2). Finally, in \sref{s:boundaryfinal} we combine the above estimates in (3).

    \begin{figure}[p]
  \centering
  \colorlet{ClL}{blue!65!black}\colorlet{ClC}{green!55!black}\colorlet{ClR}{orange!80!black}%
  \colorlet{sLl}{blue!70!black}\colorlet{sLh}{orange!90!black}\colorlet{sCl}{red!65!black}%
  \colorlet{sCh}{violet!75!black}\colorlet{sRl}{green!55!black}\colorlet{sRh}{cyan!60!black}%
  \tikzset{
    cgodot/.style={circle,fill=black,inner sep=1.6pt},
    cglvl/.style={gray!20,line width=0.4pt},
    cgdom/.style={->,>=Latex,line width=0.9pt,black},
    cgfeed/.style={->,>=Latex,line width=0.9pt,black},
    cgrare/.style={->,>=Latex,line width=1.6pt,black},
    cgtoward/.style={->,>=Latex,line width=1.9pt,black},
    cggroup/.style={dashed,line width=1pt,dash pattern=on 3pt off 2.2pt,line cap=round,line join=round},
  }%
  \newcommand{\cgGrid}{\foreach \k in {0,1,2,3,4,5}{\draw[cglvl] (-2.9,\k) -- (13.6,\k);
      \node[font=\scriptsize,gray!55,anchor=east] at (-2.9,\k){$k{=}\k$};}}%
  \newcommand{\cgPts}{%
    \node[cgodot](L0) at (-2,0){};\node[cgodot](L1) at (-1,1){};\node[cgodot](L2) at (0,2){};\node[cgodot](L3) at (3,3){};%
    \node[cgodot](C0) at (0,0){};\node[cgodot](C1) at (1,1){};\node[cgodot](C2) at (2,2){};\node[cgodot](C3) at (5,3){};%
    \node[cgodot](R0) at (2,0){};\node[cgodot](R1) at (3,1){};\node[cgodot](R2) at (4,2){};\node[cgodot](R3) at (7,3){};%
    \node[cgodot](bL) at (6,4){};\node[cgodot](bC) at (8,4){};\node[cgodot](bR) at (10,4){};%
    \node[cgodot](cL) at (9,5){};\node[cgodot](cC) at (11,5){};\node[cgodot](cR) at (13,5){};}%
  \newcommand{\cgBulk}{\draw[cgfeed] (cL)--(bL);\draw[cgfeed] (cC)--(bC);\draw[cgfeed] (cR)--(bR);%
    \draw[cgfeed] (bL)--(L3);\draw[cgfeed] (bC)--(C3);\draw[cgfeed] (bR)--(R3);}%
  \begin{tikzpicture}
    \node[draw,rounded corners,fill=gray!5,inner sep=5pt,font=\small]
      {$0 \rightleftharpoons A{+}B \qquad 2A{+}2B \rightleftharpoons 5A{+}3B$};
  \end{tikzpicture}\\[2pt]
  \resizebox{0.82\textwidth}{!}{%
  \begin{tikzpicture}[x={(0.8cm,0)},y={(0,0.66cm)}]
    \cgGrid \cgPts
    \foreach \a/\b in {L0/L1,C0/C1,R0/R1}{\draw[cgdom] (\a) to[bend left=18] (\b);\draw[cgdom] (\b) to[bend left=18] (\a);}
    \foreach \a/\b in {L2/L3,C2/C3,R2/R3}{\draw[cgdom] (\a) to[bend left=8] (\b);\draw[cgdom] (\b) to[bend left=8] (\a);}
    \cgBulk
    \node[font=\scriptsize,below=1pt of C0]{$(n,0)$};
  \end{tikzpicture}}\\[1pt]
  {\footnotesize\textbf{(a)} $\mathcal G^{(0)}$ (time $O(1)$): two bounded $2$-cycles per column (low $\pm(1,1)$, high $\pm(3,1)$). The bulk feeds into the high cycle. Nothing dissipates.}\\[3mm]
  \resizebox{0.82\textwidth}{!}{%
  \begin{tikzpicture}[x={(0.8cm,0)},y={(0,0.66cm)}]
    \cgGrid
    \draw[cggroup,sLl] (L0)--(L1);\draw[cggroup,sCl] (C0)--(C1);\draw[cggroup,sRl] (R0)--(R1);
    \draw[cggroup,sLh] (L2)--(L3);\draw[cggroup,sCh] (C2)--(C3);\draw[cggroup,sRh] (R2)--(R3);
    \cgPts \cgBulk
    \foreach \u/\d in {L1/L2,C1/C2,R1/R2}{\draw[cgrare] (\u) to[bend left=28] (\d);\draw[cgrare] (\d) to[bend left=28] (\u);}
    \foreach \x/\y/\c/\lab in {%
      -2.45/-1.5/sLl/{$\mathbf X^{(0)}_{(n-2,0)}$}, 2.45/-1.5/sCl/{$\mathbf X^{(0)}_{(n,0)}$}, 7.35/-1.5/sRl/{$\mathbf X^{(0)}_{(n+2,0)}$},
      -2.45/-2.3/sLh/{$\mathbf X^{(0)}_{(n,2)}$}, 2.45/-2.3/sCh/{$\mathbf X^{(0)}_{(n+2,2)}$}, 7.35/-2.3/sRh/{$\mathbf X^{(0)}_{(n+4,2)}$}%
    }{\draw[cggroup,\c] (\x,\y) -- ++(0.8,0) node[anchor=west,black,xshift=2pt,font=\scriptsize]{\lab};}
  \end{tikzpicture}}\\[1pt]
  {\footnotesize\textbf{(b)} $\mathcal G^{(1)}$ (time $O(n)$): each $2$-cycle is one coarse-grained set (a dashed coloured line in the legend). The leading $O(n^{-1})$ transitions (solid) link a low set to the high set two columns right and back, by $0\to A{+}B$ and $A{+}B\to0$. They shift the index $+2$ then $-2$ and cancel, so the pair is still a cycle.}\\[3mm]
  \resizebox{0.82\textwidth}{!}{%
  \begin{tikzpicture}[x={(0.8cm,0)},y={(0,0.66cm)}]
    \cgGrid
    \draw[cggroup,ClL] (L0)--(L1)--(L2)--(L3);
    \draw[cggroup,ClC] (C0)--(C1)--(C2)--(C3);
    \draw[cggroup,ClR] (R0)--(R1)--(R2)--(R3);
    \cgPts \cgBulk
    \draw[cgtoward,shorten >=2.5pt,shorten <=2.5pt] (R0) -- (C0);
    \draw[cgtoward,shorten >=2.5pt,shorten <=2.5pt] (C0) -- (L0);
    \draw[cgtoward,shorten <=2.5pt] (L0) -- (-2.8,0);
    \node[font=\scriptsize,gray!55,anchor=east] at (-3.0,0.5){origin\,$\leftarrow$};
    \foreach \x/\y/\c/\lab in {%
      -2.45/-2.1/ClL/{$\mathbf X^{(1)}_{(n-2,0)}$}, 2.45/-2.1/ClC/{$\mathbf X^{(1)}_{(n,0)}$}, 7.35/-2.1/ClR/{$\mathbf X^{(1)}_{(n+2,0)}$}%
    }{\draw[cggroup,\c] (\x,\y) -- ++(0.8,0) node[anchor=west,black,xshift=2pt,font=\scriptsize]{\lab};}
  \end{tikzpicture}}\\[1pt]
  {\footnotesize\textbf{(c)} $\mathcal G^{(2)}$ (time $O(n^2)$): merging each balanced cycle gives the sets $\mathbf X^{(1)}_{(m,0)}$, shown two columns apart. The leading exit mechanism has order $O(n^{-2})$ and shifts the index from $m$ to $m-2$, toward the origin.}
  \caption{\textbf{A network requiring two coarse-grainings to expose dissipation}, for $0\rightleftharpoons A{+}B$, $2A{+}2B\rightleftharpoons 5A{+}3B$ (schematic, with state $(x_1,x_2)=(\#A,\#B)$, column $x_1=n\gg1$, level $x_2=k$). Two stacked reversible cycles form a \emph{balanced} pair of trapping sets: the leading transitions cancel, so no net motion toward the origin is visible after one coarse-graining (b), and it emerges only after the second (c). Coarse-grained sets are drawn as dashed coloured lines joining their states. Solid arrows represent dominant transitions between states or sets at each stage; times are measured in jumps of $Y$.}
  \label{fig:coarsegraining}
\end{figure}

\subsection{Large deviations}\label{s:ldp}
As we aim to establish the dissipation estimate \eref{e:mainbound} outside of a sufficiently large compact set in $\Z^2_{\geq0}$, in this section we will study the behavior of $X$ (or, in fact, $Y$) in $\sg$ when $n \gg 1$, so that the parameter $n$ can be considered as a scaling parameter. Hence, all the inequalities in this section will be considered to be valid for large enough $n$, even when not explicitly stated.

Recalling \rref{r:active} we define complexes that are \emph{active} at level  $k \in \Z_{\geq0}$.
\begin{definition}\label{def:activecomplex}
  A complex $y$ is said to be \emph{active} at level $k \in \Z_{\ge 0}$ if $\lambda_y((n,k))>0$ for large enough $n$. 
  We denote the set of complexes that are active at level $k$ by $\C(k)$.
	We set 
	\begin{equ}
		z(k):=\arg\max\{z_1\, :\, z\in \cal{C}(k)\}\, , 
	\end{equ}
	and for $v\in \Z^2$, we set 
	\begin{equ}
		y(k,v):=\arg\max\{y_1\, :\, y\in \cal{C}(k)\quad \text{and}\quad\exists\, y\to y'\in \cal{R}(v)\}\, ,
	\end{equ}
	where in both case, in case of multiple maximizers, we choose one with the largest second coordinate. 
\end{definition}

We recall that throughout we denote by $\{Y(s)\}_{s \in \Z_{\geq0}}$ the embedded (discrete-time) Markov chain associated to the continuous-time Markov chain $X$. Then the transition probabilities for the chain $Y$ are given in \eqref{e:transprob}.

Combining the above definition with \eref{e:rates} we notice that for any $(n,k)\in\sg$ and for any $v\in\Z^2$ such that $p((n,k), (n,k)+v)>0$ we have
\begin{equs}
  p((n,k),(n,k)+v) &= \frac{\sum_{y\to y'\in\Re(v)}\kappa_{y \to y'} \frac{k!}{(k-y_2)!}\frac{n!}{(n-y_1)!}}{\sum_{z \to z' \in \Re}\kappa_{z \to z'}\frac{k!}{(k-z_2)!}\frac{n!}{(n-z_1)!}}\\
   & \leq \frac{|\Re|\kappa_+(k-z_2(k))! }{\kappa_- (k-y_2(k,v))!} \frac{(n-z_1(k))!}{(n-y_1(k,v))!} \\
	&\leq a(K+\gamma\log n)  n^{-(z_1(k) - y_1(k,v))}\label{e:ratesub}
\end{equs}
for a polynomial function $a(\cdot)$ of bounded degree that depends only on $(\mathcal S, \mathcal C, \mathcal R)$ and $\{\kappa_{r}\}_{r \in \mathcal R}$.
Similarly, upon possibly increasing the maximal degree and the coefficients of $a(\cdot)$, we can derive a lower bound with a matching exponent for $n$: 
\begin{equs}
  p((n,k),(n,k)+v)  &\geq \frac{\kappa_-(k-z_2(k))! }{|\Re|\kappa_+ (k-y_2(k,v))!} \frac{(n-z_1(k))!}{(n-y_1(k,v))!} \\
	&\geq a(K+\gamma\log n)^{-1}  n^{-(z_1(k) - y_1(k,v))}\,,\label{e:rateslb}
\end{equs}
whenever $n\geq K_0$.
We note that when $\gamma = 0$ the above bounds are \emph{precise}: the multiplicative pre-factor $a$ in front of the probabilities $p(\cdot)$ becomes an $n$-independent constant.

\begin{remark}\label{r:ldp}
  The above bounds show that, for fixed $K$, as $n\to \infty$  the transition probabilities of the Markov process $Y$ satisfy for every $v\in \Z^2$ the large deviation bound
  \begin{equ}\label{e:stepldp}
    \lim_{n \to \infty} (\log n)^{-1} \log p((n,k),(n,k)+v) = - J(v,k)
  \end{equ}
  with rate function
  \begin{equ}\label{e:steprf}
    J(v,k) := \begin{cases}(z_1(k) - y_1(k,v)) \quad &\text{if}\quad  \C(k)\cap\{y: \exists y\to y'\in \cal{R}(v)\}\neq \emptyset \\
    \infty \qquad & \text{else }\end{cases} \,.
  \end{equ}

  From this it immediately follows that one can establish similar estimates on trajectories:
  \begin{equ}\label{e:pathldp}
    \lim_{n \to \infty} (\log n)^{-1} \log \mathbb P((Y(s))_{s=0}^S \in \Gamma_S) = -\inf_{g \in \Gamma_S} I(g)
  \end{equ}
  where $\Gamma_S$ is the set of discrete-time trajectories of $Y$ on a finite horizon $S \in \Z_{\geq0}$ of the form
  \begin{equ}
  g = \{g_0 \to \dots \to g_S~:~g_s \in \NN \}
  \end{equ}
   and
  \begin{equ}\label{e:pathrf}
    I(g) := \sum_{s=0}^{S-1} J(g_{s+1}-g_s , (g_s)_2)
  \end{equ}
  Note that the above definition can be extended to any directed graph $g$ of transitions in $\mathbb Z_{\geq 0}^2$.
\end{remark}
The above remark motivates the definition of the (large deviations) \emph{cost} of a transition of the process $Y(s)$ and of a path:
\begin{definition}\label{d:cost}
  The \emph{cost} of a transition $v$ at level $k$ is $J(v,k)$ as defined in \eref{e:steprf}, and the cost of a directed graph $g$ in $\Z^2_{\geq0}$ is
  \begin{equation*}
  I(g):= \sum_{x\to x'\in E(g)} J(x'-x , x_2).
  \end{equation*}
  We say that a transition $v$ is \emph{dominant} at level $k$ if it has zero cost, \ie $J(v ,k) = 0$.
\end{definition}
\noindent Note that the costs of both transitions and paths are \emph{independent} on reaction rate constants $\kappa_{ y \to y'}$. 

\subsection{Preparatory definitions and exit estimates}\label{s:prep}
We now define boxes (or rather, a $n$-dependent frame) as
\begin{equs}
\ss(K)  & := \{(j,k ) \in\mathbb Z_{\geq0}^2~:~0\leq k\leq\lfloor K+\gamma\log n\rfloor, j \in [n- 2 K \Kz, n+2K \Kz]\},\\
\bar B_{\gamma,n}(K)  & := \{(j,k ) \in \mathbb Z_{\geq 0}^2 ~:~k \in [0, (K+ \gamma \log n)+\Kz ], j \in [n- 2 K \Kz, n+2K \Kz]\},
\end{equs}
for $n \geq 2KK_0$ and any $\gamma \geq 0$.
We denote points inside such $n$-dependent ``frames'' with a footer-$n$ notation, i.e. for any $n \in \mathbb N$ we define $(j+n,k) =: (j,k)_n \in \ss(K)$.

We now proceed to partition the set $B_{\gamma,n}(K)$ into $\{\XX_m\}$ for the upcoming analysis. This partition will depend on the value of $\gamma$, which in turn will be chosen to be positive for $(1,0)$-WE networks and equal to $0$ for $(1,0)$-SE networks.
\begin{definition}\label{d:partitionofB}
  We define the partition $\{\XX_m\}$ of $\bar B_{\gamma,n}(K)$, with $m$ varying in the set 
  $${m\in [-2KK_0, 2KK_0]\times [0, K+1]},$$
  as follows: write $H_n=\lfloor K+\gamma\log n\rfloor$. For any value of $\gamma\geq0$ we define
  \begin{equ}
    \XX_{(j,k)} := \{(j,k)_n\}\qquad \text{ for any }j \in [-2K\Kz, 2K\Kz], \quad k \in [0, K)\,.
  \end{equ}
    Furthermore, for any $j \in [-2K\Kz, 2K\Kz]$,
  \begin{description}
    \item[if $\gamma = 0$]\,\\ $\XX_{(j,K)} := \{(j,K)_n\}$ and $\XX_{(j,K+1)} := \{(j,k)_n~:~k \in [K+1, K+\Kz]\}$.
    \item[if $\gamma > 0$] \,\\$\XX_{(j,K)} := \{(j,k)_n:k \in [K,H_n]\}$\\
    and $\XX_{(j,K+1)} := \{(j,k)_n:H_n<k\le H_n+\Kz\}$.
  \end{description}
\end{definition}

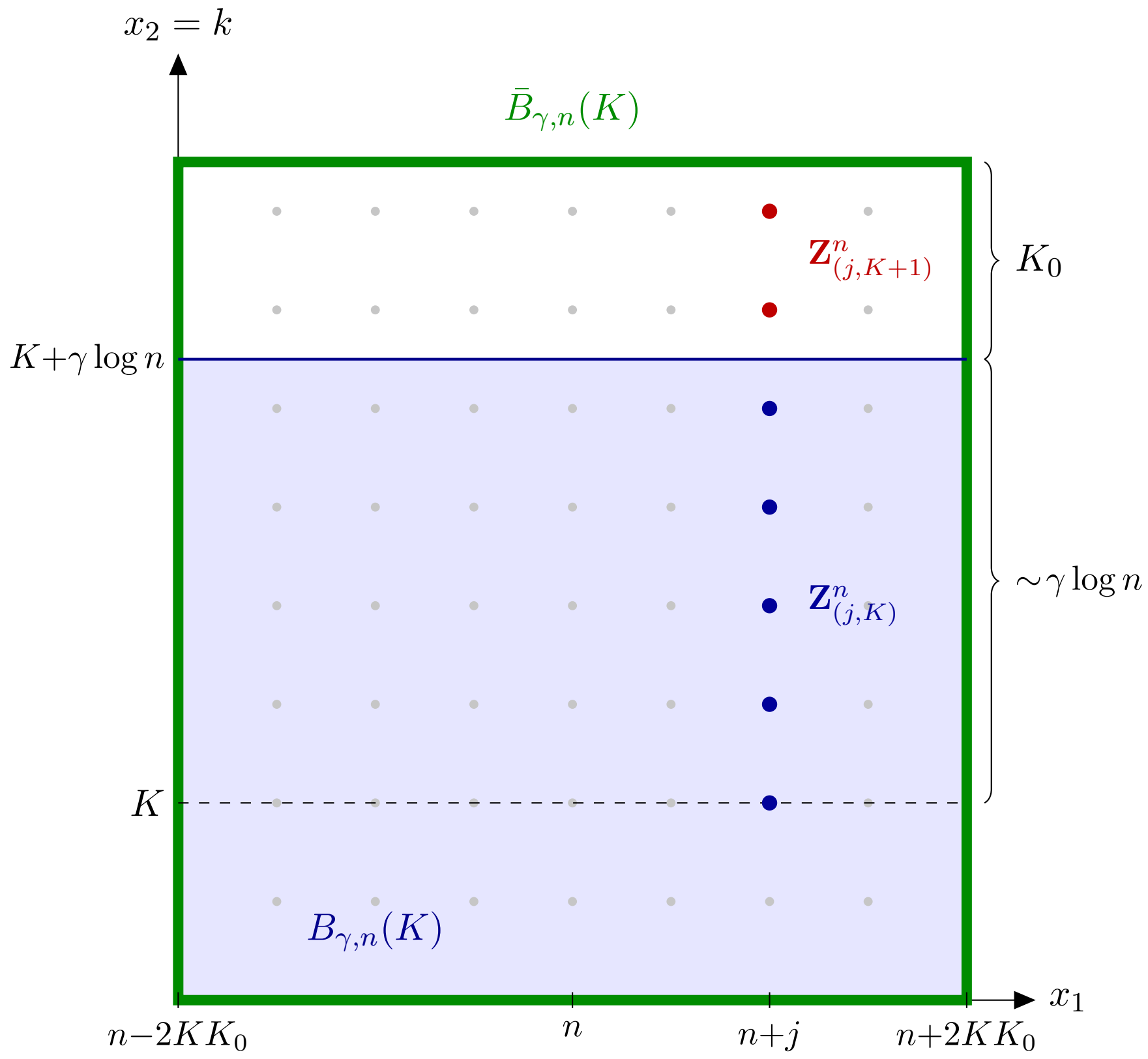
\begin{figure}[t]
  \centering
\begin{tikzpicture}[>=triangle 45]
\begin{scope}
  \def\Wm{8}            
  \def\Klev{2}          
  \def\Lflat{6.5}       
  \def\Uflat{8.5}       

  \fill[blue!10] (0,0) rectangle (\Wm,\Lflat);

  \foreach \i in {0,...,8}{
     \foreach \k in {0,...,8}{ \fill[gray!45] (\i,\k) circle (1.3pt); }
  }

  \draw[->] (0,0) -- (\Wm+0.7,0) node[right] {$x_1$};
  \draw[->] (0,0) -- (0,9.6) node[above] {$x_2=k$};

  \draw[green!55!black,line width=3pt] (0,0) rectangle (\Wm,\Uflat);

  \draw[blue!55!black,thick] (0,\Lflat) -- (\Wm,\Lflat);  
  \draw[dashed] (0,\Klev) -- (\Wm,\Klev);                 
  \node[anchor=east] at (0,\Klev) {$K$};
  \node[anchor=east] at (0,\Lflat) {\small $K{+}\gamma\log n$};

  \def\jj{6}
  \foreach \k in {2,3,4,5,6}{ \fill[blue!60!black] (\jj,\k) circle (2.2pt); } 
  \foreach \k in {7,8}{       \fill[red!75!black]  (\jj,\k) circle (2.2pt); } 
  \node[blue!60!black,anchor=west] at (\jj+0.25,4) {\small $\mathbf Z^n_{(j,K)}$};
  \node[red!75!black,anchor=west]  at (\jj+0.25,7.5) {\small $\mathbf Z^n_{(j,K+1)}$};

  \draw[decorate,decoration={brace,amplitude=4pt}] (\Wm+0.18,\Uflat) -- (\Wm+0.18,\Lflat)
        node[midway,right=5pt] {$\Kz$};
  \draw[decorate,decoration={brace,amplitude=4pt}] (\Wm+0.18,\Lflat) -- (\Wm+0.18,\Klev)
        node[midway,right=5pt] {\small $\sim\!\gamma\log n$};

  \draw (0,0.08) -- (0,-0.08) node[below] {\small $n{-}2K\Kz$};
  \draw (\Wm,0.08) -- (\Wm,-0.08) node[below] {\small $n{+}2K\Kz$};
  \draw (\Wm/2,0.08) -- (\Wm/2,-0.08) node[below] {\small $n$};
  \draw (\jj,0.08) -- (\jj,-0.08) node[below] {\small $n{+}j$};

  \node[blue!55!black] at (2,0.7) {$B_{\gamma,n}(K)$};
  \node[green!55!black,anchor=south] at (\Wm/2,\Uflat+0.18) {$\bar B_{\gamma,n}(K)$};
\end{scope}
\end{tikzpicture}
  \caption{\textbf{The local window and its partition} (drawn for the case $\gamma>0$).
  The finite window $\bar B_{\gamma,n}(K)$ (green frame), of horizontal width $4KK_0$
  centred at a column $n$ and height $K+\gamma\log n+K_0$, and the strip
  $B_{\gamma,n}(K)$ it contains (blue). At a column $n+j$, the partition
  $\{\mathbf Z^n_{(j,k)}\}$ of $\bar B_{\gamma,n}(K)$ consists of singletons for $k<K$
  together with the set $\mathbf Z^n_{(j,K)}$ (spanning $K\le k\le K+\gamma\log n$,
  about $\gamma\log n$ states, hence growing with $n$) and the $K_0$-cap
  $\mathbf Z^n_{(j,K+1)}$ of fixed height just above it.}
  \label{fig:window-partition}
\end{figure}

The window $\bar B_{\gamma,n}(K)$, the strip $B_{\gamma,n}(K)$ it contains, and the partition $\{\XX_m\}$ are illustrated in \fref{fig:window-partition}.

Finally, for $x \in \ss(K)$, $W \subseteq \bar B_{\gamma,n}(K)$ and $A \subseteq W$ we define
\begin{equ}
  Q_W(x,A) := \mathbb P(Y(\tau_W) \in A~|~Y(0) = x )
\end{equ} where $\tau_W := \inf\{s>0~:~Y(s) \in W\}$. The rest of the paragraph is dedicated to establishing bounds on the above quantity when $A$ is one of the $\XX_{(j,k)}$.

\begin{remark}\label{r:singletons} When $k < K$ (and also when $k=K$ for $\gamma = 0$) the sets defined above are singletons. Moreover, since $k<K$ implies $(j,k)\in B_0(K)$, the following estimate holds: whenever there exists an active reaction $y \to y' \in \mathcal R$ such that $(j',k') - (j,k) = y'-y$, by \eref{e:rateslb} and \eref{e:ratesub} when $\gamma = 0$ (and the lines underneath that equation) we have
  \begin{equ}
    a^{-1} q_{(j,k),(j',k')}(n) \leq Q_{\{(j,k)_n\}^c}((j,k)_n,\{(j',k')_n\}) \leq a q_{(j,k),(j',k')}(n)
  \end{equ}
  where $q_{(j,k),(j',k')}$ is a power of $n$ and $a > 0$ is a sufficiently large, $n$-independent constant.
  \end{remark}

  \subsubsection*{$(1,0)$-WE networks} The case of (1,0)-weakly endotactic network consists of those networks whose dominant reactions in the $(1,0)$ direction are \emph{all} orthogonal to the vector $w = (1,0)$. This class of networks corresponds to systems that do not dissipate their $x_1$ component, at least to first order, as long as $x_2$ is large enough for the dominant reactions to be active. In the analysis of stability and dissipation for $(1,0)$-WE networks the above estimate is not sufficient. Indeed, in that case the most likely exit from $B_{0,n}$ for a point $(j,k)_n$ with sufficiently large $K$ is going to be in $\XX_{(j,K+1)}$, and our coarse-graining analysis will not allow to observe dissipation of the trajectories. To observe such dissipation we must instead consider the dynamics of $Y$ in a larger set, where the probability of an exit through $\XX_{(j,K+1)}$ is negligible. The following lemma, proven in the appendix, shows that this is indeed the case for the set $B_{\gamma,n}$ for $K>0$ and $\gamma>0$ sufficiently large.

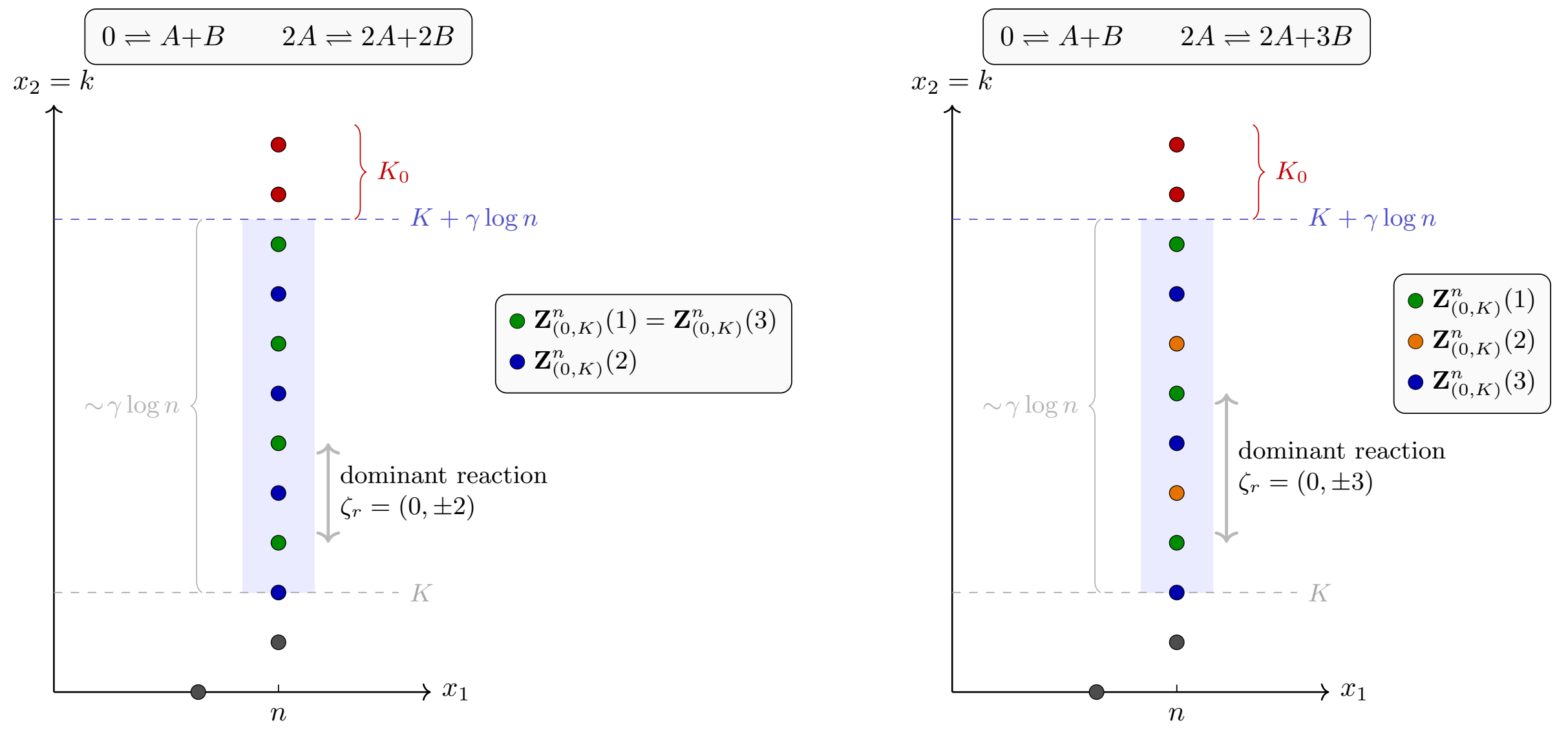
\begin{figure}[t]
  \centering
  \colorlet{SliceA}{blue!70!black}
  \colorlet{SliceB}{green!55!black}
  \colorlet{SliceC}{orange!90!black}
  \colorlet{CapCol}{red!75!black}
  \colorlet{StripCol}{blue!8}
  \newcommand{\swatch}[1]{\tikz[baseline=-0.6ex]{\filldraw[fill=#1,draw=black,line width=0.2pt] (0,0) circle (2.6pt);}}
  \resizebox{\textwidth}{!}{%
\begin{tikzpicture}[
    yscale=0.62, xscale=1.0,
    every node/.style={font=\small},
    dot/.style={circle,inner sep=0pt,minimum size=5.2pt,draw=black,line width=0.3pt},
  ]

\begin{scope}[xshift=0cm]

  \node[draw,rounded corners,fill=gray!5,inner sep=6pt,anchor=south]
        at (0,12.6)
        {$0 \rightleftharpoons A{+}B \qquad 2A \rightleftharpoons 2A{+}2B$};

  \fill[StripCol] (-0.45,2) rectangle (0.45,9.5);

  \draw[->,line width=0.6pt] (-2.8,0) -- (1.9,0) node[right]{$x_1$};
  \draw[->,line width=0.6pt] (-2.8,0) -- (-2.8,11.8) node[above]{$x_2=k$};
  \draw (0,0) node[below=2pt]{$n$} -- (0,0.15);

  \draw[dashed,gray!70] (-2.8,2) -- (1.5,2) node[right,font=\footnotesize]{$K$};
  \draw[dashed,SliceA!70] (-2.8,9.5) -- (1.5,9.5)
        node[right,font=\footnotesize,text width=1.7cm,align=left]{$K+\gamma\log n$};

  \foreach \k in {2,4,6,8}{ \node[dot,fill=SliceA] at (0,\k){}; } 
  \foreach \k in {3,5,7,9}{ \node[dot,fill=SliceB] at (0,\k){}; } 
  \foreach \k in {10,11}{ \node[dot,fill=CapCol] at (0,\k){}; }
  \node[dot,fill=gray!60!black] at (0,1){};    
  \node[dot,fill=gray!60!black] at (-1,0){};   

  \draw[<->,line width=1.1pt,gray!55]
       (0.62,3) -- (0.62,5)
       node[midway,right,font=\footnotesize,text=black,align=left]
       {dominant reaction\\ $\zeta_r=(0,\pm 2)$};

  \draw[decorate,decoration={brace,amplitude=4pt},gray!60]
       (-0.95,2) -- (-0.95,9.5)
       node[midway,left=4pt,font=\footnotesize]{$\sim\!\gamma\log n$};
  \draw[decorate,decoration={brace,amplitude=4pt,mirror},CapCol]
       (0.95,9.5) -- (0.95,11.4)
       node[midway,right=4pt,font=\footnotesize,text=CapCol]{$\Kz$};

  \node[draw,rounded corners,fill=gray!4,inner sep=5pt,anchor=west,
        align=left,font=\footnotesize] at (2.7,7) {%
     \swatch{SliceB}~$\XX_{(0,K)}(1)=\XX_{(0,K)}(3)$\\[3pt]
     \swatch{SliceA}~$\XX_{(0,K)}(2)$%
  };

  \node[anchor=north,text width=7cm,align=center,font=\small] at (0.1,-0.9)
       {\textbf{(a)} the dominant reactions are vertical, with reaction vectors $\zeta_r=(0,\pm 2)$, so $\XX_{(0,K)}$ splits into $2$ slices (even/odd heights).};

\end{scope}

\begin{scope}[xshift=11.2cm]

  \node[draw,rounded corners,fill=gray!5,inner sep=6pt,anchor=south]
        at (0,12.6)
        {$0 \rightleftharpoons A{+}B \qquad 2A \rightleftharpoons 2A{+}3B$};

  \fill[StripCol] (-0.45,2) rectangle (0.45,9.5);

  \draw[->,line width=0.6pt] (-2.8,0) -- (1.9,0) node[right]{$x_1$};
  \draw[->,line width=0.6pt] (-2.8,0) -- (-2.8,11.8) node[above]{$x_2=k$};
  \draw (0,0) node[below=2pt]{$n$} -- (0,0.15);

  \draw[dashed,gray!70] (-2.8,2) -- (1.5,2) node[right,font=\footnotesize]{$K$};
  \draw[dashed,SliceA!70] (-2.8,9.5) -- (1.5,9.5)
        node[right,font=\footnotesize,text width=1.7cm,align=left]{$K+\gamma\log n$};

  \foreach \k in {2,5,8}{ \node[dot,fill=SliceA] at (0,\k){}; } 
  \foreach \k in {3,6,9}{ \node[dot,fill=SliceB] at (0,\k){}; } 
  \foreach \k in {4,7}{   \node[dot,fill=SliceC] at (0,\k){}; } 
  \foreach \k in {10,11}{ \node[dot,fill=CapCol] at (0,\k){}; }
  \node[dot,fill=gray!60!black] at (0,1){};    
  \node[dot,fill=gray!60!black] at (-1,0){};   

  \draw[<->,line width=1.1pt,gray!55]
       (0.62,3) -- (0.62,6)
       node[midway,right,font=\footnotesize,text=black,align=left]
       {dominant reaction\\ $\zeta_r=(0,\pm 3)$};

  \draw[decorate,decoration={brace,amplitude=4pt},gray!60]
       (-0.95,2) -- (-0.95,9.5)
       node[midway,left=4pt,font=\footnotesize]{$\sim\!\gamma\log n$};
  \draw[decorate,decoration={brace,amplitude=4pt,mirror},CapCol]
       (0.95,9.5) -- (0.95,11.4)
       node[midway,right=4pt,font=\footnotesize,text=CapCol]{$\Kz$};

  \node[draw,rounded corners,fill=gray!4,inner sep=5pt,anchor=west,
        align=left,font=\footnotesize] at (2.7,7) {%
     \swatch{SliceB}~$\XX_{(0,K)}(1)$\\[3pt]
     \swatch{SliceC}~$\XX_{(0,K)}(2)$\\[3pt]
     \swatch{SliceA}~$\XX_{(0,K)}(3)$%
  };

  \node[anchor=north,text width=7cm,align=center,font=\small] at (0.1,-0.9)
       {\textbf{(b)} the dominant reactions are vertical, with reaction vectors $\zeta_r=(0,\pm 3)$, so $\XX_{(0,K)}$ splits into $3$ slices (heights mod $3$).};

\end{scope}

\end{tikzpicture}%
  }
  \caption{\textbf{The slices $\XX_{(0,K)}(\ell)$ of $\XX_{(0,K)}$, for two $(1,0)$-WE networks} (drawn for $\gamma>0$, not to scale). In each panel the set $\XX_{(0,K)}$ consists of the lattice points on the column $x_1=n$ with $K\le k\le K+\gamma\log n$ (about $\gamma\log n$ of them), marked by the coloured dots within the shaded band. Since the network is $(1,0)$-WE, the dominant reactions at these points are vertical, with reaction vector $\zeta_r$. Following them connects points whose heights differ by a multiple of $\|\zeta_r\|$, so $\XX_{(0,K)}$ partitions into $\|\zeta_r\|$ slices $\XX_{(0,K)}(\ell)$ (the colours), one per residue class of $k$---distinct $\ell$ may give the same slice (e.g.\ $\XX_{(0,K)}(1)=\XX_{(0,K)}(3)$ in (a)). In (a) the dominant reaction $2A\rightleftharpoons 2A+2B$ has $\zeta_r=(0,\pm2)$, giving two slices (even/odd $k$). In (b) $2A\rightleftharpoons 2A+3B$ has $\zeta_r=(0,\pm3)$, giving three ($k\bmod3$). The red points just above the band are the $\Kz$-cap $\XX_{(0,K+1)}$. The grey points below are those connected to $\XX_{(0,K)}$ by occurrences of dominant reactions: at $x_2=0$ and $x_2=1$ the dominant reactions are $0\to A+B$ and $A+B\to0$, respectively.}
  \label{fig:zell}
\end{figure}

  For $\gamma>0$, let
  \[
    d=\gcd\{|(y'-y)_2|:y\to y'\in\cal R_{(1,0)},\ (y'-y)_2\ne0\}.
  \]
  The degree-maximal subsystem is vertical and weakly reversible, so its nonzero increments have both signs. For $1\le\ell\le\Kz$, define
  \begin{equ}
    \XX_{(j,K)}(\ell):=\{(j,k)_n\in\XX_{(j,K)}:k\equiv K+\ell\pmod d\}.
  \end{equ}
  Only the distinct representatives $1\le\ell\le d$ are used in the partition. Larger indices give the same sets, as in \fref{fig:zell}. Above sufficiently large $K$, any two fixed heights in the same residue are joined by a dominant path staying above $K$. Indeed, the signed increments generate $d\mathbb Z$. Adding nonnegative zero combinations makes all coefficients nonnegative, and ordering positive increments before negative ones keeps the path above its lower endpoint. For finitely many fixed endpoints these paths have a common finite height and length bound.

  \begin{lemma}[Bulk exit bounds]\label{l:noexitfromtop}
    Assume the network is $(1,0)$-WE. Let $\emptyset\ne L\subseteq\{1,\dots,d\}$ and $W:=(\bigcup_{\ell\in L}\XX_{(j,K)}(\ell))^c$. For every sufficiently large $K>\Kz$ and every $\beta>0$, take $\gamma=\beta$. Then there exists $C>0$ such that, for all sufficiently large $n$, uniformly over $j\in[-2K\Kz,2K\Kz]$,
    \begin{equs}
      \sup_{\substack{K\le k\le K+\Kz\\(j,k)_n\in W^c}}Q_{W}((j,k)_n,\XX_{(j,K+1)}) &\leq C n^{-\beta}\,,\\
      \sup_{(j,k)_n \in \bigcup_{\ell \in L} \XX_{(j,K)}(\ell)}Q_{W}((j,k)_n,\XX_{(j,K+1)}) &\leq C (\log n)^{-1}
    \end{equs}
    For fixed $R$ and $H=H_n+O(1)$, the same bounds hold for exit above $H$ from any $W^c\subseteq\{|x_1-n|\le R,\ K\le x_2\le H\}$, and $\sup_{x\in W^c}\E_x\tau_W\le C\log n$.

  \end{lemma}
  In particular, letting $L = \{1, \dots, d\}$ the above result says that the probability of exiting the vertical strip from the top when entering from its bottom is (algebraically) as small as needed, and it is going to $0$ uniformly in the starting point within the strip.

    The following result states that exit probabilities and the most likely exit sets are, to leading order in $n$, shared by all points 
    that can be connected by dominant reactions. Note that this result does not hold for exits into the bulk (i.e., through $\XX_{(j,K+1)}$, covered by \lref{l:noexitfromtop} above). The proof of the lemma is given in the appendix to improve readability.

\begin{lemma}[Boundary exit bounds]\label{l:commonexit}
      Let $K$ be sufficiently large, $\gamma> 0$, $j\in [-2K\Kz, 2K\Kz]$ and fix $\ell \in \{1, \dots, d\}$. Let $\widetilde {\mathcal R}'$ be a subset of $\mathcal R$ and let $A_n$ be one of the following:
      \begin{enumerate}
      \item $A_n=\{x+y'-y \in B_{\gamma,n}(K)~:~(x+y'-y)_2\geq K, x \in \XX_{(j,K)}(\ell), y\to y'\in \widetilde {\mathcal R}'\}$.
      \item $A_n\subseteq\{x+y'-y \in B_{\gamma,n}(K)~:~(x+y'-y)_2< K, x \in \XX_{(j,K)}(\ell), y\to y'\in \widetilde {\mathcal R}'\}$.
      \end{enumerate}
      Let now $\widetilde {\mathcal R}$ be the set of reactions that allow to reach $A_n\setminus \XX_{(j,K)}(\ell)$ from $\XX_{(j,K)}(\ell)$, i.e.\
      \begin{equation*}
          \widetilde {\mathcal R}=\{y\to y'\in \mathcal R~:~\exists\, x\in \XX_{(j,K)}(\ell) \text{ with } x+y'-y\in A_n\setminus \XX_{(j,K)}(\ell)\}.
      \end{equation*}
      Then, if $\widetilde {\mathcal R}$ is not empty, there exists $C,C'>0$ such that, for all sufficiently large $n$, uniformly on $(j,k)_n\in\XX_{(j,K)}(\ell)$
      \begin{equ}
      (C \log n)^{-C'} n^{-\tilde \beta} \leq Q_W((j,k)_n , A_n) \leq (C \log n)^{C'} n^{-\tilde \beta}\,,
      \end{equ}
      where $W := \XX_{(j,K)}(\ell)^c$ and $\tilde \beta := \max_{y \to y' \in \mathcal R} y_1 - \max_{\tilde y \to \tilde y' \in \widetilde{\mathcal R}} \tilde y_1$\,.
\end{lemma}

\subsection{W-graphs and exit probabilities}\label{s:wgraphs1}
We now proceed to introduce the necessary quantities from \cite{fw98} to realize the program outlined at the beginning of the section.
The first quantity we introduce is the one of $W$-graphs. Roughly speaking, for a finite set $L$ and a ``boundary'' $W \subset L$, $W$-graphs are transition graphs that represent, for every pair of points $m \in L\setminus W$ and $m' \in W$ a unique exit path  within $L\setminus W$ from $m$ to $m'$, i.e. a path in $L\setminus W$ leading from $m$ to $m'$. Throughout, we define a graph $g$ as a pair of vertices $V(g)$ and (directed) edges $E(g) \subseteq V(g)\times V(g)$.
\begin{definition}[W-graph] \label{d:wgraph} Let $L$ be a finite set of points, $W \subset L$. Then a graph $g$  with $V(g) = L$ is a $W$-graph if all of its edges start from $L\setminus W$ and 
\begin{enumerate}
\item every $m \in L \setminus W$ is the origin of exactly one directed edge,
\item there are no closed loops in $g$.
\end{enumerate}
A condition equivalent to 2) is: for every $m \in L\setminus W$ there exists a path $m \to \dots m' \in g$ for an $m' \in W$.
\end{definition}
We will denote by $G(W)$ the family of $W$-graphs and for any pair $m \in L\setminus W , m' \in L $ by $G_{m,m'}(W)\subset G(W)$ the subset of $W$-graphs connecting $m$ to $m'$. Furthermore, for a given graph $g$, associating to each directed edge $m \to m' \in E(g)$ a weight $q_{m,m'}$ we define $\pi(g)$ as
\begin{equ}\label{e:pi}
  \pi(g)  := \prod_{m \to m' \in E(g)} q_{m, m'}\,.
\end{equ}

The following lemma (analogous to \cite[Lemma 3.3, Chapter 6]{fw98}) estimates exit probabilities from the union of sets $\{\OO_m\}_{m \in L \setminus W}$ through the cost associated to $W$-graphs, given some uniform (in the entrance point) control on the transition probabilities between those sets.

The sets $\OO_m$ are finite, and the stopping sets contain the complete one-jump exterior of each transient union. Write $\OO_U=\bigcup_{m\in U}\OO_m$ and $\sigma_m=\inf\{s>0:Y(s)\notin\OO_m\}$. We assume $\mathbb P_x(\sigma_m<\infty)=1$ for $m\notin W$, $x\in\OO_m$.
\begin{lemma}\label{l:exitestimate1}
Let $L$ be finite, $W\subset L$, and $\OO = \bigcup_{m \in L} \OO_m$ with $\OO_m \cap \OO_{m'} = \emptyset $ if $m \neq m'$ and assume that there exist $a \geq 2$ and 
$\{q_{m,m'}\}_{m,m' \in L^2}$ a family of non negative numbers such that for every $m\in L\setminus W$, $m'\in L\setminus\{m\}$ we have
\begin{equ}\label{e:exitestimate1}
  a^{-1} q_{m,m'} \leq  Q_{\OO \setminus\OO_m}(x, \OO_{m'}) \leq a q_{m,m'}\qquad \text{for all } x \in \OO_m\,,
\end{equ}
 Then, if $r := |L\setminus W|$,
\begin{equ}\label{e:exitestimate2}
a^{-4^r} \frac{\sum_{g \in G_{m,m'}(W)} \pi(g)  }{\sum_{g \in G(W)} \pi(g)  }  \leq Q_{\OO_{W}}(x, \OO_{m'}) \leq a^{4^r} \frac{\sum_{g \in G_{m,m'}(W)} \pi(g)  }{\sum_{g \in G(W)} \pi(g)  }\,.
\end{equ}
for all $x\in \OO_m$, $m \notin W$, $m' \in W$, where $\OO_W=\cup_{m''\in W} \OO_{m''}$, and provided that the denominator is positive.

For the second assertion, let $\bar m\in W$, and require \eref{e:exitestimate1} only for $m'\ne\bar m$, together with
\begin{equ}\label{e:exitestimate3}
  Q_{\OO_{L\setminus U}}(x,\OO_{\bar m})\le a^{-1}
  \qquad\text{for all }x\in\OO_m,\quad
  \{m\}\subseteq U\subseteq L\setminus W.
\end{equ}
Then \eref{e:exitestimate2} holds for $m'\in W\setminus\{\bar m\}$, with $G(W),G_{m,m'}(W)$ denoting $(W\setminus\{\bar m\})$-graphs on $L\setminus\{\bar m\}$, again provided that the denominator is positive.
\end{lemma}

The proof of the lemma above mimics the one of \cite[Lemma 3.3, Chapter 6]{fw98} and is deferred to the appendix. The same conclusion holds with $1/2$ in place of $a^{-1}$ in \eref{e:exitestimate3}. In WE, \lref{l:exittimes} supplies this bound. For the WE cap, apply the second assertion with $\bar m$ its index.
A possible choice of $\{\OO_m\}$ is given, e.g., by the partition $\{\XX_m\}$ defined in the previous section. 
We introduce the corresponding sets of indices $L_\gamma$ as follow: 
 \begin{description}
    \item[if $\gamma = 0$]\,\\ 
    $L_0:= \{(j, k, 0)\in \Z^3 : j\in [-2K\Kz, 2K\Kz], k\in [0, K]\}$
    \item[if $\gamma > 0$] \,\\
        $L_\gamma:= \{ (j,k, 0)\in \Z^3 : j\in [-2K\Kz, 2K\Kz], k\in [0, K-1]\}$\\ 
        $\cup \{(j,K,\ell)\in\Z^3:j\in[-2K\Kz,2K\Kz],\ \ell\in\{1,\ldots,d\}\}$.
  \end{description}

In the $(1,0)$-SE case, using the estimates \eref{e:ratesub} and \eref{e:rateslb}, we defined in \dref{d:cost} a cost for each edge (and each path) in the graph representing the transitions between the states in $B_{0, n}(K)$. This graph is independent of $n$, and can be seen as a graph with vertices in $L_0$. In the $(1,0)$-WE case, an equivalent graph can be defined on (the finite set) $L_\gamma$ using the result of \lref{l:commonexit}. We extend the definition of $I$ to this graph. Combining these results to \lref{l:exitestimate1}, we can show the following corollary. 

Include all possible one-jump exits in the stopping set. In SE we record which state the process reaches when it jumps above the strip. Each such state has its own index in $W$, and $\XX_\dagger=\emptyset$. In WE all jumps above the strip count as one exit event. Their landing states form $\XX_\dagger$, with $\gamma$ chosen as in \lref{l:exittimes}. Give the remaining exit sets indices in $W$ and enlarge $L_\gamma$ to contain $W$. We use $\inf\emptyset=\infty$ and $\log0=-\infty$.
\begin{corollary}\label{c:fw}
    Let $W\subset L_\gamma$, $m\in L_\gamma\setminus W$ and $m'\in W$, and set $\XX_W:= \cup_{m''\in W}\XX_{m''}$.     
    If $\inf_{g\in G(W)}I(g)<\infty$, then uniformly for $x\in\XX_m$, if the Markov chain $Y(s)$ starts at $x$, 
  \begin{equ}
    \lim_{n \to \infty} (\log n)^{-1} \log Q_{\XX_W\cup\XX_\dagger}(x,\XX_{m'}) = -\mathcal W(m,m')\,,
  \end{equ}
  where 
  \begin{equ}
    \mathcal W(m,m') := \inf_{g \in G_{m,m'}(W)} I(g) - \inf_{g' \in G(W)} I(g')
  \end{equ}
\end{corollary}

In other words, the probability that $Y(s)$ started in $\XX_m \in \ss(K)$ will  first exit $\ss(K) \setminus \XX_W$ at $\XX_{m'} \in \XX_W$ as $n \to \infty$ has logarithmic rate given by minus the minimal (total) large deviations cost of $W$-graphs connecting $m$ to $m'$ plus the minimal (total) large deviations cost of all possible $W$-graphs. In this sense, $W$-graphs capture the globally cheapest paths to leave from every point of $L_\gamma\setminus W$ into $W$.
\begin{proof}[Proof of \cref{c:fw}]
  It is sufficient to identify $q_{m,m'}$ with monomials in $n$ and $a$'s with multiplicative logarithmic prefactors from \eref{e:ratesub}, \eref{e:rateslb} (and in the (1,0)-WE case \lref{l:commonexit}). Then, the result follows by \lref{l:exitestimate1} upon taking a log and the $n \to \infty$ limit.
\end{proof}

\subsection{Combining $W$-graphs}\label{s:wgraphs2}

We now proceed to construct the family of graphs $G$ with vertices in $\ss(K)$ in which we will look for $W$-graphs. The main idea here is to identify a ``vertical'' set of points in $\ss(K)$ that can be used as level lines for a Lyapunov function as the one described in Section 2. More specifically, we want to identify a set of points in $\ss(K)$ with one element $m_i \in \ss(K)$ for each value $i \in \{0, \dots, K\}$ of the second coordinate, i.e., $m_i = (k_i, i)_n$. These sets will induce a total order in $\sg$.

Until the end of this section, we consider a weighted, directed, finite-degree graph $\G$ with weights $\{J_{x\to x'}\}_{x \to x' \in E(\G)}$. We define the total cost of such a graph by  
\[
    I(\G):= \sum_{x\to x'\in E(\G)} J_{x\to x'}.
\] 

The following lemma shows that the coarse-graining procedure preserves the uniformity of the exit probability estimates in the starting condition $x$: if the exit probabilities do not depend, to leading order, on the initial condition for a family of sets $\{\OO_1, \dots, \OO_M\}$ and these sets are strongly connected by their dominant transitions, then the exit probabilities from the \emph{union} of these sets does not depend on the initial condition, either.
To interpret the statement of the lemma in this sense, recall that the exit probabilities are related to the cost of $W$-graph by~\cref{c:fw}, choosing as $J_{v,v'}$ the cost defined in~\dref{d:cost}. 

To state the lemma we introduce a few definitions. Given a set $\OO\subset V(\G)$, $x\in \OO^c$ is said to be \emph{reachable from} $\OO$ if for some $v\in \OO$, $v\to x\in E(\G)$. 
For every $\OO\subset V(\G)$, the minimum exit cost over all starting points in $\OO$ to a point $w\in \OO^c$ reachable from $\OO$ is:
\begin{equ}
 I_w^*(\OO) := \min_{v \in \OO}\min_{g \in G_{v,w}(V(\G)\setminus \OO)}I(g)\,.
\end{equ}
The minimal-cost exit points of the set $\OO$ are defined as 
\begin{equ}\label{e:nearestNeighbor}
	\mathcal{E}(\OO) := \arg\min\{I_w^*(\OO): w\in V(\G) \setminus \OO \text{  reachable from  } \OO\}\,.  
\end{equ}
For two disjoints sets $\OO, \OO^*$ of $V(\G)$, we say that $\OO$ is connected with minimal-cost to $\OO'$, noted $\OO\to \OO'$, if $\mathcal{E}(\OO)\cap \OO'\neq \emptyset$. 

\begin{lemma}\label{l:coarsening}
  Let $\OO_1, \dots \OO_M$ be disjoint sets of $V(\G)$ and assume that, for all $i \in \{1, \dots, M\}$ and for all $y \in\OO_i^c$ reachable from $\OO_i$, we can find a graph
	\[
		g\in \bigcap_{x\in \OO_i} G_{x,y}(V(\G)\setminus \OO_i)
	\]
	with cost $I(g)=I_y^*(\OO_i)$.
 
	Assume further that the cycle-graph of the minimal-cost connections among $\{\OO_1, \dots, \OO_M\}$ is strongly connected. 
	Then for all $y \in V(\G)\setminus\bigcup_i \OO_i$ reachable from $\bigcup_i \OO_i$, we can find a $V(\G)\setminus \bigcup_i\OO_i$-graph $g'$ with cost $I(g')=I_y^*(\bigcup_i \OO_i)$ and such that $$g'\in \bigcap_{x\in \bigcup_i \OO_i} G_{x,y}\left(V(\G)\setminus\bigcup_i \OO_i\right)$$. 
\end{lemma}

We summarize in the following lemma some properties of partitions of $W$-graphs, which will be useful in proving Lemma~\ref{l:coarsening}. When restricting a $W$-graph $g$ to a subset $\OO_i \subset V(\G)\setminus W$ we define  $E(g_i):=\{v \to v' \in E(g)~:~v \in \OO_i\}$ and $V(g_i):=V(\G)$. 
\begin{lemma}\label{l:unionwgraphs}
  Let  $\{\OO_1, \dots, \OO_M\}$ be a partition of $V(\G) \setminus W$. Each subgraph $g_i$ of a $W$-graph $g$ is itself a $V(\G) \setminus \OO_i$-graph. Inversely, if each $g_i$ is a $V(\G) \setminus \OO_i$-graph and $g = \bigcup_i g_i$ does not contain any directed cycle then $g$ is a $W$-graph. In both cases we have that $I(g)   = \sum_i I(g_i)$.
\end{lemma}

\begin{proof}[Proof of \lref{l:unionwgraphs}]
  We prove the second statement, as the first and third one are direct consequences of the definition of $W$-graph. From \dref{d:wgraph} (1) we know that each $v \in  \OO_i$ is the origin of exactly one directed edge. Since the union of a $V(\G)\setminus \OO_i$-graph to another graph $g'$ will only add incoming edges to $g'$,  \dref{d:wgraph} (1) also holds in $g$ for all  $v \in \bigcup_i \OO_i$. (2) follows from the assumption of the Lemma.
\end{proof}

\begin{proof}[Proof of \lref{l:coarsening}]	
	Let $y\in V(\G)\setminus\bigcup_i \OO_i$ reachable from $\bigcup_i\OO_i$, $g$ an $V(G){\setminus}\bigcup_i \OO_i$-graph, $y$-minimal and let $i(y)\in \{1, \ldots, M\}$ so that a parent $x$ of $y$ is in $\OO_{i(y)}$.
	We construct a $V(\G){\setminus}\bigcup_i\OO_i$-graph $g'$, that satisfies the conditions of the theorem, as follow: choose a directed spanning tree toward $\OO_{i(y)}$ in the cycle-graph. For every $j \in \{1,\dots, M\}\setminus \{i(y)\}$, let $$\OO_j \to \OO_{j'} \to \dots \to \OO_{j^{(n)}} =  \OO_{i(y)}$$ be the path in this tree and denote by $J$ the ordered set  of indices on this path. 
	Now let $y_j\in \mathcal{E}(\OO_j)\cap\OO_{j'}$ an exit point of $\OO_j$ with lowest cost, with $\mathcal{E}$ defined in Relation~\eqref{e:nearestNeighbor}. We introduce a $g'_j \in \bigcap_{x \in \OO_j}G_{x,y_j}(V(\G) \setminus \OO_j)$, $y_j$-minimal, and concatenate it with  $g'_{j'} \in \bigcap_{x \in \OO_{j'}} G_{x,y_{j'}}(V(\G) \setminus \OO_{j'})$, $y_{j'}$-minimal with $y_{j'}\in \mathcal E(\OO_{j'})\cap \OO_{j^{(2)}}$, and so on until a minimal $g'_{i(y)}=g'_{j^{(n)}} \in \bigcap_{x \in \OO_{j^{(n)}}} G_{x,y}(V(\G) \setminus \OO_{j^{(n)}})$, which all exist by assumption. Every path in $g^{\prime}=\bigcup_j g_j^{\prime}$ follows the chosen tree and then reaches $y$, so $g^{\prime}$ is a $V(\G)\setminus\bigcup_i\OO_i$-graph. 

	For $j\in \{1, \ldots, M\}\setminus\{i(y)\}$, $g_j'$ has minimal exit cost of $\OO_j$, and $g_j$ being a $V(\G)\setminus \OO_j$-graph by~\lref{l:unionwgraphs}, we have $I(g_j)\geq I(g_j')$. Besides, $g_{i(y)}\in G_{x,y}(V(\G)\setminus\OO_{i(y)})$, and therefore $I(g_{i(y)})\geq I_y^*(\OO_{i(y)})=I(g_{i(y)}')$. 
	This leads to 
	\begin{equ}
		I(g')   = \sum_{i =1}^M I(g_i') \leq  \sum_{i =1}^M I(g_i) = I(g)=I^*_y\left(\bigcup_{i}\OO_i \right)
	\end{equ}
	and by construction, $g'\in G_{x,y}(V(\G)\setminus\bigcup_i \OO_i)$ for all $x\in \bigcup_{i}\OO_i$. 
\end{proof}

\section{Boundary layer: designing a local Lyapunov function}
\label{s:boundary2}

\subsection{Application to weakly reversible CRNs}\label{s:fibers1}

While the lumping procedure introduced in the previous section could be applied to any family of Markov chains in a finite set whose rates satisfy \eref{e:ratesub} and \eref{e:rateslb}, we will now study the properties of this procedure resulting from the weakly reversible structure of a CRN. In particular, following the intuition developed in \exref{ex:mainexample}, we will show that this procedure can only result in a sequence of most likely transitions between lumped states that moves the process towards the origin and never away from it, thereby showing that the dynamics of this family of networks is \emph{dissipative} close to the boundary. More precisely, we will show that there can never be an iteration of the lumping process that includes a sequence of $\mathcal O(1)$ transitions between a lumped set to another that is directly to its right (read: to another ``copy'' of the same set translated away from the origin). We state this result in the following key lemma, connecting the structure of weakly reversible CRNs to the dissipative nature of the associates mass action kinetics close to the boundary.

The initial sets $\OO_m^{(0)}$ are the translates by $(j,0)$, $j\in\mathbb Z$, of the sets in column $0$ listed below.
\[
\begin{array}{ll}
 (1,0)\text{-SE}: & \{(0,k)\},\quad 0\le k\le K,\\[2pt]
 (1,0)\text{-WE}: & \{(0,k)\},\quad 0\le k<K_c,\\
 & \{(0,k): k\ge K_c,\ k\equiv K_c+\ell\pmod d\},\quad 1\le\ell\le d .
\end{array}
\]
In WE fix $K_c$ sufficiently large and take $K\ge K_c+2\Kz$. In SE retain the singleton exit targets with $K<k\le K+\Kz$. Write $s$ for the number of sets in column $0$, including these targets.

\begin{lemma}\label{l:noforward}
  Let $\tilde {\mathcal G}(m)$ be the graph resulting from the $m$-step coarse-graining procedure from last section started on $\mathcal G_0$. Then if there exists $\ell > 0$ such that $\OO' = \OO + (\ell,0)$, we have $\OO \not \Rightarrow_{\tilde {\mathcal G}(m)} \OO'$.
\end{lemma}

To make the connection between the structure of the CRN and its lumped dynamics through the proof of \lref{l:noforward}, we remind the reader about the following key theorem, resulting from combining \cite{horn72} with \cite{anderson10}.
\begin{theorem}\label{t:complexbalance}
  Let $(\mathcal S, \mathcal C, \mathcal R)$ be a weakly reversible CRN. There exists a choice of reaction rate constants $\{\kappa_r\}_{r \in \R}$ and a $c \in \mathbb R_{>0}^d$ such that the associated process $X$ has a stationary distribution $\mu$ on $\mathbb Z_{\geq 0}^d$ with mass function
  \begin{equ}
    \mu(x) = \prod_{i = 1}^d e^{-c_i} \frac{c_i^{x_i}}{x_i!}\qquad \qquad x \in \mathbb Z_{\geq 0}^d.
  \end{equ}
\end{theorem}

This theorem will be pivotal in our proof of dissipation. Indeed, our proof of \lref{l:noforward} will work by contradiction: since the estimates we establish are independent on the choice of reaction rate constants $\{\kappa_r\}$, were the process not dissipative for a choice of $\{\kappa_r\}$ then the same would hold for the choice in \tref{t:complexbalance}, which in turn would contradict the superexponential decay of Poisson tails. The proof is given below.

\begin{proof}[Proof of \lref{l:noforward}]
  Assume that there exists $\ell\in \mathbb N\setminus\{0\}$ such that $\OO' = (\ell,0) + \OO$ and $\OO \Rightarrow_{\mathcal G(m)} \OO'$. Since the transition costs are independent of the positive rate constants, choose the constants in \tref{t:complexbalance} with $c=(1,1)$. For a state $x$, a set $D$ and a lumped set $\OO$, put
\[
\begin{aligned}
 a(x,D)&:=\sum_{y\in D\setminus\{x\}}q(x,y)
        =\sum_{r\in\Re:\,x+\zeta_r\in D\setminus\{x\}}\lambda_r(x),\\
 \rho(\OO)&:=\sum_{z\in\OO}\mu(z)a(z,\OO^c),\\
 \rho(x)&:=\rho(\{x\})=\mu(x)\sum_{y\ne x}q(x,y).
\end{aligned}
\]
In stationarity, visits to $\OO$ start and end at rate $\rho(\OO)$.
Count the visits whose exit is into $\OO'$.
Their rate is the middle term in
\begin{equ}\label{e:cellflux}
 \rho(\OO)\inf_{x\in\OO}\mathbb P_x(Y(\tau_{\OO^c})\in\OO')
 \leq\sum_{x\in\OO}\mu(x)a(x,\OO')\leq\min\{\rho(\OO),\rho(\OO')\}.
\end{equ}
Here the last bound uses balance for $\OO'$.

For a reaction from $(n+j,k)$, the Poisson weights give
\[
 \mu(n+j,k)\lambda_r(n+j,k)
 =\frac{e^{-2}\kappa_r}{(n+j-y_{r,1})!\,(k-y_{r,2})!}.
\]
After removal of $1/n!$, the exponent of this reaction is $y_{r,1}-j$.
Let $e(\OO)$ be the largest such exponent among the exits from the untranslated set $\OO$.
In SE include cap exits. In WE use genuine exits from the untruncated set.
The vertical factorial sums converge, and the flux at the logarithmic cap in WE is $o(n^{-A}/n!)$ for every fixed $A$.
Thus $\rho(\OO+(n,0))\sim b_{\OO}n^{e(\OO)}/n!$, with $b_{\OO}>0$ for each set with a genuine exit.
A set of bounded width without a genuine exit cannot meet a component unbounded in the first coordinate, by weak reversibility.

In WE we apply \eref{e:cellflux} to the translates of $\OO,\OO'$ by $(n,0)$ restricted to $x_2\leq\lfloor K_c+\gamma\log n\rfloor$, where \cref{c:fw} applies with $K=K_c$ and the cap flux estimated above leaves the exponents unchanged.
On a dominant exit the probability in \eref{e:cellflux} is $n^{-o(1)}$, uniformly in the starting point.
The exponents are integers. Hence, on the formal strip,
\begin{equ}\label{e:fluxorder}
 \OO\to\OO' \quad\Longrightarrow\quad
 \max_{\substack{(j,k)\in\OO,\ r\ {\rm active}\\(j,k)+\zeta_r\in\OO'}}(y_{r,1}-j)
 =e(\OO)\leq e(\OO').
\end{equ}
Along the proposed path the exponent cannot decrease. Yet
$e(\OO+(\ell,0))=e(\OO)-\ell$, a contradiction.

A merge cannot contain both $\OO_m^{(0)}$ and $\OO_m^{(0)}+(\ell,0)$ for $\ell>0$, since the set containing $\OO_m^{(0)}$ before the merge would reach its right translate.
Translated components that contain the same $\OO_m^{(0)}$ coincide. Thus each lumped set contains at most $s$ sets of the initial partition, and each round with a merge reduces the number of sets up to horizontal translation.
There are at most $s-1$ such rounds.
At a merge the old sets have the same exponent by \eref{e:fluxorder}. Join them by crossing reactions of that exponent.
Every exit from the union was an exit from an old set, so the exponent can only decrease.
By induction, the sets $\OO_m^{(0)}\subseteq\OO$ are connected by reaction edges of exponent at least $e(\OO)$.
A simple path between two of these sets has at most $s-1$ edges. Thus the horizontal width is at most $(s-1)\Kz$.
Only finitely many shapes occur, so one sufficiently large $\gamma$ suffices in \lref{l:exittimes}.
\end{proof}

\subsection{Construction of fibers} \label{s:fibers2}

\begin{definition}\label{d:translation} A directed graph $\mathcal G$ embedded in $\sz$ is \emph{translation invariant} if there exists $p \in \mathbb N\setminus \{0\}$ and $n_0 \in \mathbb {N}$ such that for any $v, v' \in \sz$ with $v_1 , v_1' \geq n_0$
  \begin{equ}
    v \to v' \in E(\mathcal G) \Leftrightarrow v+(p,0) \to v' + (p,0) \in E(\mathcal G)
  \end{equ}
  The smallest $p \in \mathbb N$ for which the above is true is called the \emph{period} of the graph $\G$.
\end{definition}

\begin{definition}\label{d:fibers}
  A set of \emph{fibers} is a family of sets $\{\mathcal F_n\}_{n =n_0}^\infty$ of the form
  \begin{equ}\label{e:fiber}
    \mathcal F_n := \{(n+f(i),i)~:~i \in \{0, \dots, h(n)\}\}
  \end{equ}
  for a \emph{skeleton} $f~:\Z_{\ge0} \to \mathbb Z$ and $n_0$ such that $\mathcal F_n \subset \sg$ for all $n\geq n_0$.
  We further say that such a set of fibers \emph{orders}
  $\sg$ \abbr{wrt}
  a translation invariant graph $\mathcal G$ embedded in $\sg$ if for any $n,n' \geq n_0$, having $v \Rightarrow v'$ for $v\in \mathcal F_n, v' \in \mathcal F_{n'}$ implies $n\geq n'$.
\end{definition}
In other words, $\{\mathcal F_n\}_{n =n_0}^\infty$ orders
$\sg$ \abbr{wrt}
$\G$ if there is no directed path connecting (any point of) $\mathcal F_{n'}$ to (any point of) $\mathcal F_n$ if $n > n'$. We use the same notation for the restrictions of these fibers to $S$.
\begin{lemma}\label{l:fibers}
  Let $\G$ be a directed acyclic graph that is embedded in the strip $\sz$ for a $K\in \mathbb Z_{\geq 0}$. Further assume that $\G$ is translation invariant (with period $1$), has uniformly bounded horizontal edge displacements, and that for any $v,v' \in V(\mathcal G)$ joined by a nonempty path $v \Rightarrow_{\G} v'$,  $v_2 = v_2'$ implies $v_1 > v_1'$. Then there exists a set of fibers $\{\mathcal F_n\}_{n =n_0}^\infty$ ordering
  $\sz$ \abbr{wrt}
  $\mathcal G$.\\

\end{lemma}

The proof of this lemma adapts the depth-first search algorithm for the topological sorting problem (\cite{tarjan76}) to our periodic setting. Here and throughout, for a vertex $v \in V(\G)$ we define the set of children of $v$ as
\begin{equ}
  C(v):= \{v' \in V(\G)~:~v \to v' \in E(\G)\}
\end{equ}

We prove this result by constructing $\mathcal F_n$ for a fixed $n>n_0$ through an iterative procedure. Note that this is equivalent to finding the skeleton $f$, which in turn sufficient to construct the full set of fibers $\{\mathcal F_n\}$. The intuition of the proof is as follows: we will define a guess for the set $\mathcal F_n$ that consists of the union of two sets $\mathcal E \cup \mathcal E'$ where, morally speaking, $\mathcal E'$ is a list of vertices whose children are, by construction, in $\mathcal E \cup \mathcal E'$ or to its left\footnote{meaning: they can be written as $(m,i)$ and for all $(m',i) \in \mathcal E\cup \mathcal E'$ we have $m'>m$.}  and $\mathcal E$ is a list of vertices whose children may or may not be in $\mathcal E \cup \mathcal E'$ or to its left.
Our procedure then iterates over the vertices of the graph and sequentially updates the sets $\mathcal E, \mathcal E'$ respecting the above intuition, until $\mathcal E = \emptyset$ and $\mathcal E'$ is a fiber, so that by the definition of $\mathcal E'$ the fiber $F_n := \mathcal E'$ orders $\G$.

We start our procedure defining our first ansatz for $\mathcal F_n$ as the set $\mathcal E = \{(n+f_0(k),k)~:~k \in \{0,1,\dots, K\}\}$ for a choice of initial skeleton $f_0~:~\{0,1,\dots, K\} \to \mathbb Z$ and define the set $\mathcal E' := \emptyset$. We also define the set $\tilde {\mathcal E} = \emptyset$ with the operational purpose of keeping track of the vertices that have already been ``visited'' by the iterative procedure but whose children are not (yet) guaranteed to be in $\mathcal E \cup \mathcal E'$ or to its left. 
At the beginning of each epoch choose the element $v\in\mathcal E$ with the highest second coordinate and apply the following recursive procedure:
\begin{enumerate}
  \item Add $v$ to $\tilde{\mathcal E}$. For each $v'\in C(v)$, let $v''$ be the unique element of $\mathcal E\cup\mathcal E'$ with $v''_2=v'_2$. If $v'_1>v''_1$, remove $v''$ from whichever of these sets contains it, add $v'$ to $\mathcal E$, and apply the procedure to $v'$.
  \item After all children have been considered, remove $v$ from $\mathcal E$ and $\tilde{\mathcal E}$ and add it to $\mathcal E'$.
\end{enumerate}
Repeat until $\mathcal E=\emptyset$.

In words, the above procedure performs a sort of depth-first search in the graph $\G$. For every node $v \in \mathcal E$ it marks the node as ``seen'' (including it in $\tilde{\mathcal E}$) and proceeds to consider all its children $v' \in C(v)$: if the child is in $\mathcal E \cup \mathcal E'$ or to its left this means that our present guess $\mathcal E \cup \mathcal E'$ orders $\sz$ \abbr{wrt} the edge $v \to v'$, else the ``guess'' $\mathcal E \cup \mathcal E'$ needs to be updated so that this condition is satisfied by adding $v'$ to $\mathcal E$, $\tilde{\mathcal E}$ and starting to consider $C(v')$. Once all recursive child visits from a vertex $v$ are complete, all its offspring is in $\mathcal E \cup \mathcal E'$ or to its left, and therefore $v$ will be added to $\mathcal E'$ and removed from $\mathcal E, \tilde{\mathcal E}$.

\begin{proof}

We now show that this procedure terminates, i.e. that after a finite number of steps $\mathcal E = \emptyset$. We do this by considering the state of $\mathcal E, \mathcal E'$ every time we have $\tilde{\mathcal E} = \emptyset$. This corresponds to a moment in the procedure outlined above where we have explored \emph{all} the children of the node that was first added to $\tilde{\mathcal E}$ the last time we had $\tilde{\mathcal E} = \emptyset$. We refer to these moments of the procedure as \emph{epochs} we denote by $\mathcal E_k, \mathcal E_k'$ the sets $\mathcal E, \mathcal E'$ at the end of the $k$-th epoch. By the fact that $|\mathcal E \cup \mathcal E'| = K+1 < \infty$, to prove the claim that the procedure terminates in finitely many steps it is sufficient to show that $|\mathcal E_k'|$ is monotonically increasing in $k$ and that the number of operations in each epoch is finite.
We proceed to prove both claims, starting from the second.

Let $R < \infty$ be a uniform bound on the degree of each vertex in $\G$. To prove the finiteness claim, unfold the recursive calls into a tree. This tree has degree at most $R$ and depth at most $K+1$, and therefore contains at most $1+R+\cdots+R^K$ vertices. It remains to show the bound on the depth. Now,
\begin{enumerate}
  \item Note that by the definition of our recursion, an edge $v \to v'$ with $v'$ in $\tilde{\mathcal E} \subseteq \mathcal E \subseteq \mathcal E \cup \mathcal E'$ is neglected by the procedure. In particular, if a node $(m,i) \in \tilde{\mathcal E}$ none of its offspring (children of any degree in the tree) will be of the form $(m',i)$ for $m' \leq m$.
  \item Further, recall that by our assumption $v \Rightarrow_{\G} v'$ with  $v_2 = v_2'$ implies $v_1 > v_1'$, so that if $(m,i) \in \tilde{\mathcal E}$ none of its offspring (children of any degree in the tree) will be of the form $(m',i)$ for $m' > m$.
\end{enumerate}
From the points above it follows that if $(m,i) \in \tilde{\mathcal E}$ none of its offspring (children of any degree in the tree) will be of the form $(m',i)$. Since the vertices of the tree are a subset of $\sz$, i.e., points of the form $(m,i)$ for $i \in \{0,1,\dots, K\}$ we have proved our claim.

We now proceed to prove the monotonicity of $|\mathcal E_k'|$. First, we note that $|\mathcal E_k'|$ is a nondecreasing function of $k$: a height already represented in $\mathcal E_k'$ may be reopened during the next epoch, but its representative can only move right and is returned to $\mathcal E'$ before that epoch ends. We further note that once a point is added to $\tilde{\mathcal E}$ that point (or another point to its right) will be added to $\mathcal E'$ by the end of the epoch. Since at the beginning of an epoch $k$ we add $v \not \in \mathcal E_k'$ to $\tilde{\mathcal E}$ we must have that $v \in \mathcal E_{k+1}'$ so that $|\mathcal E_{k+1}'| > |\mathcal E_{k}'|$.
The representatives only move right, so the edge inequalities at other processed vertices are preserved. At termination these inequalities give the required ordering, after increasing $n_0$ so that no path between retained fibers can leave the periodic region and return. Indeed, if $J_G$ bounds horizontal edge displacements, any path increases its first coordinate by at most $KJ_G$, since every return to a previous height strictly decreases that coordinate. For $f_0=0$, every update follows a path from an initial vertex, and the final edge inequalities give $f(k)\ge d(\pi)$ for every path $\pi$ in the periodic height graph ending at height $k$, where $d(\pi)$ is its horizontal displacement. Thus $f(k)$ is the maximal such displacement. Deleting repeated-height segments increases displacement, so a maximizing path visits each height at most once.
\end{proof}

We now show that, choosing initial guess for $\mathcal F_n$ as the vertical set $\mathcal E := \{(n,i)~:~i \in \{0,1,\dots, K\}\}\subset V(\G)$, under the following condition the procedure for constructing $\mathcal F_n$ leaves the upper part of $\mathcal E$ unchanged. On the coarse graph of \sref{s:fibers1}, we instead start all upper bands to the right of the lower points and denote their common column by $n$.
We use one vertex per residue band in each column in WE and retain the cap targets as singletons in SE. Within each lumped set we take the connecting edges of \lref{l:noforward} in both directions, and for each dominant exit we take a crossing edge from \eref{e:fluxorder}.
\begin{lemma}\label{l:straightontop} Let $\mathcal F_n$ be constructed as detailed above. If $f_0\equiv0$ and there exists $K_{\G}\in\{0,\ldots,K-1\}$ such that
\begin{equ}
  v \to v' \not \in E(\G) \qquad \text{if $v_2>K_{\G}$ and either } \begin{cases}v_1 \leq v_1'\text{ and } v_2<v_2'\\
  v_1 < v_1'\text{ and } v_2>v_2'\end{cases}\,,
\end{equ}
then letting $\rho := \max_{v \to v' \in E(\G)} \|v'-v\|_1$, if $K > (K_{\G} + 3) \rho^2 $ we have that
\begin{equ}
  (n, i) \in \mathcal F_n\qquad \text{for all } i \in \{K-2 \rho, K-2 \rho +1, \dots, K-1\}\,.
\end{equ}

For the coarse graph on an infinite boundary portion of $S$, there is a $K_*$ such that, for $K\ge K_*+2\Kz$ and $L=(s-1)\Kz$, the choice $f_0(k)=-L\mathbf1_{\{k<K_*\}}$ gives
\begin{equ}\label{e:celllabel}
 -L\le f(k)\le0,\qquad f(k)=0\quad(k\ge K_*),\qquad N(x)=x_1-f(x_2).
\end{equ}
If a path of dominant transitions starts in $\mathcal F_n$ and ends in $\mathcal F_{n'}$, then $n'\le n$.
\end{lemma}

\begin{proof}
For $f_0\equiv0$, let
\[
 I_0=\{i\in\{0,\ldots,K_{\G}\}:(n,i)\in\mathcal F_n\}
    =\{i\in\{0,\ldots,K_{\G}\}:f(i)=0\}.
\]
Write $R=\max_v|C(v)|$ for the maximal degree of the graph under consideration, and recall $\rho=\max_{v\to v'\in E(\G)}\|v'-v\|_1$.
If $f(k)>0$, choose a shortest path $\gamma:(n,i)\Rightarrow(n+f(k),k)$ attaining the maximal displacement $d(\gamma)=f(k)$, as in \lref{l:fibers}. Its first step moves right, since otherwise deleting that step and translating the remaining path to column $n$ gives a shorter path of at least the same displacement. Hence $i\le K_{\G}$, and
\[
 f(k)\ge f(i)+d(\gamma)=f(i)+f(k),\qquad f(i)\ge0.
\]
Thus $i\in I_0$. If $f\equiv0$, then $I_0=\{0,\ldots,K_{\G}\}$, so $I_0\ne\emptyset$ in either case.
A return to the same height moves strictly left by \lref{l:fibers}, so $\gamma$ visits each height at most once. Since only steps from heights at most $K_{\G}$ can increase the first coordinate,
\[
 f(k)\le\rho\,\#\{t:0\le t<|\gamma|,\ (\gamma_t)_2\le K_{\G}\}
       \le (K_{\G}+1)\rho .
\]

On the coarse graph a path from $\OO_m^{(0)}$ to $\OO_m^{(0)}+(\ell,0)$ has $\ell\le0$ by \lref{l:noforward}. Deleting such returns leaves at most $s-1$ edges. From each initial vertex there are at most $1+R+\cdots+R^{s-1}$ such paths, and
\[
 d(\gamma)\le(s-1)\Kz=L,\qquad
 f(k)=\max_{\gamma:i\Rightarrow k}\{f_0(i)+d(\gamma)\}.
\]
Indeed, the maximum bounds every update, so the rightward updates terminate. The reverse inequality follows by summing $f(k')-f(k)\ge j'-j$ along the edges $(j,k)\to(j',k')$ of any path. This also gives $N(j',k')\le N(j,k)$. For $i<K_*$, $f_0(i)+d(\gamma)\le -L+L=0$.

For $i\ge K_*$, let $j$ be the starting column and put $m=\max_{r:\zeta_{r,1}\ne0}y_{r,1}$. Weak reversibility gives $y_{r,1},y'_{r,1}\le m$ whenever $\zeta_{r,1}\ne0$, since a return path must leave the product column. Let $\OO$ contain an upper point in column $j$. On an infinite boundary portion the lattice contains $(-b,0)$ for some $b>0$ (in rank one the span is horizontal). By weak reversibility there are reactions $r_1,\ldots,r_q$ with $\sum_{p=1}^q\zeta_{r_p}=(-b,0)$. Applying them in this order moves the process $b$ columns to the left and returns it to its starting height. First take enough degree-$m$ leftward steps that all subsequent repetitions of $r_1,\ldots,r_q$ stay in columns at most $j-m$. This path leaves $\OO$, whose width is bounded, and every step up to exit has exponent at least $m-j$. Its vertical excursion is bounded. In SE choose $K_*$ large enough to enable the path, counting a cap crossing as an exit. In WE start sufficiently high in the same residue band, taking $K_*=K_c$ and then $n$ large enough to stay below the cap. In either case $e(\OO)\ge m-j$.

By \eref{e:fluxorder} and the choice of the connecting edges, every edge of $\gamma$ then has exponent at least $m-j$. A first crossing into a column $j'>j$ would give, for a reaction from column $a$,
\[
 m-j\le y_{r,1}-a=y'_{r,1}-(a+\zeta_{r,1})\le m-j'<m-j.
\]
Thus $d(\gamma)\le0$ for $i\ge K_*$. The path bound now gives $f_0(k)\le f(k)\le0$, which is \eref{e:celllabel}.

Now consider, for $f_0\equiv0$, a point $(m,i)\in\mathcal F_n$ with $i>K_{\G}$, $m>n$. Take the path $\gamma$ from $(n,i_0)$ with $i_0\in I_0$ constructed above. Let $\ell$ be its last timestep with $(\gamma_\ell)_2\leq K_{\G}$. After this step, each upward step decreases the first component by at least $1$ and increases the second by at most $\rho-1$. If $q$ is the number of these upward steps, then
\[
 i\leq K_{\G}+\rho+q(\rho-1),\qquad
 0<m-n\leq(K_{\G}+1)\rho-q.
\]
Hence
\[
 i\leq K_{\G}+\rho+\bigl((K_{\G}+1)\rho-1\bigr)(\rho-1)
 \leq(K_{\G}+3)\rho^2-2\rho<K-2\rho,
\]
which proves the claim.
\end{proof}

\subsection{Exit time estimates and proof of \pref{p:finalboundary}}\label{s:boundaryfinal}

Fix $R<\infty$ and $K$ sufficiently large. Take $H\ge1$ fixed in SE and $H=\lfloor K+\gamma\log n\rfloor+O(1)$ with $\gamma$ fixed and sufficiently large in WE.
\begin{lemma}\label{l:exittimes}
 Let $\OO_n\subset S\cap\{|x_1-n|\le R,\ 0\le x_2\le H\}$ and let $\tau_0(n)$ be its exit time. There exists $\bar\beta>0$ such that, for all sufficiently large $n$,
 \begin{equ}\label{e:uniontime}
  \sup_{x\in\OO_n}\E_x\tau_0(n)\le n^{\bar\beta}.
 \end{equ}
 In WE, as $n\to\infty$,
 \begin{equ}\label{e:unioncap}
  \sup_{x\in\OO_n}\mathbb P_x\bigl(X_2(\tau_0(n))>H\bigr)=o(1).
 \end{equ}
 For the fibers $\mathcal F_n$ of \lref{l:straightontop}, in WE,
 \begin{equ}\label{e:highreturn}
  \sup_{\substack{x\in\mathcal F_n\\x_2\ge K}}
  \mathbb P_x\!\left(\tau_n\le\inf\{t\ge0:X_2(t)<K\}\right)=o(1),
  \qquad \tau_n=\inf\{t>0:X(t)\notin\mathcal F_n\}.
 \end{equ}
\end{lemma}
\begin{proof}
Let $\{\OO_{n,j}\}$ be a partition of $\OO_n$, let $\tau_j(n)$ be the exit time from $\OO_{n,j}$ and let $Q$ be the kernel at successive exits, absorbed on $\OO_n^c$. We may assume $\sup_{x\in S}x_1=\infty$, since otherwise $\OO_n$ is eventually empty. For each fixed $n$, every point of the finite set $\OO_n\subsetneq S$ can reach its complement, so the mean exit times are finite.
For \eref{e:uniontime} it suffices to find a partition and a fixed $m\ge1$ satisfying
 \begin{equ}\label{e:exittimes1}
   \sup_{x \in \OO_{n,j}}\mathbb E_x{\tau_j(n)} \leq n^{\beta_j}\qquad \inf_{x \in \OO_n}Q^m(x, \OO_n^c) > n^{-\beta_0}
 \end{equ}
with $\sup_j\beta_j<\infty$.
This is basic renewal theory. By the Markov property, the expected exit time must satisfy
  \begin{equs}
     \sup_{x \in \OO_{n}}\mathbb E_x{\tau_0(n)} & \leq m \sup_j \sup_{x \in \OO_{n,j}}\mathbb E_x{\tau_j(n)}   + \sup_{x \in \OO_{n}} \sum_{x' \in \OO_n} Q^m(x, x') \mathbb E_{x'}{\tau_0(n)}\\
    & \leq m \sup_j \sup_{x \in \OO_{n,j}}\mathbb E_x{\tau_j(n)} + (1-\inf_{x \in \OO_n}Q^m(x, \OO_n^c)) \sup_{x' \in \OO_n}\mathbb E_{x'}{\tau_0(n)}
  \end{equs}
  Rearranging and plugging in the bounds \eref{e:exittimes1} yields
  \begin{equ}
    \sup_{x \in \OO_{n}}\mathbb E_x{\tau_0(n)} \leq m (\inf_{x \in \OO_n}Q^m(x, \OO_n^c))^{-1} \sup_j \sup_{x \in \OO_{n,j}}\mathbb E_x{\tau_j(n)} \leq m \cdot n^{\beta_0 + \sup_j \beta_j}.
  \end{equ}
The same calculation applies to $Y$, with exit time $\tau_{\OO_n^c}$ in jumps.

In WE take $A=\OO_n\cap\{x_2\ge K\}$ as one part, if nonempty, and take singletons below $K$. In SE take only singletons.

For each singleton height $k$, fix $(a,k)\in S$. By irreducibility there is a path from $(a,k)$ to a point with first coordinate greater than $a+2R$. Its translate from $(n+j,k)$ is enabled for large $n$ and leaves the window. There are finitely many such heights, so by \eref{e:transprob} and \eref{e:rates} these paths take at most $m-1$ jumps and have probability greater than $n^{-\beta_0}$ for some fixed $m,\beta_0$. In WE the first exit from $A$ either leaves $\OO_n$ or reaches a singleton. This gives the second bound in \eref{e:exittimes1}. A singleton takes one jump to leave, and \lref{l:noexitfromtop} gives $\sup_{x\in A}\E_x\tau_{A^c}\le CH$ in WE. This gives the first bound in jumps with $\beta_j=1$. Since every active state's total rate is at least $\kappa_-$,
\[
 \sup_{x\in\OO_n}\E_x\tau_0(n)
 \le \kappa_-^{-1}\sup_{x\in\OO_n}\E_x\tau_{\OO_n^c}
 \le \kappa_-^{-1}m n^{\beta_0+1}.
\]
This proves \eref{e:uniontime} for $\bar\beta>\beta_0+1$.

In WE each excursion entered from below $K$ starts at height at most $K-1+\Kz$. There are at most $\E_x\tau_{\OO_n^c}$ such entries in expectation. Apply the uniform bound of \lref{l:noexitfromtop} to the initial excursion and its bound for $K\le k\le K+\Kz$ to each subsequent excursion. With $\beta=\gamma$ this gives
\[
 \mathbb P_x\bigl(X_2(\tau_0(n))>H\bigr)
 \le C(\log n)^{-1}+Cn^{-\gamma}\E_x\tau_{\OO_n^c}=o(1)
\]
for $\gamma>\beta_0+1$. Finally take $\OO_n=\mathcal F_n$ and $H=h_n$. A vertical jump below $K$ stays in the fiber since $f=0$ for $k\ge K-2K_0$. Before descending below $K$, the fiber can therefore be left only by crossing the cap or by a horizontal jump. In WE each horizontal reaction has first source coordinate below the maximum. By \eref{e:ratesub} its probability per jump is at most $Ca(H)/n$. Together with \lref{l:noexitfromtop}, this bounds the probability in \eref{e:highreturn} by $C(\log n)^{-1}+CHa(H)/n=o(1)$. For caps in neighboring columns differing by $O(1)$, stop the initial excursion also at its first horizontal jump and use the lowest cap for excursions entered from below $K$.
\end{proof}
We now have all the tools to prove the main result of the section.

\begin{proof}[Proof of \pref{p:finalboundary}]
  In this proof we realize the program outlined at the beginning of \sref{s:boundary}.\\[10pt]
 We consider the horizontal boundary in $S$. The argument for the vertical boundary is identical. If $x_1$ is bounded on $S$, or the stoichiometric span is one-dimensional and nonhorizontal, this boundary portion is finite. Otherwise the stoichiometric lattice contains a nonzero horizontal vector.
 We start by considering the set $\bar B_{\gamma,n}(K)$ with $\gamma = 0$ for $(1,0)$-SE networks and $\gamma > 0$ (to be fixed later) otherwise. For each initial set $\OO_m^{(0)}$ in column $j$, the corresponding set in $S\cap B_{\gamma,n}(K)$ is
 \begin{equ}
  \bigl((n,0)+\OO_m^{(0)}\bigr)\cap S\cap B_{\gamma,n}(K)
  =\{(j,k)_n\in S:(j,k)\in\OO_m^{(0)},\ 0\le k\le h_n\},
  \qquad |j|\le2K\Kz .
 \end{equ}
We note that the exit estimates required to apply \lref{l:exitestimate1} hold, for large enough $n$, either by \rref{r:singletons} or by \lref{l:commonexit} and \lref{l:exittimes}. Hence, by \lref{l:exitestimate1}, \cref{c:fw}, $W$-graphs control the exit probabilities from any union of these sets. We then iterate the coarse-graining procedure described in point (1) of our program: maximal strongly connected components of the graph connecting $\{\OO_m^{(k)}\}$ with the most likely transitions between those sets are lumped to form $\{\OO_m^{(k+1)}\}$. This procedure combines $\omega$-limit sets on shorter timescales into larger ``lumps'' and identifies the most likely exit points from such ``lumps''. Here, \lref{l:coarsening} guarantees that the required properties of such coarse-grained sets (mainly independence of most likely exits on initial condition within the set) are preserved by the coarse-graining procedure. \\[10pt]
We stop the procedure described above when we have obtained a Directed Acyclic Graph (DAG) expressing the most likely transitions connecting the sets $\{\OO_m^{(k)}\}$. Termination and bounded width follow from \lref{l:noforward}. By \lref{l:straightontop}, the index $N$ in \eref{e:celllabel} cannot increase along a dominant transition. Extend $f$ by zero at omitted heights. Since $s-1\le K+\Kz$, we have $\|f\|_\infty+\Kz\le s\Kz\le 2K\Kz$. In WE cut fiber $n$ at $h_n$. Neighboring caps differ by $O(1)$, which is covered by enlarging the fixed collars in \lref{l:exittimes} and \lref{l:commonexit}.

Let $\tau_n=\inf\{t>0:X(t)\notin\mathcal F_n\}$ and $\Delta=N(X(\tau_n))-n$. Since the skeleton and reaction jumps are bounded, $|\Delta|\le D$ for a fixed $D$. Every vertex has a zero-cost exit and the dominant graph within each fiber is acyclic, so the minimum cost in the denominator is zero. Every path to a higher label has positive cost. Thus \cref{c:fw} gives probability $o(1)$ for ordinary exits to a higher label. In WE, cap exits are $o(1)$ by \lref{l:exittimes}, after choosing $\gamma$ large enough. In SE this includes cap exits to a higher label, since the cap targets were included in the ordering. Before reaching $[K-K_0,K]$ no cap exit is possible. In this band the fiber is vertical and a maximal-degree leftward reaction exits to a smaller label with probability at least $c_K>0$. Hence
\[
 \mathbb P_x(\Delta<0)\geq c_K\{1-\mathbb P_x(\Delta>0)\}\qquad((1,0)\text{-SE}).
\]
Consequently, in both cases, for all sufficiently large $n$,
\begin{equ}\label{e:endofstory}
 \limsup_{n\to\infty}\sup_{x\in\mathcal F_n~:~x_2\leq K}
 \mathbb P_x(\Delta>0)=0,
 \qquad \inf_{x\in\mathcal F_n~:~x_2\leq K}\mathbb P_x(\Delta<0)\geq c>0.
\end{equ}
In WE, \eref{e:highreturn} and the strong Markov property at the first entrance into $\{x_2<K\}$ extend \eref{e:endofstory} to every $x\in\mathcal F_n$.

We now proceed to prove the pivotal bound \eref{e:mainbound2}.
With $V$ and $u$ as in \eref{e:lyapunovfunction} and \eref{e:entropy-u}, write $B(n)=(u(n)+u(h_n))^{\alpha+1}/(\alpha+1)$ and set
\[
 a_n=u(n)^\alpha\log n,\qquad
 r_n=(n\log n)^\alpha\log(h_n+2)=o(a_n).
\]
Since $h_n=O(\log n)$ and $|h_{n+\ell}-h_n|\leq1$ for bounded $\ell$ and large $n$, Taylor expansion gives, uniformly for $|\ell|\leq D$,
\[
 B(n+\ell)-B(n)=\ell a_n+O(r_n).
\]
For an exit point $y$ in the boundary layer, $V(y)=B(N(y))$ where we recall $N(y)=y_1-f(y_2)$. For an exit into the bulk, since $f=0$ on the upper part of the fiber, $y_1=N(y)$ and $0<y_2-h_{N(y)}\leq K_0+1$. Hence, in both cases,
\[
 0\leq V(y)-B(N(y))\leq Cr_n.
\]
Since $V(x)=B(n)$ on $\mathcal F_n$, stopped Dynkin's formula and the uniform extension of \eref{e:endofstory} give
\begin{equs}
 \E_x\left[\int_0^{\tau_n}\L V(X(s))\,ds\right]
 &=\E_x V(X(\tau_n))-B(n)\\
 &\leq a_n\E_x\Delta+Cr_n\\
 &\leq-a_n\bigl(\mathbb P_x(\Delta<0)
                   -D\mathbb P_x(\Delta>0)\bigr)+Cr_n\\
 &\leq-c_1(n\log n)^\alpha\log n.\label{e:final}
\end{equs}
By \eref{e:uniontime}, $\sup_{x\in\mathcal F_n}\E_x\tau_n\le n^{\bar\beta}$ in both cases. Choosing $\alpha>\bar\beta$ and comparing with \eref{e:final} completes the proof.
\end{proof}

\section{Bulk estimate and patching}\label{s:patch}
This section is dedicated to the proof of the dissipation inequality in the bulk of the state space, stated in Proposition~\ref{p:totalbulk0}. 

The proof is carried out in three steps. 
First we study the regions of the bulk where the CRN is strongly endotactic. Here, the Lyapunov function defined in~\eqref{e:lyapunovfunction} is dissipated by a dominant reaction, and therefore after one jump of the process, see Lemma~\ref{l:NonHH}. The proof uses the notion of tier-sequence, see~\cite{anderson2020tier}. 
Then we study the weakly endotactic regions, where the dominant reactions do not dissipate the Lyapunov function. The process has to be observed in a faster timescale, allowing sub-dominant reactions to occur. The result is given in Proposition~\ref{2DScalingDif1}, and uses scaling methods and standard results of stochastic calculus. 
Finally, we check that the dissipation inequality still holds at the states near the frontier with the boundary set, where the Lyapunov function changes expression, see Lemma~\ref{l:s4-patching}.

\subsection*{Notations}
The bulk estimates are independent of Sections~\ref{s:boundary} and \ref{s:boundary2}, and require some specific notations. 

For $w\in S^1$, we introduce the set of $w$-dominant complexes:
\[
 \mathcal C_w=\mathop{\rm arg\,max}_{y\in\mathcal C}\langle y,w\rangle.
\]
At a state $(N,N^{w_2/w_1})$ with $N\gg1$, the reactions starting from
complexes in $\mathcal C_w$ are dominant. The dissipation inequality is
shown for states such that $x_1\geq x_2$. By symmetry of the role of the
two variables, the whole inequality can then be deduced. We introduce
\[
 \Omega_\gamma(K)=\{x\in\mathbb Z_{\geq0}^2:x_1\geq x_2\}
                            \setminus B_\gamma(K),
 \qquad S_>=\{w\in S^1:w_1>0\quad\text{and}\quad w_1\geq w_2\geq 0\}.
\]
For a two-dimensional weakly reversible CRN
$\mathcal X=(\mathcal S,\mathcal C,\mathcal R)$, we set
\[
 \mathcal P_{\mathcal X}
 =\{w\in S_>:\mathcal R_w\cap\bar{\mathcal R}_w=\emptyset\},
\]
where $\mathcal R_w$ and $\bar{\mathcal R}_w$ are introduced in relation~\eqref{e:dominant}. This set is finite: each such $w$ is orthogonal to a nonzero dominant reaction vector. For $c>0$ and
$w\in S_>$, the set of states $\mathcal H_c(w)$ is defined as 
\begin{equation}\label{e:s4-H}
 \mathcal H_c(w)=\left\{x\in\mathbb Z_{>0}^2: e^{-c}\leq\frac{x_2}{x_1^{w_2/w_1}}\leq e^c\right\}.
\end{equation}

With $u$ as in \eref{e:entropy-u}, put $v(x)=u(x_1)+u(x_2)$ and write the entropy function used in the bulk as
\begin{equation}\label{e:s4-VS}
 V_{\mathcal S}(x)=\frac{v(x)^{\alpha+1}}{\alpha+1}.
\end{equation}
Here $\alpha\geq1$ is taken sufficiently large for the boundary estimates, see Proposition~\ref{p:finalboundary}. Outside the two boundary layers and $A_0$, \eref{e:lyapunovfunction} gives $V=V_{\mathcal S}$.

\subsection{Strongly endotactic regions}
In the regions where one of the dominant reactions dissipates the
entropy function, the system is point-wise dissipating.

\begin{lemma}\label{l:NonHH}
There exist some constants $\{c_w:w\in\mathcal P_{\mathcal X}\}$,
$\beta>0$, $K\geq2K_0$ and a finite set $\mathcal K$ such that
\begin{equation}\label{e:NonHH}
 \mathcal L V_{\mathcal S}(x)\leq-\beta
\end{equation}
for every
\begin{equation}\label{e:StronglyEndotacticSet}
 x\in\Omega_0(K)\setminus
       \left(\mathcal K\cup
           \bigcup_{w\in\mathcal P_{\mathcal X}}\mathcal H_{c_w}(w)\right).
\end{equation}
\end{lemma}
This lemma sets the width of the weakly endotactic regions. The proof is done by contradiction and uses the notion of \emph{tier sequence}, see~\cite{anderson2020tier}. We first prove the inequality for a CRN with one linkage class, then deduce the result when the CRN has several linkage classes.

\begin{proof}
Note that, for a reaction vector $a$ and $x_1\geq x_2>K>2K_0$, Taylor expansion gives
\begin{equation}\label{e:Taylor}
 \frac{V_{\mathcal S}(x+a)-V_{\mathcal S}(x)}{v(x)^\alpha}
 =\langle a,\log x\rangle+R_a(x),\qquad
 |R_a(x)|\leq C\left(\frac1{x_1}+\frac1{x_2}\right)
                        + C_\alpha\frac{(\log x_1)^2}{v(x)}.
\end{equation}
Here $C$ is independent of $\alpha$, and the last term tends to zero as $x_1\to\infty$ for each fixed $\alpha$. Also, uniformly over the reactions, for $K$ large enough,
\begin{equation}\label{e:LambdaBoundcc}
 \lambda_r(x)=\kappa_r x^{y_r}
 \left(1+O(x_1^{-1}+x_2^{-1})\right),\qquad
 c x^{y_r}\leq\lambda_r(x)\leq Cx^{y_r}.
\end{equation}

We first prove the result for a CRN with one linkage class. In this case $|\mathcal P_{\mathcal X}|\leq1$. 

Choose the constants $c_w$ below and, in the horizontal weakly endotactic case, an auxiliary $c_{(1,0)}$ with $\log K \geq c_{(1,0)}$, and assume by contradiction that there is a sequence $(x_p)$ such that $\|x_p\|\to\infty$, for all $p$, $\mathcal L V_{\mathcal S}(x_p)\geq-1$, and
\[
 x_p\in\Omega_0(K)\setminus \bigcup_{w\in\mathcal P_{\mathcal X}}\mathcal H_{c_w}(w).
\]
For all $p$, define $\varepsilon_p=\log(x_{p,2})/\log(x_{p,1})$. Up to a subsequence, this sequence converges, and we set
\[
 \varepsilon=\lim_{p\to\infty}\varepsilon_p, \qquad \omega=(1,\varepsilon)\quad\text{and}\quad  \bar\omega=\omega/\|\omega\|.
\]

\smallskip
\noindent\emph{If $\bar\omega\in\mathcal P_{\mathcal X}$:}
The sequence converges to the side of a weakly endotactic region. Up to taking a further subsequence, we can assume that it stays on one side, i.e., that the sequence
\[
 \nu_p=(\varepsilon_p-\varepsilon)\log(x_{p,1}) =\log(x_{p,2})-\varepsilon\log(x_{p,1})
\]
verifies either $\nu_p\geq c_{\bar\omega}$ for all $p$, or $\nu_p\leq-c_{\bar\omega}$ for all $p$. Assume that we are in the first situation. The second one is dealt with similarly. When $\varepsilon=0$, $\nu_p=\log x_{p,2}>\log K\geq c_{(1,0)}$, so only the first situation occurs. 

Since there is one linkage class and $\bar\omega\in \cal{P}_{\cal{X}}$, $\mathcal C_{\bar\omega}=\mathcal C$. Set $\rho_\omega=y_1+\varepsilon y_2$, which is the same for every $y\in\mathcal C$, and
\[
 y_2^*=\max_{y\in\mathcal C}y_2\quad\text{and}\quad 
 g(t)=t\sum_{r\in\mathcal R}\kappa_r e^{t y_{r,2}}\zeta_{r,2}.
\]
Relation~\eqref{e:Taylor} gives
\begin{equation}\label{e:s4-g-remainder}
 \frac{\mathcal L V_{\mathcal S}(x_p)}{v(x_p)^\alpha x_{p,1}^{\rho_\omega}}
 =g(\nu_p)+E_p,\qquad \text{with} \quad |E_p|\leq \left(C\frac{1+\log x_{p,2}}{x_{p,2}} +C_\alpha\frac{\log x_{p,1}}{x_{p,1}}\right) \sum_{y\in\mathcal C}e^{\nu_p y_2}.
\end{equation}
Every reaction from the complex with second coordinate $y_2^*$ decrease $x_2$. Moreover, they dominate the sum in $g(t)$ for $t$ large. Therefore, for $c_{\bar\omega}$ sufficiently large, for some constant $c>0$, $g(t)\leq-c t e^{ty_2^*}$ holds for all $t\geq c_{\bar\omega}$. Taking $K>0$ large enough, we can have $|E_p|\leq c t e^{ty_2^*}/2$ for $p$ sufficiently large. This leads to
\[
 \mathcal L V_{\mathcal S}(x_p) \leq-\frac{c}{2} v(x_p)^\alpha x_{p,1}^{\rho_\omega} |\nu_p|e^{\nu_p y_2^*}\longrightarrow-\infty.
\]
This is a contradiction. 

\smallskip\noindent\emph{If $\bar\omega\notin\mathcal P_{\mathcal X}$:}
Using the conventions of~\cite{anderson2020tier}, up to taking a further subsequence, $(x_p)$ is a tier sequence: every ratio $x_p^y/x_p^{y'}$ has a limit in $[0,\infty]$. We introduce the set $T^1_{x_p}$ of dominant complexes along this sequence. For $y\in T^1_{x_p}$, all ratios $x_p^{y'}/x_p^y$ are bounded, and when $y'\notin T^1_{x_p}$ this ratio tends to zero. Ratios between complexes in $T^1_{x_p}$ converge to finite positive numbers.

We decompose
\[
 \mathcal L V_{\mathcal S}(x_p)=D_1(x_p)+D_2(x_p)+D_3(x_p),
\]
where $D_1$ contains the reactions $y\to y'$ with $y\notin T^1_{x_p}$, i.e.,
\[
 D_1(x)=\sum_{\substack{y\to y'\in\mathcal R\\y\notin T^1_{x_p}}}
       \lambda_{y\to y'}(x)
          \bigl(V_{\mathcal S}(x+y'-y)-V_{\mathcal S}(x)\bigr),
\]
$D_2$ contains the reactions with $y,y'\in T^1_{x_p}$, and $D_3$
the reactions with $y\in T^1_{x_p}$ and $y'\notin T^1_{x_p}$.

Fix $y^*\in T^1_{x_p}$. Relations~\eqref{e:Taylor} and~\eqref{e:LambdaBoundcc} and the definition of the set $T^1_{x_p}$ give directly
\begin{equation}\label{e:s4-D12}
 \limsup_{p\to\infty}\frac{D_1(x_p)}{x_p^{y^*}v(x_p)^\alpha}\leq 0,
 \qquad
 \limsup_{p\to\infty}\frac{D_2(x_p)}{x_p^{y^*}v(x_p)^\alpha}<\infty.
\end{equation}

The definition of the first tier gives $T^1_{x_p}\subset\mathcal C_{\bar\omega}$. Since $\bar\omega\notin\mathcal P_{\mathcal X}$, $\mathcal C_{\bar\omega}\ne\mathcal C$. Therefore $T^1_{x_p}$ is a proper subset of $\mathcal C$. By weak reversibility and since the CRN has only one linkage class, there is a reaction from this set to its complement, i.e., $D_3$ is nonempty. For every such reaction,
\[
 \frac{V_{\mathcal S}(x_p+y'-y)-V_{\mathcal S}(x_p)}{v(x_p)^\alpha}
 =\log\frac{x_p^{y'}}{x_p^y}+O(1)\longrightarrow-\infty.
\]
Its rate divided by $x_p^{y^*}$ is bounded away from zero. Hence
\begin{equation}\label{e:s4-D3}
 \frac{D_3(x_p)}{x_p^{y^*}v(x_p)^\alpha}\longrightarrow-\infty,
\end{equation}
and summing up all these estimates, we get $\mathcal L V_{\mathcal S}(x_p)\longrightarrow-\infty$. This leads to a contradiction and completes the proof for a CRN with one linkage class.

\smallskip\noindent\emph{Several linkage classes.}
If $\mathcal X$ is composed of $\ell$ linkage classes, it can be seen as the union of $\ell$ CRNs $(\mathcal X_1,\ldots,\mathcal X_\ell)$, where $\mathcal X_i=(\mathcal S,\mathcal C_i,\mathcal R_i)$ is the CRN corresponding to linkage class $i$. We decompose the infinitesimal generator as follows:
\begin{equation}\label{e:DecompositionIG}
 \mathcal Lf(x)=\sum_{i=1}^\ell\mathcal L_i f(x),\qquad
 \mathcal L_i f(x)=\sum_{r\in\mathcal R_i}\lambda_r(x)
                                      (f(x+\zeta_r)-f(x)).
\end{equation}
The single linkage class case gives a common $K$ and a finite set such that $\mathcal L_iV_{\mathcal S}(x)\leq-1$ outside the cone $\mathcal H_{c_i}(w_i)$ if $\mathcal P_{\mathcal X_i}=\{w_i\}$, and everywhere in $\Omega_0(K)$ outside that finite set if $\mathcal P_{\mathcal X_i}=\emptyset$.

To deduce~\eqref{e:NonHH}, note that
\[
 \mathcal P_{\mathcal X}\subset\bigcup_{i=1}^\ell\mathcal P_{\mathcal X_i}.
\]
For $w\in\mathcal P_{\mathcal X}$ choose $c_w=\max\{c_i:w_i=w\}$. Applying~\eqref{e:DecompositionIG} proves the desired inequality except for states in cones $\mathcal H_{c_i}(w_i)$ with $w_i\notin\mathcal P_{\mathcal X}$.

To get the inequality in these regions, we carry out once again the tier-sequence method. By contradiction, introduce a sequence $(x_p)\in\mathcal H_{c_i}(w_i)$ such that $\|x_p\|\to\infty$ and $\mathcal L V_{\mathcal S}(x_p)\geq-1$. Up to a subsequence it is a tier sequence for the full set $\mathcal C$. The bounds of the cone show that $T^1_{x_p}=\mathcal C_{w_i}$. Since $w_i\notin\mathcal P_{\mathcal X}$, $\mathcal R_{w_i}\cap\bar{\mathcal R}_{w_i}\neq\emptyset$ and therefore, $D_3$ is not empty. 
The bounds~\eqref{e:s4-D12} now
applied to the full network, together with~\eqref{e:s4-D3}, give a contradiction as previously.
This concludes the proof.
\end{proof}

\subsection{A scaling argument for the weakly endotactic regions}\label{s:s4-scaling}
In this part, we find a stopping time for the states in the weakly
endotactic regions. For $\bar\omega\in\mathcal P_{\mathcal X}$ with $\bar\omega_2>0$, we use the normalization $\omega=(1,a)$, where $a={\bar\omega}_2/{\bar\omega}_1\in(0,1]$. 
If $(1,0)\in\mathcal P_{\mathcal X}$, take $K>e^{c_{(1,0)}}$. Then $\mathcal H_{c_{(1,0)}}((1,0))\subset B_0(K)$, and this cone belongs to the boundary region.
We set
\begin{equation}\label{e:s4-rho}
 \rho_0=\max_{y\to y'\in\mathcal R_{\bar\omega}}\langle y,\omega\rangle \quad\text{and}\quad
 \rho_d=\max_{y\to y'\in \bar{\mathcal{R}}_{\bar\omega}}\langle y,\omega\rangle.
\end{equation}
\begin{remark}\label{r:Decrease}
    Since the CRN is ${\bar\omega}$-WE, $0<\rho_d<\rho_0$, and for $y\to y'$ such that $\langle y, \omega\rangle =\rho_d$, $\langle y'-y,\omega\rangle\leq 0$. 
\end{remark}

We first treat the case $a<1$.  
\begin{proposition}\label{2DScalingDif1}
If $0<a<1$ and $c>0$, one can find stopping times $\tau_N>0$ and
some $0<b<1$ such that
\begin{equation}\label{e:s4-scaling1}
 \begin{split}
 \limsup_{N\to\infty}\sup_{(N,k_N)\in \mathcal H_c({\bar\omega})}
 \frac{\mathbb E_{(N,k_N)}V_{\mathcal S}(X(\tau_N))}
      {(N\log N)^{\alpha+1}}
 &\leq\frac{b^{\alpha+1}}{\alpha+1},\\
 \limsup_{N\to\infty}\sup_{(N,k_N)\in \mathcal H_c({\bar\omega})}
             \frac{\mathbb E_{(N,k_N)}\tau_N}{N}&=0.
 \end{split}
\end{equation}
Moreover, the process up to $\tau_N$ stays in a set $\cal{H}_{c'}({\bar\omega})$ for some $c'>c$. 
\end{proposition}

We focus on the states $x_N=(N,k_N)\in \cal{H}_c(\bar\omega)$, i.e., such that
\[
 e^{-c}N^a\leq k_N\leq e^cN^a,
\]
and we call $(X^N(t))$ the process starting at $x_N$.
We introduce 
\begin{equation}
    S^N(t):= \croc{X^N(t), \omega}=X_1^N(t)+aX_2^N(t). 
\end{equation}
Note that if $X_2^N$ remains of order $N^a$ and $X_1^N$ of order $N$, then
\[
 \frac{v(X^N(t))}{N\log N}
       =\frac{X_1^N(t) +X_2^N(t)}{N}+o(1)=\frac{S^N(t)}{N}+o(1).
\]
We first show that $(S^N(t)/N)$ decreases significantly and then deduce the dissipation of the entropy. The fastest reactions changing $(S^N(t))$ are the reactions in $\bar{\cal{R}}_{\bar\omega}$, which happen at rate $N^{\rho_d}$, and decrease the process, see remark~\ref{r:Decrease}.
We need to see $O(N)$ such reactions to get a significant decrease of $(S^N(t)/N)$: the relevant timescale is therefore $t\mapsto t/N^{\rho_d-1}$.

The proof is carried out in three steps. We first show that $(X_2^N(t))$ stays of order $N^a$, using a coupling with a queueing process. We then show that $(S^N(t/N^{\rho_d-1})/N)$ cannot increase by a fixed amount. Finally, we show that it has decreased after a time $t_0>0$.

We use the following elementary bound for the queue comparison.
\begin{lemma}\label{l:s4-queue}
Let $Z(0)=0$, and let $Z$ be the chain on $\mathbb Z_{\geq0}$ with
transitions
\[
 z\longrightarrow z+K_0\quad\text{at rate }\lambda,
 \qquad z\longrightarrow z-1\quad\text{at rate }\mu
                          \text{ if }z\geq1.
\]
If $K_0\lambda<\mu$, there exist $C,\theta>0$ such that
\begin{equation}\label{e:s4-queue}
 \mathbb P\left(\sup_{s\leq T}Z(s)\geq q\right)
                      \leq C(1+T)e^{-\theta q},\qquad T\geq0, q\geq K_0.
\end{equation}
\end{lemma}
\begin{proof}
Choose $\theta>0$ sufficiently small that $\lambda(e^{\theta K_0}-1)+\mu(e^{-\theta}-1)<0$. During an excursion away from zero, $e^{\theta Z}$ is a supermartingale. Optional stopping bounds the probability that an excursion starting at $K_0$ reaches level $q$ by $e^{\theta K_0-\theta q}$. The expected number of excursions starting before time $T$ is at most $\lambda T$. Applying the strong Markov property at their starting times and summing these bounds concludes the proof. 
\end{proof}

\begin{proof}[Proof of Proposition~\ref{2DScalingDif1}]
Set $\eta=\rho_d-1>-1$. The process $X^N$ can be written as
\begin{equation}\label{2DSDE}
 X^N(t)=x_N+\sum_{r\in\mathcal R}\zeta_r
       \int_0^t\int_0^\infty
       \mathbf1_{\{u\leq\lambda_r(X^N(s-))\}}\,
                                      \mathcal P_r(du,ds),
\end{equation}
where $(\mathcal P_r)$ are independent Poisson random measures
with intensity $ds\,du$. 

We set $S^N(t)=\langle X^N(t),\omega\rangle$ and
\[
 (\overline S^N(t),\overline X_1^N(t),\overline X_2^N(t))
 =\left(\frac{S^N(t/N^\eta)}N,
        \frac{X_1^N(t/N^\eta)}N,
        \frac{X_2^N(t/N^\eta)}{N^a}\right).
\]
Fix $m_1=1/2$ and $M_1=2$. Let $m_2,M_2>0$ be chosen below, with
$2m_2<e^{-c}$ and $M_2/2>e^c$. Define
\[
 H_{i+}^N=\inf\{t\geq0:\overline X_i^N(t)\geq M_i\},\qquad
 H_{i-}^N=\inf\{t\geq0:\overline X_i^N(t)\leq m_i\},\qquad
 H_i^N=H_{i+}^N\wedge H_{i-}^N.
\]
For a constant $t_0>0$, to be specified later, set
\begin{equation}\label{2DstoppingTime}
 \overline\tau_N=t_0\wedge H_1^N\wedge H_2^N,
 \qquad \tau_N=\tau_N/N^\eta.
\end{equation}
We first prove that, for some $0<b<1$,
\begin{equation}\label{e:s4-prob-drop}
 \inf_{(N,k_N)\in \mathcal H_c({\bar\omega})}
 \mathbb P_{(N,k_N)}(\overline S^N(\overline\tau_N)\leq b)
                                                        \longrightarrow1.
\end{equation}

\smallskip\noindent\emph{Control of $H_2^N$: uniform bound of $(\overline X_2^N(t))$.}
We show that
\begin{equation}\label{e:s4-H2}
 \mathbb P(H_2^N\leq t_0\wedge H_1^N)\longrightarrow 0, 
\end{equation}
by upper and lower bounding the process $(X_2^N(t))$ using a coupling argument. We only write the upper bound argument, the lower bound is done the same way.  

For $m_1N\leq x_1\leq M_1N$ and $m_2N^a\leq x_2\leq M_2N^a$, the rate of a reaction $r$ starting from $y_r$ satisfies
\begin{equation}\label{e:s4-scaled-rates}
 \lambda_r(x)=\kappa_r N^{\croc{y_r, \omega}}\left(\frac{x_1}{N}\right)^{y_{r,1}}\left(\frac{x_2}{N^a}\right)^{y_{r,2}}(1+O(N^{-a})).
\end{equation}
After $m_2,M_2$ are fixed, the total rate of the non-dominant reactions is $o(N^{\rho_0})$ uniformly on the box.

Let $y_2^*=\max\{y_2:y\in\mathcal C_{\bar\omega}\}$. By ${\bar\omega}$-WE, a reaction $y\to y'\in \cal{R}_{\bar\omega}$ such that $y_2=y_2^*$ verifies $y_2'-y_2<0$. We call $r^*$ such a reaction. On $M_2N^a/2\leq x_2\leq M_2N^a$, the rate of the reaction $r_*$ is at least $\mu N^{\rho_0}$, where
\[
 \mu=c(M_2/2)^{y_2^*}.
\]
Similarly, a dominant reaction that increases $x_2$ starts from a complex with $y_2\leq y_2^*-1$. For all sufficiently large $N$, the total rate of the reactions icreasing $x_2$, including the slower reactions, is at most $\lambda N^{\rho_0}$ throughout the box, where we can take
\[
 \lambda=2C M_2^{y_2^*-1}.
\]
Here $c,C>0$ depend on $m_1,M_1$ and the network, and not on $M_2$.
We can choose $M_2$ sufficiently large so that $K_0\lambda<\mu$.

Let $Z_+$ be the queue of Lemma~\ref{l:s4-queue} with rates $\lambda,\mu$. Using similar arguments as in~\cite{LR24-2}, we can construct a coupling that verifies, up to $\tau_N$,
\begin{equation}\label{e:s4-upper-coupling}
 X_2^N(s)\leq \lceil M_2N^a/2\rceil+K_0 +Z_+(sN^{\rho_0}). 
\end{equation}
Therefore
\[
 \mathbb P(H_{2+}^N\leq t_0\wedge H_1^N\wedge H_{2-}^N)
 \leq\mathbb P\left(\sup_{s\leq t_0N^{\rho_0-\eta}}Z_+(s)
                         \geq M_2N^a/2-K_0-1\right)\longrightarrow0, 
\]
which concludes the proof of the upper-bound. For the lower bound, a similar argument can be used choosing $m_2$ small, and this proves~\eqref{e:s4-H2}.

\smallskip\noindent\emph{Control of $H_{1+}^N$: upper bound of $(\overline X_1^N(t))$.}
The rates of the reactions are uniformly bounded in the box, for $N$ large enough, as 
\[
    \lambda_r(x)\leq C_rN^{\croc{y_r, \omega}}, 
\]
for some constant $C_r$ depending on the $m_i, M_i$, but independent of $N$. 
Using the SDE~\eqref{2DSDE}, we get
\begin{multline}\label{2DUpperBoundSN}
 \overline S^N(t\wedge\overline\tau_N) \leq\overline S^N(0)+R_N,\\
 \text{with}\quad R_N:=\frac{1}{N}\sum_{\substack{r\in \cal{R}:\\ \croc{\zeta_r, \omega}>0}} \croc{\zeta_r, \omega} \cal{P}_r([0, C_rN^{\croc{y_r, \omega}}]\times [0, t_0N^{-\eta}]). 
\end{multline}
Any reaction $y\to y'$ increasing $(S^N(t))$ verifies $\croc{y, \omega}<\rho_d$, see remark~\ref{r:Decrease}. Therefore, $R_N\to 0$ in probability, and $(\overline{S}^N(t))$ doesn't increase. 

Since $\overline X_1^N\leq\overline S^N$ and $\overline S^N(0)=1+O(N^{a-1})$, this gives
\[
 \mathbb P(H_{1+}^N\leq t_0\wedge H_{1-}^N\wedge H_2^N)\longrightarrow0.
\]

\smallskip\noindent\emph{Decrease of $(\overline S^N(t))$.}
Choose a reaction $r_d$ with $\croc{y_{r_d}, \omega}=\rho_d$ and $\croc{\zeta_{r_d}, \omega}=-\delta_d<0$. On the stopped box, its rate is at least $c_dN^{\rho_d}$ for a constant $c_d>0$. The random variable
\[
 Y_N=\frac1N\cal{P}_{r_d}
          \bigl([0,c_dN^{\rho_d}]\times [0,t_0N^{-\eta}]\bigr)
\]
has mean $c_dt_0$ and variance $c_dt_0/N$. On the event $\overline\tau_N=t_0$, we therefore have
\begin{equation}\label{2DeqDecreaseSN}
 \overline S^N(\overline\tau_N)
       \leq\overline S^N(0)+R_N-\delta_dY_N.
\end{equation}
We choose $t_0$ so that $0<\delta_dc_dt_0<1-m_1$. If $\overline\tau_N=H_{1-}^N$ and no other exit occurs, then
$\overline S^N(\overline\tau_N)\leq m_1+aM_2N^{a-1}+O(1/N)$. The other exits have probability tending to zero by the first two steps. Using $R_N\to0$ in probability and $Y_N\to c_dt_0$ in probability, we obtain~\eqref{e:s4-prob-drop}, for example with $b=1-\delta_dc_dt_0/3$.

\smallskip\noindent\emph{Conclusion of the proof.}
For states in the stopped box, uniformly,
\[
 \frac{v(y)}{N\log N}
    =\frac{y_1}{N}+o(1)
    =\frac{\langle y,\omega\rangle}{N}+o(1).
\]
The scaled processes are bounded until $\tau_N$, we can therefore raise $\overline{S}^N(t)$ at the power $\alpha+1$. It follows that
\[
 \limsup_{N\to\infty}\sup_{x_N}
 \frac{\mathbb E_{x_N}V_{\mathcal S}(X^N(\tau_N))}
      {(N\log N)^{\alpha+1}}
 \leq\frac{b^{\alpha+1}}{\alpha+1}.
\]
The second limit in~\eqref{e:s4-scaling1} is straightforward, since $\tau_N\leq t_0N^{1-\rho_d}$ and $\rho_d>0$.
\end{proof}

When $a=1$, both $(X_1^N(t))$ and $(X_2^N(t))$ move at the same timescale and have to be controlled by a coupling argument. 
\begin{proposition}\label{2DScaling1}
If $a=1$ and $c>0$, one can find stopping times $\tau_N>0$ and some $a_0>0$ such that
\begin{equation}\label{e:s4-scaling-diagonal}
 \begin{split}
 \limsup_{N\to\infty}\sup_{e^{-c}\leq k_N/N\leq1}
 \frac{\mathbb E_{(N,k_N)}V_{\mathcal S}(X(\tau_N))
                      -V_{\mathcal S}(N,k_N)}{(N\log N)^{\alpha+1}}
 &\leq-a_0,\\
 \limsup_{N\to\infty}\sup_{e^{-c}\leq k_N/N\leq1}
                  \frac{\mathbb E_{(N,k_N)}\tau_N}{N}&=0.
 \end{split}
\end{equation}
Moreover, the stopped process, stays in a set $\cal{H}_{c'}({\bar\omega})$ for some $c'>c$. 
\end{proposition} 
The proof uses the same arguments as the one from Proposition~\ref{2DScalingDif1}. 
\begin{proof}
We keep the previous notations. Setting $s_N=1+k_N/N$, we take $M_1=M_2=M$ and
$m_1=m_2=m$, with $M>0$ large and $m>0$ small. We fix $\delta\in(0,e^{-c}/4)$, choose
$m<e^{-c}/4$, and add the stopping time
\[
 H_{S-}^N=\inf\{t:\overline S^N(t)\leq s_N-\delta\},\qquad
 \overline\tau_N=t_0\wedge H_1^N\wedge H_2^N\wedge H_{S-}^N,
 \quad \tau_N=\overline\tau_N/N^{\rho_d-1}.
\]

\smallskip\noindent\emph{Bound of the process.}
Using the same arguments as in the proof of Proposition~\ref{2DScalingDif1}, we can show that both processes $(X_1^N(t))$ and $(X_2^N(t))$ stay in their box with high probability. The upper bound is given by relation~\eqref{2DUpperBoundSN}, taking $M=3$. Their lower bound is given by a coupling with a queue as the lower bound version of relation~\eqref{e:s4-upper-coupling}. 

More precisely, for a lower bound of $(\overline{X}_2^N(t))$ for example, for $t\leq\overline\tau_N$ we have $\overline S^N(t)>s_N-\delta$. Therefore, if $X_2^N(t)\leq2mN$, then
\[
 X_1^N(t)\geq(s_N-\delta-2m)N\geq N.
\]
Let $y^*=\arg\min\{y_2, \, y\in \cal{C}_{\bar\omega}\}$.  By $\bar\omega$-WE, a reaction $y^*\to (y^*)'$ has $(y^*)'\in \cal{C}_{\bar\omega}$, and therefore has $y_2^*>(y_2^*)'$.  By~\eqref{e:s4-scaled-rates} with $a=1$, for $N\leq x_1\leq3N$ and $mN\leq x_2\leq2mN$ its rate is at least $\mu N^{\rho_0}$, with $\mu=c\,m^{y_2^-}$. A dominant reaction that decreases $x_2$ starts from a complex $y\in\mathcal C_{\bar\omega}$ with $y_2\geq y_2^-+1$. For all sufficiently large $N$, the total rate of the reactions decreasing $x_2$, including the slower reactions, is therefore at most $\lambda N^{\rho_0}$ on the same set, with $\lambda=2C(2m)^{y_2^-+1}$. Here $c,C>0$ depend on the network and not on $m$. We choose $m$ sufficiently small so that $K_0\lambda<\mu$. Let $Z_-$ be the queue of Lemma~\ref{l:s4-queue} with rates $\lambda,\mu$. As for~\eqref{e:s4-upper-coupling}, we can construct a coupling that verifies, up to $\tau_N$,
\[
 X_2^N(s)\geq\lfloor2mN\rfloor-K_0-Z_-(sN^{\rho_0}).
\]
Since $k_N\geq e^{-c}N>4mN$, Lemma~\ref{l:s4-queue} gives
\[
 \mathbb P(H_{2-}^N\leq\overline\tau_N)
 \leq\mathbb P\left(\sup_{s\leq t_0N^{\rho_0-\eta}}Z_-(s)\geq mN-K_0-1\right)\longrightarrow0, 
\]
which concludes the proof of the lower bound of $(\overline{X}_2^N(t))$. The lower bound of $(\overline{X}_1^N(t))$ is shown the same way, switching the two species. 

\smallskip\noindent\emph{Decrease of $(\overline S^N(t))$ and of the entropy.} On this box, we choose the decreasing reaction $r_d$ and the constant $\delta_d,c_d$ as before, and take $0<\delta_dc_dt_0<\delta$. At $t_0$, we use~\eqref{2DeqDecreaseSN}. At $H_{S-}^N$, the sum has decreased by $\delta$. Therefore, for $a_* =\delta_dc_dt_0/3$, uniformly in the initial state,
\[
 \mathbb P(\overline S^N(\overline\tau_N)\leq s_N-a_*)\longrightarrow1.
\]
Since $v(y)/(N\log N)=(y_1+y_2)/N+o(1)$, we have
\[
 \frac{\mathbb E_{x_N}V_{\mathcal S}(X(\tau_N))-V_{\mathcal S}(x_N)}
 {(N\log N)^{\alpha+1}}
 \leq \frac{(s_N-a_*)^{\alpha+1}-s_N^{\alpha+1}}{\alpha+1}+o(1)<-a_0,
\]
for $a_0>0$ small enough. 
\end{proof}

We can now conclude with the dissipation inequality needed for Proposition~\ref{p:totalbulk0}:
\begin{corollary}\label{2DCorollaryScaling}
For $\bar\omega\in\mathcal P_{\mathcal X}$ with $\bar\omega_2>0$ and $c>0$, for every sufficiently
large $x=(N,k_N)\in \mathcal H_c({\bar\omega})$ with $k_N\leq N$,
one can find an integrable stopping time $\tau_N>0$ such that
\begin{equation}\label{2DCorollaryEq}
 \mathbb E_xV_{\mathcal S}(X(\tau_N))-V_{\mathcal S}(x)
                                       \leq-\mathbb E_x\tau_N.
\end{equation}
\end{corollary}
\begin{proof}
The corollary is a direct consequence of Proposition~\ref{2DScalingDif1} and~\ref{2DScaling1}
\end{proof}

\subsection{Between the bulk and the boundary}\label{s:s4-gluing}
The dissipation inequality remains to be shown for the states in the set 
\[
 \mathcal I=\{(n,k)\in S:n\geq k,\ K+\gamma\log n<k\leq K+\gamma\log n +K_0\}, 
\]
for which one jump can lead to the boundary definition of the function, which prevents the use of Lemma~\ref{l:NonHH}. 

\begin{lemma}\label{l:s4-patching}
For $K$ sufficiently large, there exist $\beta>0$ and a finite set
$\mathcal K$ such that
\begin{equation}\label{e:s4-patching}
 \mathcal LV(n,k)\leq-\beta,\qquad (n,k)\in\mathcal I\setminus\mathcal K.
\end{equation}
The same inequality holds with the coordinates exchanged.
\end{lemma}

\begin{proof}
\smallskip\noindent\emph{The $(1,0)$-WE case:} Here, $\gamma>0$. We introduce $y^*=\arg\max\{y_2: y\in \cal{C}_{(1,0)}\}$. Note that $y_1^*=\max_{y\in \cal{C}} y_1$, that by $(1,0)$-WE, any reaction in $\cal{R}_{(1,0)}$ is vertical, and that the reaction starting from $y^*$ must decrease the second coordinate.  
We set $\kappa_-=\min_{r\in \cal{R}} \kappa_r$ and $\kappa_+=\max_{r\in \cal{R}}\kappa_r$. For a state $(n, k)\in \cal{I}$, we have:
\begin{align*}
    \cal{L}(V)(n,k)\leq &\frac{\kappa_-}{2} n^{y_1^*}k^{y_2^*}(V(n, k-1)-V(n,k))\\
        &+ \kappa_+ n^{y_1^*}k^{y_2^*-1}|\cal{R}_{(1,0)}|(V(n, k+K_0)-V(n,k))\\
        &+\kappa_+ n^{y_1^*-1}k^{K_0}|\cal{R}|(V(n+K_0, k+K_0)-V(n,k))
\end{align*}
Using Taylor expansion, we can find $C, c>0$ such that, for $n$ large enough and $(n,k)\in \cal{I}$, 
\[
    \begin{aligned}
        \frac{V(n, k-1)-V(n,k)}{v(n,k)^\alpha}&\leq -c\log(\log n),\\
        \frac{V(n, k+K_0)-V(n,k)}{v(n,k)^\alpha}&\leq C\log(k),\\
        \frac{V(n+K_0, k+K_0)-V(n,k)}{v(n,k)^\alpha}&\leq C\log(n).
    \end{aligned}
\]
These four inequalities show that for $n$ large enough, $k\leq K+K_0+\gamma\log n$, we have
\[
    \frac{\cal{L}(V)(n,k)}{v(n,k)^\alpha}\leq -\beta, 
\]
for some $\beta>0$, which leads to inequality~\eqref{e:s4-patching}.

\smallskip\noindent\emph{The $(1,0)$-SE case:} Here, $\gamma=0$. We set $y^*\to (y^*)'\in \bar{\cal{R}}_{(1,0)}$. $y^*\to (y^*)'\in \cal{R}_{(1,0)}$ by $(1,0)$-SE, therefore $y_1^*=\max_{y\in \cal{C}}y_1$, and $y_1^*>(y_1^*)'$. Moreover, for $y\to y'\in \cal{R}_{(1,0)}$, $y_1\geq y_1'$. 
Therefore, 
\begin{align*}
    \cal{L}(V)(n,k)\leq &\frac{\kappa_-}{2} n^{y_1^*}k^{y_2^*}(V(n-1, k+K_0)-V(n,k))\\
        &+ \kappa_+ n^{y_1^*}k^{K_0}|\cal{R}_{(1,0)}|(V(n, k+K_0)-V(n,k))\\
        &+\kappa_+ n^{y_1^*-1}k^{K_0}|\cal{R}|(V(n+K_0, k+K_0)-V(n,k))
\end{align*}
Since $K<k\leq K+K_0+1$, Taylor expansion gives, uniformly on this finite set of heights,
\[
\begin{aligned}
 \frac{V(n-1,k+K_0)-V(n,k)}{v(n,k)^\alpha}
   &=-\log n+O_K(1),\\
 \frac{V(n,k+K_0)-V(n,k)}{v(n,k)^\alpha}
   &=O_K(1),\\
 \frac{V(n+K_0,k+K_0)-V(n,k)}{v(n,k)^\alpha}
   &=O_K(\log n).
\end{aligned}
\]
The negative term in the preceding generator bound therefore dominates the other two terms for $n$ sufficiently large, which proves~\eqref{e:s4-patching}.
\end{proof}

\begin{proof}[Proof of Proposition~\ref{p:totalbulk0}]
For $x$ in the region of patching, take the first jump time $t_1$: Lemma~\ref{l:s4-patching} gives
\[
 \mathbb E_x(V(X(t_1))-V(x))
       =\frac{\mathcal LV(x)}{\lambda(x)}
       =\mathcal LV(x)\mathbb E_xt_1\leq-\beta\mathbb E_xt_1.
\]
At a bulk state outside this region, all one-jump endpoints remain in the bulk, so $\mathcal LV(x)=\mathcal LV_{\mathcal S}(x)$. If $x$ is outside the weakly endotactic cones, we apply Lemma~\ref{l:NonHH} and again take $\sigma_x=t_1$.

For states in $\mathcal H_{c_w}(w)$ with $w_2>0$, use the stopping time of
Corollary~\ref{2DCorollaryScaling}. The stopped process stays in a larger cone $\cal{H}_{c_w'}(w)$, and therefore outside the logarithmic boundary strips when $N$ is sufficiently large. On this whole stopped set $V=V_{\mathcal S}$, and~\eqref{2DCorollaryEq} gives the desired inequality.
\end{proof}

\section{Notation table}

For convenience, we include here a table with the most common notation utilized throughout the paper, in order of appearance in the text.

\begin{center}
\begin{longtable}{|l|l|}
\caption{Table containing the most common notation used throughout the paper.} \label{tab:notation} \\

\hline \multicolumn{1}{|c|}{\textbf{Notation}} & \multicolumn{1}{c|}{\textbf{Definition}} \\ \hline 
\endfirsthead

\multicolumn{2}{c}%
{{\bfseries \tablename\ \thetable{} -- continued from previous page}} \\
\hline \multicolumn{1}{|c|}{\textbf{Notation}} & \multicolumn{1}{c|}{\textbf{Definition}} \\ \hline 
\endhead

\hline \multicolumn{2}{|r|}{{Continued on next page}} \\ \hline
\endfoot

\hline \hline
\endlastfoot

$\Re_w$ & $\arg\max\{\dtp{y,w}~:~y \to y' \in \Re\}$ \\
$\bar{\Re}_w$ & $\arg\max\{\dtp{y,w}~:~y \to y' \in \Re, \dtp{y'-y,w} \neq 0\}$ \\
$\Re(v)$ & $\{y\to y'\in\Re\,:\,y'-y=v\}$ \\
$\kappa_+$ & $\max_{y\to y'\in\Re} \kappa_{y\to y'}$ \\
$\kappa_-$ & $\min_{y\to y'\in\Re} \kappa_{y\to y'}$ \\
$B_\gamma(K)$ & $\{(n,k) \in \mathbb Z_{\geq 0}^2~:~ n\ge1,\ k \leq K + \gamma \log n\}$\\
$B_\gamma'(K)$ & $\{(k,n) \in \mathbb Z_{\geq 0}^2~:~ n\ge1,\ k \leq K + \gamma \log n\}$ \\
$K_0$ & $1+\max_{r}\{\|y_r\|_1+\|\zeta_r\|_1\}$ \\
$C(k)$ & $\{y : y_2\leq k, \exists y\to y'\in\Re\}$ \\
$z(k)$ & $\arg \max\{z_1 : z \in C(k)\}$ \\
$y(k, v)$ & $\arg \max\{y_1 : y \in C(k)\text{ and } \exists y \to y' \in \Re(v)\}$ \\
$\ss(K)$  & $\{(j,k ) \in\mathbb Z_{\geq0}^2~:~0\leq k\leq\lfloor K+\gamma\log n\rfloor, j \in [n- 2 K \Kz, n+2K \Kz]\}$\\
$\bar B_{\gamma,n}(K)$  & $\{(j,k ) \in \mathbb Z_{\geq 0}^2 ~:~k \in [0, (K+ \gamma \log n)+\Kz ], j \in [n- 2 K \Kz, n+2K \Kz]\}$ \\
$Q_W(x,A)$ & $\mathbb P(Y(\tau_W) \in A~|~Y(0) = x )$ \\
\end{longtable}
\end{center}

\section*{Statement on AI use}

This paper is the result of a project that lasted more than 5 years. The bulk estimates in Section 4 can be found in the PhD thesis of LL, while the boundary estimates were established in a joint work among all the remaining authors. AI was used to check the proofs and to fix minor technical mistakes (mainly resulting in edited proofs of Lemmas 3.6 and 3.7), but did not contribute conceptually to the development of the ideas contained in this paper. While proofreading, the authors enforced a minimal change policy with respect to the original human written manuscript, and take full responsibility of the results presented in this paper.

\section*{Acknowledgements}

AA and LL acknowledge the partial support of the Faculty of science of the University of Bern, in particular the MJP project ``Stochastic models of aging''. AA acknowledges the partial support of the Swiss National Science Foundation
project grant 2000-1-243054

\bibliographystyle{plain}

\appendix

\section{Exit estimates for Lemmas~\ref{l:noexitfromtop} and~\ref{l:commonexit}}\label{s:exitproofs}

For both lemmas, it is sufficient to prove the desired result for $j=0$: by substituting $n\to n+j$ for $|j|\leq |2K\Kz|$ in the result, we obtain uniformity in $j$ upon possibly increasing the constants. 

We set throughout the proofs of the two lemmas the following notations: 
\begin{equ}\label{e:ystarWE}
    y_1^*:= \max_{y\in\C}y_1, \quad \cal{R}_0=\{y\to y'\in\Re: \, y_1=y_1^*\} \quad \text{and}\quad y^*\to (y^*)'\in \arg\max_{y\to y'\in\cal{R}_0}^{\mathrm{lex}}(y_2,y'_2). 
\end{equ}
noting that $y_1^*=z_1(\Kz)$ and $\cal{R}_0=\cal{R}_{(1,0)}$ where $\cal{R}_{(1,0)}$ is defined in Relation~\eqref{e:dominant}. 
Note also that the exit probability from a set $W^c$ through a set $A$ is a solution $\psi~:~\bar B_{n,\gamma}(K) \to \mathbb R$ to the boundary value problem 
\begin{equ}\label{e:poissonineq}
    \begin{cases}
        \mathcal L \psi (x) = 0\qquad &\text{for } x \in W^c\\
        \psi(x) = \mathbbm 1_{A}(x)\qquad &\text{for } x \in W
    \end{cases}
\end{equ}
where $\mathcal L $ is as in \eqref{eq:generator} \cite{norris98}. Hence, solving the inequality
\begin{equ}\label{e:poissonub}
    \begin{cases}
        \mathcal L \psi_+ (x) \leq 0\qquad &\text{for } x \in W^c\\
        \psi_+(x) \geq \mathbbm 1_{A}(x)\qquad &\text{for } x \in W
    \end{cases}
\end{equ}
gives an upper bound to $\psi(x)$ since for $x\in W^c$, 
\begin{equs}
    \psi_{+}(x) & = \mathbb E_{x}(\psi_+(X(\tau))) - \mathbb E_{x}\left[\int_0^\tau \mathcal L \psi_+(X(s)) ds\right] \\
    & \geq \mathbb E_{x}(\psi_+(X(\tau))) \geq \psi(x), \\
\end{equs}
where $\tau$ is the exit time of $W^c$. 
Similarly, solving the inequality 
\begin{equ}\label{e:poissonlb}
    \begin{cases}
        \mathcal L \psi_- (x) \geq 0\qquad &\text{for } x \in W^c\\
        \psi_-(x) \leq \mathbbm 1_{A}(x)\qquad &\text{for } x \in W
    \end{cases}
\end{equ}
gives a lower bound to $\psi(x)$.

\begin{proof}[Proof of \lref{l:noexitfromtop}]
    Fix $\gamma>0$, put $\beta=\gamma$ and first take $H=\lfloor K+\gamma\log n\rfloor$. For the cap probability we use \eref{e:poissonub} on $W^c$, with boundary value $1$ above $H$ and $0$ on the other exits. We choose $\psi_+$ as
    \begin{equ}
        \psi_+((j,k)_n)=
        \begin{cases}
            \psi_{n,+}(k) &\text{if }K\le k\le H,\\
            \mathbbm 1_{\{k>H\}} &\text{otherwise},
        \end{cases}
    \end{equ}
    with $\psi_{n,+}: k\in \Z_{\geq 0}\mapsto C(\log n)^{-1}n^{-\beta}\exp[\beta k/\gamma]$, for some $C>0$ specified later. This function dominates the boundary values on $W$.

    To bound the generator at $x=(j,k)_n\in W^c$, we write $j=0$. The rate bounds below hold uniformly for $|j|\le R$ by replacing $n$ with $n+j$ in the rates. Since $\psi_+$ depends only on $k$,
    \begin{equation}\label{eq:alvrfleejhcbdjhsac}
        \mathcal L \psi_{+}(x)=\sum_{y\to y'\in\Re}\lambda_{y\to y'}(x)
           \bigl(\psi_{+}(x+(0,y'_2-y_2))-\psi_{+}(x)\bigr).
    \end{equation}
    By \eref{e:ystarWE} $y^*\to (y^*)'$ is a dominant reaction at $(n,k)$, and since the CRN is $(1,0)$-WE then necessarily $y^*_1=(y^*)'_1$. Moreover, $(y_2^*)'<y_2^*$ else by weak reversibility there would exist a reaction from $(y^*)'$ leading to a contradiction with definition \eqref{e:ystarWE}. Since $\psi_{n,+}$ is monotone increasing and strictly positive,
    \begin{equation}\label{eq:lerbvkherbvjhbdvcjhvsdcj}
        \psi_{+}(x+(y^*)'-y^*)-\psi_{+}(x)< 0.
    \end{equation}

    For $n,K$ large enough,
    \begin{equ}\label{e:InLambda}
        \lambda_{y^*}((n, k))\geq \frac{1}{2}n^{y_1^*}k^{y_2^*}\quad \text{and} \quad \lambda_y((n,k))\leq n^{y_1}k^{y_2},\, \text{for all  } y\in \cal{C}. 
    \end{equ}
    Combining \eqref{eq:lerbvkherbvjhbdvcjhvsdcj} and \eqref{e:InLambda}, from \eqref{eq:alvrfleejhcbdjhsac} we get 
\begin{multline*}
    \mathcal L \psi_{+}(x)
        \leq n^{y_1^*}k^{y_2^*}\Bigg[\frac{1}{2} \kappa_-\bigg(\psi_{n,+}(k+(y_2^*)'-y_2^*) - \psi_{n,+}(k)\bigg)\\
        +|\cal{R}|\kappa_+\bigg(k^{-1}+k^{\Kz-y_2^*}n^{-1}\bigg)\big(\psi_{+}(x+(0,\Kz))-\psi_{+}(x)\big)\Bigg],
\end{multline*}    
  because for each reaction $y\to y'$ we have
  \begin{equation*}
   \psi_{+}(x+y'-y)-\psi_{+}(x)\leq \psi_{+}(x+(0,\Kz))-\psi_{+}(x)
  \end{equation*}
 and one of the following statements holds. In WE every reaction with $y_1\ne y_1'$ falls in the last case.
  \begin{itemize}
   \item $y=y^*$ and therefore $y_2'\leq (y_2^*)'$, so
   \begin{equation*}
       \psi_{n,+}(k+y_2'-y_2^*) - \psi_{n,+}(k)\leq \psi_{n,+}(k+(y_2^*)'-y_2^*) - \psi_{n,+}(k)<0.
   \end{equation*}
   \item $y_1=y_1^*$ and $y_2\leq y_2^*-1$, hence $\lambda_y((n,k))\leq n^{y_1^*}k^{y_2^*}k^{-1}$.
   \item $y_1<y_1^*$, hence $\lambda_y((n,k))\leq n^{y_1^*}k^{y_2}n^{-1}\leq n^{y_1^*}k^{y_2^*}k^{\Kz-y_2^*}n^{-1}\leq n^{y_1^*}k^{y_2^*}k^{\Kz}n^{-1}$.
  \end{itemize}

    The same three cases and \eref{e:InLambda} give $\mathbb P_x(Y(1)_2\ge k)\le C(k^{-1}+H^{\Kz}/n)$. For $K,n$ sufficiently large, the bounded jump sizes then give $\E_x[Y(1)_2-k]\le-c_0<0$. Summing up to $\tau_W\wedge l$ gives $c_0\E_x(\tau_W\wedge l)\le k$, hence $\E_x\tau_W\le CH$.

    We show \eref{e:poissonub} in two steps. First, we deal with the states $x=(n, k)$ on the top of the strip, that can jump above $H$, i.e. such that $k\in [K+\gamma\log n-\Kz, K+\gamma\log n]$. Here, by using
    \begin{equation*}
        \psi_{+}(x+(0,\Kz))-\psi_{+}(x)\leq 1,
    \end{equation*}
    we get
\begin{equs}
    \mathcal L \psi_{+}(x)
        &\leq n^{y_1^*}k^{y_2^*}\pq{\frac{1}{2} \kappa_-\bigg(\psi_{n,+}(k+(y_2^*)'-y_2^*) - \psi_{n,+}(k)\bigg)+|\cal{R}|\kappa_+\bigg(k^{-1}+k^{\Kz}n^{-1}\bigg)}\\
        &\leq n^{y_1^*}k^{y_2^*}\pq{\frac{1}{2}\kappa_-C(\log n)^{-1}\bigg(\exp[-\beta/\gamma]-1\bigg)+|\cal{R}|\kappa_+\bigg(\gamma^{-1}(\log n)^{-1}+(K+\gamma \log n)^{\Kz}n^{-1}\bigg)}
\end{equs}    
  Hence, taking $C$ large enough, inequality~\eref{e:poissonub} holds for $n$ large enough. 

  When $x=(n,k)$ is not on the top of the strip, i.e. when $K\leq k\leq K+\gamma\log n-\Kz$, none of the reactions can jump above $H$. Hence we have
  \begin{equation*}
  \psi_{+}(x+(0,\Kz))-\psi_{+}(x)=\psi_{n,+}(k+\Kz) - \psi_{n,+}(k)\leq\psi_{n,+}(k)\exp[\beta\Kz/\gamma]
  \end{equation*}
  and therefore
\begin{multline*}
    \mathcal L \psi_{+}(x)\leq n^{y_1^*}k^{y_2^*}\psi_{n,+}(k)\Bigg[\frac{1}{2}\kappa_-\bigg(\exp[-\beta/\gamma]-1\bigg)\\
    +|\cal{R}|\kappa_+\bigg(K^{-1}+(K+\gamma \log n)^{\Kz}n^{-1}\bigg)\exp[\beta \Kz/\gamma]\Bigg]
\end{multline*}    
and inequality~\eref{e:poissonub} holds for $n$ large enough, provided that
\begin{equation*}
 K>\frac{2 |\cal{R}|\kappa_+ \exp[\beta \Kz/\gamma]}{\kappa_-(1-\exp[-\beta/\gamma])}.
\end{equation*}
Since $\beta=\gamma$, the bound $\psi_{n,+}(k)$ is at most $CH^{-1}e^{k-H}$. A bounded change in $H$ changes only the constant in the estimate at the top of the strip.

Thus the probability of exit above $H$ is $O(n^{-\beta})$ for $K\le k\le K+\Kz$ and $O((\log n)^{-1})$ for all $k\le H$. For $W^c=\bigcup_{\ell\in L}\XX_{(j,K)}(\ell)$, an exit into $\XX_{(j,K+1)}$ is an exit above $H_n$.
\end{proof}

\begin{proof}[Proof of \lref{l:commonexit}]
    Throughout the proof, we use the following notations. 
    We take a reaction $\tilde{y}^*\to (\tilde{y}^*)'\in\arg\max_{y\to y'\in\widetilde{\cal R}}(y_1,y_2)$, with lexicographic order. Note that $\tilde \beta= y_1^*-\tilde y_1^*$. Put $A=A_n$, $A^0=A\setminus\XX_{(0,K)}(\ell)$ and $c^*=1/2$. The constant $c\ge1$ may increase from line to line. 
    We divide the strip $\XX_{(0,K)}(\ell)$ into three sets as follow: 
    \begin{equ}\label{e:DefSetMUD}
        \begin{cases}
            M_n:= \{(n,k)\in \XX_{(0, K)}(\ell): \, K+2\Kz< k<K+\gamma\log n-2\Kz\}\\
            D_n:= \{(n,k)\in \XX_{(0, K)}(\ell): \, k\leq K+2\Kz\}\\
            U_n:= \{(n,k)\in \XX_{(0, K)}(\ell): \, k\geq K+\gamma\log n-2\Kz\}\\
        \end{cases} \, .
    \end{equ}
    
    {\em Assumption (1):} Starting with the upper bound, we consider the function $\psi_{n,+}(k):= C(\log n)^{C'}n^{-\tilde\beta}k$, and set $\psi_+$ as 
    \begin{equ}
        \psi_+(x=(j,k)_n)=
        \begin{cases}
            1 \quad &\text{if} \quad x\in A^0\\
            0 \quad &\text{if}\quad  x\notin A\cup \XX_{(0, K)}(\ell)\cup \XX_{(0, K+1)}\\
            \psi_{n,+}(k) \quad &\text{if}\quad x\in \XX_{(0, K)}(\ell)\cup \XX_{(0, K+1)}
        \end{cases}\, . 
    \end{equ}
    We now show that~\eqref{e:poissonub} holds for this function. The second inequality is straightforward since $A\cap \XX_{(0, K+1)}=\emptyset$. Let $(n, k)\in \XX_{(0, K)}(\ell)$. Using Relation~\eqref{e:InLambda}, and similar arguments in the previous proof, we get 
    \begin{equs}
        \cal{L}\psi_{+}(n,k)\leq& -\kappa_-c^*n^{y_1^*}k^{y_2^*}C(\log n)^{C'}n^{\tilde{y}_1^*-y_1^*}(y_2^*-(y_2^*)')\\
            &\quad +|\cal{R}|\kappa_+cn^{y_1^*}k^{y_2^*-1}\Kz C(\log n)^{C'}n^{\tilde{y}_1^*-y_1^*}\\
            &\quad +|\widetilde{\cal{R}}|\kappa_+cn^{\tilde y_1^*}k^{\tilde y_2^*}\\
            &\quad +|\cal{R}|\kappa_+cn^{y_1^*-1}k^{\Kz}\Kz C(\log n)^{C'}n^{\tilde{y}_1^*-y_1^*}\\
            \leq & -n^{\tilde y_1^*}k^{y_2^*}\pq{C(\log n)^{C'}(\kappa_-c^*-\kappa_+ck^{-1}-\kappa_+ck^{\Kz}n^{-1})-\kappa_+ck^{\Kz}}
    \end{equs}
    where the second line concerns the reactions in $\cal{R}_0$ that do not start from $y^*$, whose positive contributions are bounded by $\Kz C(\log n)^{C'}n^{-\tilde\beta}$, also on the top cap. The third line corresponds to the exits to $A^0$ through a reaction of $\widetilde{\cal{R}}$, and finally the last term corresponds to the other reactions, that by definition of $\widetilde{\cal{R}}$, can't get out of the strip through $A^0$. Taking $C'$ large enough allows to compensate for the last term in the last line, since $k\leq K +\gamma\log n $, and we can conclude to~\eqref{e:poissonub} when $n$ and $K$ are large enough. 

    For the lower bound, we first show a uniform lower bound on $Q_{M_n^c}(x, A^0)$, where $M_n$ is defined in Relation~\eqref{e:DefSetMUD}. Put $\zeta_*=(\tilde y^*)'-\tilde y^*$. By the definition of $\widetilde{\cal R}$, some $\zeta'=y'-y$, $y\to y'\in\widetilde{\cal R}'$, satisfies
    \[
      \begin{gathered}
      \zeta_{*,1}=\zeta'_1,\qquad \zeta_{*,2}-\zeta'_2\in d\Z,
      \qquad |\zeta_{*,2}-\zeta'_2|\le2\Kz,\\
      x\in M_n\Longrightarrow x+\zeta_*\in A^0.
      \end{gathered}
    \]
    We set 
    \begin{equ}
        \psi_{-}(x=(j,k)_n)=
        \begin{cases}
            \mathbbm 1_{A^0}(x) \quad &\text{if} \quad x\in M_n^c\\
            \psi_{n,-}(k) \quad &\text{if} \quad x\in M_n. 
        \end{cases}
    \end{equ}   
    with $\psi_{n,-}(k):= C^{-1}(\log n)^{-C'}n^{-\tilde\beta}$ and $C, C'>0$ specified later. We now show that~\eqref{e:poissonlb} holds. For $(n,k)\in M_n$ and $n$ large enough, $\psi_{n,-}(k)\in [0,1]$. Since $\psi_-$ is constant on $M_n$, only the reactions getting out of $M_n$ change the generator. At the state $(n,k)$, 
    \begin{equs}
        \cal{L}\psi_-(n,k)
        \geq &  -|\cal{R}_0|\kappa_+n^{y_1^*}k^{y_2^*}C^{-1}(\log n)^{-C'}n^{\tilde{y}_1^*-y_1^*}\\
        &-|\cal{R}|\kappa_+n^{y_1^*-1}k^{\Kz}C^{-1}(\log n)^{-C'}n^{\tilde{y}_1^*-y_1^*}\\
        &+\kappa_-c^*n^{\tilde y_1^*}k^{\tilde y_2^*}
    \end{equs}
    where the last term is due to the reaction that gets out to $A^0$. 
    Taking $C'>\Kz$, we can conclude to~\eqref{e:poissonlb} when $n$ is large enough. 
    Therefore, for $x\in M_n$, 
    \[
        Q_{\XX_{(0,K)}(\ell)^c}(x, A^0)\geq Q_{M_n^c}(x, A^0)\geq \psi_-(x)=C^{-1}(\log n)^{-C'}n^{-\tilde\beta}\, . 
    \]

    For $x\in U_n$, a bounded number of $y^*\to(y^*)'$ jumps reaches a point $x^*(x)\in M_n$. For $x\in D_n$, choose a fixed dominant path to $M_n$, staying above $K$. Writing $\mathbb P_{C_n}$ for the transition matrix killed outside $C_n=\XX_{(0,K)}(\ell)$, we have
    \[
      \begin{gathered}
      Q_W(x,A^0)\geq\mathbb P_{C_n}^{l(x)}(x,x^*(x))Q_W(x^*(x),A^0),\\
      \inf_{x\in U_n\cup D_n}\mathbb P_{C_n}^{l(x)}(x,x^*(x))\ge c_0>0.
      \end{gathered}
    \]
    The last bound follows from~\eref{e:InLambda} on $U_n$, and from the fixed height and length bounds of the dominant paths on $D_n$.

    {\em Assumption (2):} The proof of the upper bound is left to the reader. If $\widetilde{\cal{R}}\cap \cal{R}_0\neq \emptyset$, it is straightforward since $\tilde{\beta}=0$, and if $\widetilde{\cal{R}}\cap \cal{R}_0= \emptyset$, it can be proven analogously to the upper bound in {\em Assumption (1)}. 
    We start by showing that the probability to exit $M_n\cup U_n$ through $D_n$ is bounded away from $0$, i.e. show a uniform lower bound on $Q_{W\cup D_n}(x, D_n)$ for $x\in M_n\cup U_n$. 
    For some $\alpha\in (0,1)$ specified later, we define the functions 
    \begin{equ}
        \psi_{n,-}(k)=1-\alpha -\frac{(1-2\alpha)k}{K+\gamma\log n -K_0} \quad \text{and}\quad \psi_-(x=(j,k)_n)=
        \begin{cases}
            \psi_{n,-}(k)\quad &\text{if}\quad x\in \XX_{(0,K)}(\ell)\\
            0\quad &\text{if} \quad x\notin \XX_{(0,K)}(\ell). 
        \end{cases}
    \end{equ}
    Note that $\psi_-$ is non negative. We are going to show Relation~\eqref{e:poissonlb} where $W^c$ is here $M_n\cup U_n$ and $A$ is $D_n$. The second inequality is straightforward. 
    For $(n,k)\in M_n$, we have 
    \begin{equs}
        \cal{L}\psi_-(n,k)&\geq \kappa_-c^*n^{y_1^*}k^{y_2^*}\frac{1-2\alpha}{K+\gamma\log n-\Kz}(y_2^*-(y_2^*)')\\
        &\, -|\cal{R}_0|\kappa_+cn^{y_1^*}k^{y_2^*-1}\frac{1-2\alpha}{K+\gamma\log n-\Kz}\Kz\\
        &\, -|\cal{R}|\kappa_+cn^{y_1^*-1}k^{\Kz}\psi_{n,-}(k)\\
        &\geq n^{y_1^*}k^{y_2^*}\left[\frac{1-2\alpha}{K+\gamma\log n-\Kz}\left(\kappa_-c^*-\frac{\Kz|\cal{R}_0|\kappa_+c}{k}\right)-\frac{|\cal{R}|k^{\Kz-y_2^*}\kappa_+c}{n}\right]
    \end{equs}
    and the last term is positive for $n\gg 1$ and $K$ large enough. The second line contains the contribution of the upward reactions in $\cal{R}_0$ that do not start at $y^*$, which stay in $\XX_{(0, K)}(\ell)$ by definition of the set, and because we started from $M_n$ (and therefore can't exit from above or below).

    For $(n,k)\in U_n$, $\psi_{n,-}(k)<2\alpha$, therefore 
    \begin{equs}
        \cal{L}\psi_-(n,k)&\geq \kappa_-c^*n^{y_1^*}k^{y_2^*}\frac{1-2\alpha}{K+\gamma\log n-\Kz}(y_2^*-(y_2^*)')\\
        &\, -|\cal{R}_0|\kappa_+cn^{y_1^*}k^{y_2^*-1}2\alpha
        -|\cal{R}|\kappa_+cn^{y_1^*-1}k^{\Kz}2\alpha\\
        &\geq \frac{n^{y_1^*}k^{y_2^*}}{K+\gamma\log n-\Kz}\left[\kappa_-c^*(1-2\alpha)-4|\cal{R}_0|\kappa_+c\alpha-\frac{|\cal{R}|k^{\Kz-y_2^*}\kappa_+c(K+\gamma\log n-\Kz)}{n}\right]
    \end{equs}
    where in the first inequality, we used the fact that $y_2^*>(y_2^*)'$, and therefore $(n,k)+(y^*)'-y^*\in \XX_{(0,K)}(\ell)$.
    Taking $K\ge2\Kz$, $n\gg 1$ and $\alpha$ small enough, we can conclude to 
    \begin{equ}
        \inf_{x\in M_n\cup U_n}Q_{W\cup D_n}(x, D_n)\geq \inf_{x\in M_n\cup U_n}\psi_-(x)\geq \alpha/2. 
    \end{equ}

    We now focus on $D_n$. Let $x^*\in \XX_{(0,K)}(\ell)$ be a state for which a reaction in $\arg \max_{y\to y'\in \widetilde{\cal{R}}} y_1$ leads to $A$. $x^*$ is in $D_n$ by definition of $\Kz$. For $x_0\in D_n$, we have 
    \[
        Q_W(x_0, A)\geq \inf_{x'\in D_n} \mathbb{P}_{C_n}^{l(x')}(x', x^*)\mathbb{P}(x^*, A)
    \]
    where $l(x')$ is the length of a fixed dominant path staying above $K$. 
    The first term is bounded away from zero by the fixed height and length bounds of these paths, and we chose $x^*$ such that  $\mathbb P(x^*, A)\geq C^{-1} n^{-\tilde{\beta}}$. This concludes the proof for $x\in D_n$, and for $x\in M_n\cup U_n$, note that 
    \[
        Q_W(x,A)\geq Q_{W\cup D_n}(x, D_n)\min_{x\in D_n}Q_W(x,A), 
    \]
and combining the previous results, we can conclude. 
\end{proof}

\section{Proof of \lref{l:exitestimate1} }

\begin{proof}[Proof of \lref{l:exitestimate1}] The first part of the lemma corresponds exactly to \cite[Lemma 3.3, Chapter 6]{fw98}. Below, we mirror this proof to establish the second part of the lemma, with the original proof being recovered by setting $\XX_{\bar m} = \emptyset$. Throughout the proof, we write $\XX_j=\OO_j$, $Q_A=Q_{\OO_A}$, and define $\bar L := L \setminus \{\bar m\}$.

  This proof is carried out by induction on $r$. The case $r = 1$ trivially holds with the sum in the numerator being carried out on a singleton and, for any $x \in \XX_m$,
  \begin{equ}\label{e:denominator}
    a^{-2}  \leq  a^{-1}\sum_{m' \in \bar L\setminus \{m\}} Q_W(x, \XX_{m'}) \leq \sum_{g \in G(\bar L \setminus \{m\})} \pi(g)   \leq  a \sum_{m' \in \bar L\setminus \{m\}} Q_W(x, \XX_{m'}) \leq a \leq  a^2\,.
  \end{equ}
  where in the first inequality we have used \eref{e:exitestimate3},
and $a \geq 2$, giving
\begin{equ}
  \sum_{m' \in \bar L\setminus \{m\}} Q_W(x, \XX_{m'}) = 1- Q_W(x, \XX_{\bar m}) \geq \tfrac12 \geq a^{-1}\,.
  \end{equ}

We now proceed to consider the induction step: we assume that \eref{e:exitestimate2} holds for all $W$ such that $|L\setminus W|=r$, for all $m \in L \setminus W$ and for all $m' \in W\setminus \{\bar m\}$. Consider now a set $W$ such that $|L\setminus W| = r+1$, fix $m \in L\setminus W, m' \in W \setminus \{\bar m\}$ for which we want to establish \eref{e:exitestimate2} and define
\begin{equ}
  \mathbf F := \bigcup_{m'' \in L\setminus W, m'' \neq m} \mathbf Y_{m''}
\end{equ}
where
\begin{equ}\label{e:Ym}
\mathbf Y_m := \{x \in \XX_m~:~ \exists \,m'' \in L\setminus W,\,  y \in \XX_{m''} \text{ s.t. } \mathbb P_n(y,x) > 0\}\,.
\end{equ}
Then we expand
\begin{equs}\label{e:exp1}
  Q_W(x, \XX_{m'}) & = Q_{L\setminus \{m\}}(x, \XX_{m'}) + \int_{\mathbf F} Q_{L\setminus \{m\}}(x, dy) Q_{W \cup \{m\}}(y, \XX_{m'})\\
  &  \quad + \int_{\mathbf F} Q_{L\setminus \{m\}}(x, dy) \int_{\mathbf Y_m} Q_{W \cup \{m\}}(y, d y') Q_{L\setminus \{m\}}(y', \XX_{m'})\\
    & \quad  + \int_{\mathbf F} Q_{L\setminus \{m\}}(x, dy) \int_{\mathbf Y_m} Q_{W \cup \{m\}}(y, d y')  \int_{\mathbf F} Q_{L\setminus \{m\}}(y', dy'') Q_{W \cup \{m\}}(y'', \XX_{m'})\\
  & \quad + \int_{\mathbf F} Q_{L\setminus \{m\}}(x, dy) \int_{\mathbf Y_m} Q_{W \cup \{m\}}(y, d y')\int_{\mathbf F} Q_{L\setminus \{m\}}(y', dy'') \\
  &\qquad \qquad \qquad \qquad \qquad \qquad\qquad  \cdot \int_{\mathbf Y_m} Q_{W \cup \{m\}}(y'', d y''') Q_{L\setminus \{m\}}(y''', \XX_{m'}) + \dots
\end{equs}

We now proceed to prove the upper bound and the lower bound separately, starting from the former.
{\bf Upper bound:}
Now, defining
\begin{equs}
  \overline A_{m,m''} & := \sup_{y \in \XX_m}Q_{L\setminus \{m\}}(y, \XX_{m''})& m'' \neq m\qquad\quad\, \\
    \overline B_{m,m''} & := \sup_{y \in \XX_m}\int_{\mathbf F}Q_{L\setminus \{m\}}(y, d y') Q_{W \cup \{m\}}(y', \XX_{m''}) \qquad & m'' \in  W \cup \{m\}
\end{equs}
we can bound \eref{e:exp1} as
\begin{equs}
Q_W(x, \XX_{m'}) \leq \overline A_{m,m'} + \overline B_{m,m'} + \overline B_{m,m}\overline A_{m,m'} + \overline B_{m,m} \overline B_{m,m'} + \dots
\end{equs}
Further, choose a simple path $m=m_0,m_1,\ldots,m_\ell\in W\setminus\{\bar m\}$ of positive $q$-weights, with $m_i\notin W$ for $i<\ell$. Such paths exist by positivity of the forest denominator and give almost-sure absorption in $W$. Observing that
\begin{equ}
  0\leq\overline B_{m,m}\leq1-\delta_m<1,
  \qquad \delta_m:=a^{-\ell}\prod_{i=0}^{\ell-1}q_{m_i,m_{i+1}}>0,
\end{equ}
we can write
\begin{equs}
 Q_W(x, \XX_{m'}) \leq \frac{\overline A_{m,m'} + \overline B_{m,m'}}{1 - \overline B_{m,m}}
\end{equs}

First, we have proved for $r=1$ that
\begin{equ}
  \overline A_{m,m'} \leq a q_{m,m'} \leq a^4\frac{q_{m,m'}}{\sum_{m'' \in \bar L \setminus\{m\} } q_{m,m''}}\,.
\end{equ}
Then, recalling that all exit graphs (\eg $G(W \cup \{m\}), G_{m'',m'}(W \cup \{m\})$) are intended as graphs on the reduced vertex set $\bar L$, we have
\begin{equs}
  \overline B_{m,m'} & \leq \sum_{m'' \in L \setminus W, m'' \neq m} a^4\frac{q_{m,m''}}{\sum_{m'' \in \bar L \setminus\{m\} } q_{m,m''}} a^{4^r} \frac{\sum_{g \in G_{m'',m'}(W \cup \{m\})}\pi(g)  }{\sum_{g \in G(W \cup \{m\})}\pi(g)  } \\
  &= a^{4^r + 4} \frac{D_{m,m'}}{\sum_{m'' \in \bar L \setminus\{m\} } q_{m,m''} E_m}
\end{equs}
where $D_{m,m'}$ is the sum of products $\pi(g)  $ over those graphs in $G_{m,m'}(W)$ (again on the reduced vertex set $\bar L$) where the arrow beginning at $m$ does not lead immediately to $m'$, while $E_m$ is the sum of the same products over all $(W \cup \{m\})$-graphs on $\bar L$. Combining the bounds above we obtain
\begin{equ}\label{e:exp2}
  \overline A_{m,m'} + \overline B_{m,m'} \leq a^{4^r + 4} \frac{q_{m,m'} E_m + D_{m,m'}}{\sum_{m'' \in \bar L \setminus\{m\} } q_{m,m''} E_m}
\end{equ}
and notice that the numerator corresponds to the sum of $\pi(g )$ over all graphs on $\bar L$ belonging to $G_{m,m'}(W)$.

We conclude the proof by estimating the denominator $1-\overline B_{m,m}$. Since the denominator is positive by assumption, we have
\begin{equ}
  Q_{W \cup \{m\}}(y,\XX_m) = 1- Q_{W \cup \{m\}}(y,\bigcup_{m'' \in W}\XX_{m''})
\end{equ}
so that
\begin{equs}
  1-\overline B_{m,m}  & = \inf_{x \in \XX_{m}}\left[1-\int_{\mathbf F} Q_{L\setminus \{m\}}(x, dy) \left(1- Q_{W \cup \{m\}}\pc{y,\bigcup_{m'' \in W}\XX_{m''}}\right)\right]\\
  & = \inf_{x \in \XX_m}\left[Q_{L \setminus \{m\}}\pc{x,\bigcup_{m'' \in W}\XX_{m''}} + \int_{\mathbf F} Q_{L\setminus \{m\}}(x, dy) Q_{W \cup \{m\}}\pc{y,\bigcup_{m'' \in W}\XX_{m''}}\right]
\end{equs}
The terms on the righthand side can respectively be bounded from below:
\begin{equs}
  \inf_{x \in \XX_m} Q_{L \setminus \{m\}}\pc{x,\bigcup_{m'' \in W}\XX_{m''}}\geq a^{-1}\sum_{m''\in W\setminus\{\bar m\}}q_{m,m''}\geq a^{-3}\frac{\sum_{m''\in W\setminus \{\bar m\}}q_{m,m''}}{\sum_{m''\in \bar L \setminus \{m\}}q_{m,m''}}\\
\end{equs}
and
\begin{equs}
  \int_{\mathbf F} &Q_{L\setminus \{m\}}(x, dy) Q_{W \cup \{m\}}\pc{y,\bigcup_{m'' \in W}\XX_{m''}}\\
  &\geq \sum_{m''' \in L\setminus W, m''' \neq m}\sum_{m''\in W\setminus\{\bar m\}} a^{-3}\frac{q_{m,m'''}}{\sum_{j \in \bar L\setminus \{m\}} q_{m,j}} a^{-4^r} \frac{\sum_{g \in G_{m''',m''}(W \cup \{m\})}\pi(g)  }{\sum_{g \in G(W \cup \{m\})}\pi(g)  } \\
  &= a^{-3-4^r} \frac{H_{m}}{\sum_{m'' \in \bar L \setminus \{m\}} q_{m,m''} E_m}
\end{equs}
where $H_m$ is the sum of the products $\pi(g)  $ over those graphs in $G(W)$ in which the arrow beginning at $m$ leads to a point in $L\setminus W$, so that
\begin{equ}
  1-\overline B_{m,m} \geq a^{-3-4^r} \frac{H_{m} + \sum_{m'' \in W\setminus \{\bar m\}} q_{m,m''} E_m}{\sum_{m'' \in \bar L \setminus \{m\}} q_{m,m''} E_m}
\end{equ}
Noting that the denominator corresponds to the one of \eref{e:exp2} and that the numerator we have the sum of $\pi(g)  $ over all graphs $g\in G(W)$ we finally obtain
\begin{equ}
  \frac{\overline A_{m,m'} + \overline B_{m,m'}}{1 - \overline B_{m,m}} \leq a^{8+2\cdot 4^r}\frac{q_{m,m'} E_m + D_{m,m'}}{H_{m} + \sum_{m'' \in W\setminus \{\bar m\}} q_{m,m''} E_m}
\end{equ}
from which we recover the desired upper bound.\\[10pt]

{\bf Lower bound:} The proof of the lower bound proceeds analogously to the upper bound: defining
\begin{equs}
  \underline A_{m,m''} & := \inf_{y \in \XX_m}Q_{L\setminus \{m\}}(y, \XX_{m''})& m'' \neq m\qquad\quad\, \\
    \underline B_{m,m''} & := \inf_{y \in \XX_m}\int_{\mathbf F}Q_{L\setminus \{m\}}(y, d y') Q_{W \cup \{m\}}(y', \XX_{m''}) \qquad & m'' \in  W \cup \{m\}
\end{equs}
from \eref{e:exp1} we write
\begin{equs}\label{e:lb1}
  Q_W(x, \XX_{m'})  \geq \pc{\sum_{k=0}^\infty Q_{mm}^k(x,\XX_m) }(\underline A_{m,m'} + \underline B_{m,m'}) \\
\end{equs}
where
\begin{equs}
  Q_{mm}^k(x,\XX_m) & = \int_{\mathbf F} Q_{L \setminus \{m\}}(x,dy) \int_{\XX_m}Q_{W \cup \{m\}}(y,dy')\\
  & \qquad \qquad \qquad\dots \int_{\mathbf F} Q_{L \setminus \{m\}}(y^{(2k-3)},dy^{(2k-2)}) \int_{\XX_m}Q_{W \cup \{m\}}(y^{(2k-2)},dy^{(2k-1)})
\end{equs}
To estimate the first term on the right hand side of \eref{e:lb1} we write for $x \in \XX_m$, as in \eref{e:exp1}
\begin{equs}
  1 - Q_W(x, \XX_{\bar m}) &= \sum_{m' \in W \setminus \{\bar m\}} Q_W(x, \XX_{m'})\\
  & \leq \sum_{m' \in W \setminus \{\bar m\}}  \pc{\sum_{k=0}^\infty Q_{mm}^k(x,\XX_m) (\overline A_{m,m'} + \overline B_{m,m' })}\\
  & = \pc{\sum_{m' \in W \setminus \{\bar m\}}(\overline A_{m,m'} + \overline B_{m,m' })}\pc{\sum_{k=0}^\infty Q_{mm}^k(x,\XX_m) }
\end{equs}

Therefore, we have
\begin{equ}
  \pc{\sum_{k=0}^\infty Q_{mm}^k(x,\XX_m) }^{-1} \leq \frac{\sum_{m' \in W \setminus \{\bar m\}}(\overline A_{m,m'} + \overline B_{m,m' })}{1 - Q_W(x, \XX_{\bar m})} \leq a \sum_{m' \in W \setminus \{\bar m\}}(\overline A_{m,m'} + \overline B_{m,m' })
\end{equ}
where in the last inequality we have used \eref{e:exitestimate3} and $a \geq 2$.
Combining the above with \eref{e:lb1} we obtain
\begin{equs}\label{e:lb4}
  Q_W(x, \XX_{m'}) \geq a^{-1}\frac{\underline A_{m,m'} + \underline B_{m,m'}}{ \sum_{m'' \in W \setminus \{\bar m\}}(\overline A_{m,m''} + \overline B_{m,m''})}
\end{equs}
and proceed to bound the numerator and the denominator. Replicating the computations in the upper bound case we have, analogously to \eref{e:exp2}
\begin{equ}\label{e:lb2}
  \underline A_{m,m'} + \underline B_{m,m'} \geq a^{-4^r - 3} \frac{q_{m,m'} E_m + D_{m,m'}}{\sum_{m'' \in \bar L \setminus\{m\} } q_{m,m''} E_m}
\end{equ}
where $E_m, D_{m,m'}$ were defined in the upper bound part of the proof. For the denominator we have
\begin{equs}
 \sum_{m''\in W\setminus\{\bar m\}}\overline A_{m,m''} \leq a \sum_{m''\in W\setminus\{\bar m\}}q_{m,m''}\leq a^{4}\frac{\sum_{m''\in W\setminus \{\bar m\}}q_{m,m''}}{\sum_{m''\in \bar L \setminus \{m\}}q_{m,m''}}\\
\end{equs}
and, recalling the definition of $H_m$, by the induction assumption
\begin{equs}
  \sum_{m''\in W\setminus\{\bar m\}}\overline B_{m,m''} &\leq \sum_{m''' \in L\setminus W, m''' \neq m}\sum_{m''\in W\setminus\{\bar m\}} a^{4}\frac{q_{m,m'''}}{\sum_{j \in \bar L\setminus \{m\}} q_{m,j}} a^{4^r} \frac{\sum_{g \in G_{m''',m''}(W \cup \{m\})}\pi(g)  }{\sum_{g \in G(W \cup \{m\})}\pi(g)  } \\
  &= a^{4+4^r} \frac{H_{m}}{\sum_{m'' \in \bar L \setminus \{m\}} q_{m,m''} E_m}
\end{equs}
so that
\begin{equ}\label{e:lb3}
   \sum_{m'' \in W \setminus \{\bar m\}}(\overline A_{m,m''} + \overline B_{m,m'' }) \leq a^{4+4^r} \frac{H_{m}+ \sum_{m'' \in W\setminus \{\bar m\}} q_{m,m''} E_m}{\sum_{m'' \in \bar L \setminus \{m\}} q_{m,m''} E_m}
\end{equ}
Finally, combining\eref{e:lb4} with \eref{e:lb2} and \eref{e:lb3} yields
\begin{equ}
  Q_W(x, \XX_{m'}) \geq a^{-8-2\cdot 4^r}\frac{q_{m,m'} E_m + D_{m,m'}}{H_{m} + \sum_{m'' \in W\setminus \{\bar m\}} q_{m,m''} E_m}
\end{equ}
which proves the claim.
\end{proof}

\end{document}